\documentclass[11pt,twoside]{article}
\usepackage{amsmath,amssymb}
\usepackage{graphicx}
\usepackage[footnotesize]{subfigure}
\usepackage[footnotesize]{caption2}
\usepackage{authblk}
  \usepackage{color}
\allowdisplaybreaks[4]

\renewcommand{\theequation}{\thesection.\arabic{equation}}

\newtheorem{lemma}{{\bf Lemma}}[section]
\newtheorem{example}{{\bf Example}}[section]

\newtheorem{theorem}{{\bf Theorem}}[section]
\newtheorem{remark}{{\bf Remark}}[section]

\makeatletter
\evensidemargin
\oddsidemargin
\makeatother

\title{
\Large\bf Dynamics of  the  Tyson-Hong-Thron-Novak circadian\\ oscillator model
}

\author{ {\sc Shuang Chen$^{a,b,}$\footnote{\rm{Email}: {schen@hust.edu.cn}}},
{\sc Jinqiao Duan$^{c,}$\footnote{\rm{Email}: {duan@iit.edu}}},
{\sc Ji Li$^{a,}$\footnote{\rm{Email}: {liji@hust.edu.cn}}}
\\
{\small $^{a}$ School of Mathematics and Statistics, Huazhong University of Sciences and Technology}\\
{\small Wuhan, Hubei 430074, P. R. China}\\
{\small $^{b}$ Center for Mathematical Sciences, Huazhong University of Sciences and Technology}\\
{\small Wuhan, Hubei 430074, P. R. China}\\
{\small $^{c}$ Department of Applied Mathematics, Illinois Institute of Technology}\\
{\small Chicago, IL 60616, USA}
}

\date{}

\begin{document}
\maketitle
\begin{abstract}

We study the dynamics of a circadian oscillator model
which was proposed by Tyson, Hong, Thron and Novak.
This model describes a molecular mechanism for the circadian rhythm in Drosophila.
After giving a detailed study of the equilibria,
we investigate  the effects of the rates of mRNA degradation and synthesis.
When the  rate of mRNA degradation is   high enough,
we prove that there are no periodic orbits in this model.
When the  rate of mRNA degradation is sufficiently  low,
this model is transformed into a slow-fast system.
Then based on the Geometric Singular Perturbation Theory,
we prove the existence of canard explosion, relaxation oscillations,
homoclinic orbits, heteroclinic orbits and saddle-node bifurcations
as  the rates of mRNA degradation and synthesis change.
Finally, we give the biological interpretation of the obtained results
and  point out that this model  can be transformed into a Li\'enard-like equation,
which could be helpful to investigate the dynamics of the general case.


\vskip 0.2cm
{\bf AMS(2020) Subject Classification}: 34C26; 34C37; 34E17.
\end{abstract}
\baselineskip 15pt
\parskip 10pt

\thispagestyle{empty}
\setcounter{page}{1}


\section{Introduction}\label{sec-intr}
\setcounter{equation}{0}
\setcounter{lemma}{0}
\setcounter{theorem}{0}
\setcounter{remark}{0}

Circadian rhythms of physiology with a period about 24
hours have been found in many organisms,
for example, in fruit flies, plants and vertebrate animals.
These circadian clocks allow us to adapt to the alternation of day and night.
In order to grasp the mechanisms for the generation of circadian rhythms,
numerous theoretical models ranging from generic autonomous oscillators to molecular-based models
have been proposed in the past tens of years.
See, for example, \cite{Dunlap-99,Forger-17,Gonze-11,Keener-Sneyd-98,Leloup-Goldbeter-98} and the references therein.

Based on the dimerization and proteolysis of PER and TIM proteins in Drosophila,
Tyson, Hong, Thron and Novak \cite{Tyson-etal-99}  in 1999 set up a three-dimensional circadian oscillator model
\begin{eqnarray}
\label{3D-model-1}
\begin{aligned}
\frac{d M}{d t} &= \frac{\nu_m}{1+(P_2/P_c)^{2}}-k_m M,
\\
\frac{d P_1}{d t} &= \nu_p M-\frac{k_1P_1}{J_p+P_1+rP_2}-k_3P_1-2k_a P_1^{2}+2k_d P_2,
\\
\frac{d P_2}{d t} &= k_a P_1^{2}-k_d P_2-\frac{k_2P_2}{J_p+P_1+rP_2}-k_3P_2,
\end{aligned}
\end{eqnarray}
where the system states $M$, $P_{1}$ and $P_{2}$ denote the concentration of mRNA, monomer and dimer, respectively.
The biological descriptions of the model parameters are shown in Table \ref{tab-bm} (see also in \cite[Table 1]{Tyson-etal-99}).
\begin{table}[!htbp]
\centering
\begin{tabular}{cl}
  \hline
{\rm Parameter} & {\rm Biological description}  \\
\hline
 $v_{m}$ & {\rm the maximum rate of mRNA synthesis}\\
  $k_{m}$ & {\rm the first-order rate of mRNA degradation}\\
$P_{c}$ &  {\rm the value of dimer at the half-maximum transcription rate}\\
$v_{p}$ &  {\rm the rate for translation of mRNA into the monomer}\\
$k_{1}$ &  {\rm the maximum rate for monomer  phosphorylation}\\
$k_{2}$ &  {\rm the maximum rate for  dimer phosphorylation}\\
$k_{3}$ &  {\rm the first-order degradation rate  of the monomer and dimer}\\
$J_{P}$ &  {\rm the Michaelis constant for protein kinase DBT}\\
$k_{a}$ &  {\rm the rate of dimerization}\\
$k_{d}$ &  {\rm the rate of dissociation of the dimer}\\
$r$     &  {\rm the ratio of enzyme-substrate dissociation constants for the monomer and dimer}\\
\hline
\end{tabular}
\vskip 0.2cm
\centering\caption{The biological descriptions of the model parameters.}
\label{tab-bm}
\end{table}
Let the ratio $r=2$ and $k_{1}>k_{2}$.
Additionally, assume that the dimerization reactions $k_a$ and $k_d$ are  sufficiently large compared to other rate parameters,
Tyson, Hong, Thron and Novak \cite{Tyson-etal-99} applied the quasi-steady-state approximation
(see, for instance, \cite{Boieetal-16,Goeke-Walcher-Zerz-15})
to reduce the three-dimensional system (\ref{3D-model-1}) into  a simpler two-dimensional system
\begin{eqnarray}
\label{2D-model-1}
\begin{aligned}
\frac{d M}{d t} &= \frac{4\nu_m P_{c}^{2}}{4P_{c}^{2}+(P-h(P))^{2}}-k_m M,
\\
\frac{d P}{d t} &= \nu_p M-\frac{(k_1-k_2)h(P)+k_2P}{J_p+P}-k_3P,
\end{aligned}
\end{eqnarray}
where $P=P_1+2 P_2$ denotes the total amount of PER protein,
the constant $K=k_{a}/k_{d}$ and the function $h$ is given by
\begin{eqnarray*}
h(P)=\frac{\sqrt{1+8KP}-1}{4K}, \ \ \ \ \ P\geq 0.
\end{eqnarray*}
Here, system (\ref{2D-model-1}) is called the  two-dimensional Tyson-Hong-Thron-Novak circadian oscillator model
(the THTN model for short).

Although the THTN model  has the lower dimension
than that of the original system (\ref{3D-model-1}),
there are two obstacles in analyzing its dynamics,
that is, the THTN model possesses multiple parameters
and is  topologically equivalent to a high-order polynomial system.
In order to explore the properties of the THTN model,
Tyson et al. \cite{Tyson-etal-99} numerically studied the periods of limit cycles
by varying $(K,k_{1})$ and fixing other parameters,
and found that the THTN model  has a limit cycle with a period of about 24 hours
in a large parameters domain of $(K, k_{1})$.
Simon and Volford \cite{Simon-Volford-06} used the parametric representation method to
study the properties of equilibria and bifurcation curves
by varying $(v_{p}, k_{1}$) and fixing other parameters.
Goussis and Najm \cite{Goussis-Najm-06} numerically compared
the differences of periodic solutions in  the original system (\ref{3D-model-1})
and the THTN model.
Jiang et al. \cite{Jiang-etal-17} numerically studied the effects of several model parameters on the
the periods of circadian oscillations,
and pointed out that it is greatly reasonable to apply the THTN model to study the periodic behaviors in the original system (\ref{3D-model-1}).

In the actual experiment,
it is greatly important to investigate the effects of the model parameters
on the periodic behaviors in circadian oscillator models.
Our goal is to investigate the effects of
the rates of mRNA degradation and synthesis on the periodic behaviors in the THTN model.
In particular, we focus on the cases that
the  rate of mRNA degradation is much  high or  low,
that is, the rate $k_{m}$ is sufficiently large or small.
The analysis of the THTN model  with  general $k_{m}$ is a more complicated problem,
it will be studied in future work.
In the final section, we also point out that the THTN model is topologically equivalent to
a Li\'enard-like equation.
This structure is helpful to study the global dynamics of the THTN model with general $k_{m}$
and the effects of the model parameters on the periods of circadian oscillators.

When the rate of mRNA degradation  is  high enough,
this case is called the  high degradation rate case for simplicity.
We first obtain the existence of a bounded attractor by applying Gronwall's Inequality.
Then we further prove that there are no periodic orbits in the THTN model and
all orbits starting from the initial values in the domain
with biological meaning are attracted to locally stable foci or nodes,
except for the stable manifolds of saddles.
This indicates that  circadian oscillations  could disappear when
the rate of  mRNA degradation is  high.

When the rate of mRNA degradation  is  low enough,
this case is called the  low degradation rate  case. In this case,
the THTN model is topologically equivalent to a standard slow-fast system,
which is clearly separated into one slow variable and one fast variable.
By varying the rate $k_{m}$ of  mRNA degradation  and the ratio of the rate $v_{m}$ of mRNA synthesis
to the rate $k_{m}$ of mRNA degradation,
we further analyze the periodic phenomena in the  low degradation rate case.
The analysis for this case is based on the geometric singular perturbation theory.
For convenience, we introduce some basic notions on geometric singular perturbation theory in section \ref{sec-GSPT}.
Under the assumption that the critical manifold is $S$-shaped,
then two non-hyperbolic points such as canard points and jump points \cite{Krupa-Szmolyan-01SIMA} could appear.
Consequently, the desired circadian oscillators should appear
in the form of canard cycles and relaxation oscillations \cite{Krupa-Szmolyan-01JDE},
which are obtained by establishing the normal forms near the canard points
and applying the results obtained by \cite{Fenichel-79,Krupa-Szmolyan-01SIMA,Krupa-Szmolyan-01JDE}.
Besides these oscillations,
we also investigate the saddle-node bifurcations via the normal form near saddle-node points,
and prove the existence of homoclinic orbits and heteroclinic orbits by  the Fenichel Theorem \cite[Theorem 9.1]{Fenichel-79}
and the results in \cite{Krupa-Szmolyan-01SIMA}.

This paper is organized as follows.
In section \ref{sec-GSPT},
we introduce basic notions on geometric singular perturbation theory as preparations.
In section \ref{sec-singuar-point},
we provide a complete classification of the equilibria with no parameters fixed.
In sections \ref{sec-dynamics-1} and \ref{sec-dynamics-2},
we analyze the dynamics of the THTN model in the  high degradation rate case and the  low degradation rate case, respectively.
We also give some remarks on the further study in the final section.

\section{Geometric singular perturbation theory}
\label{sec-GSPT}
\setcounter{equation}{0}
\setcounter{lemma}{0}
\setcounter{theorem}{0}
\setcounter{remark}{0}

Multiple time scale systems frequently appear in many practical applications,
such as population dynamics, cellular physiology, fluid mechanics and so on
(see, for instance, \cite{Deng-03,DuLiLi-18,Bob-Liu-07,Hek-2010,Jones-95,Kuehn-15,Li-Zhu-13,Rubin-Terman,Wang-Zhang-19,Wiggins-94}).
These systems usually admits a clear separation in two time scales,
one slow time scale and one fast time scale,
which are also called the slow-fast systems.
Following the pioneering work \cite{Fenichel-79} of Fenichel in 1979,
geometric singular perturbation theory has been developed to be an efficient method to study multiple time scale dynamics.

Now we introduce some basic notions on geometric singular perturbation theory for planar slow-fast systems.
Consider a planar slow-fast system of the form
\begin{eqnarray}
\label{fast-1}
\begin{aligned}
\frac{d x}{d t} &= x'  = f(x,y,\mu,\varepsilon),
\\
\frac{d y}{d t}&= y'  = \varepsilon g(x,y,\mu,\varepsilon),
\end{aligned}
\end{eqnarray}
where $(x,y)\in \mathbb{R}^{2}$,
$\mu\in\mathbb{R}^{m}$ with $m\geq 1$,
a small parameter $\varepsilon$ with  $0<\varepsilon \ll1$,
and the functions $f$ and $g$ are $C^{k}$ with $k\geq 3$.
Letting $\tau=\varepsilon t$,
system (\ref{slow-1}) is rescaled to
\begin{eqnarray}
\label{slow-1}
\begin{aligned}
\varepsilon\frac{d x}{d \tau} &=  \varepsilon \dot x = f(x,y,\mu,\varepsilon),
\\
\frac{d y}{d \tau} &= \dot y  = g(x,y,\mu,\varepsilon).
\end{aligned}
\end{eqnarray}
In the limiting case $\varepsilon=0$,
system (\ref{fast-1}) becomes the layer equation
\begin{eqnarray}
\label{layer-1}
\begin{aligned}
x'  &= f(x,y,\mu,0),
\\
y'  &= 0,
\end{aligned}
\end{eqnarray}
and system (\ref{slow-1}) becomes the reduced equation
\begin{eqnarray}
\label{reduce-1}
\begin{aligned}
0&=  f(x,y,\mu,0),
\\
\dot y  & = g(x,y,\mu,0).
\end{aligned}
\end{eqnarray}
For the layer equation (\ref{layer-1}) with a fixed $\mu\in\mathbb{R}^{m}$,
its equilibria set $\mathcal{C}_{\mu,0}:=\{(x,y)\in \mathbb{R}^{2}: f(x,y,\mu)=0\}$
is the phase state of the reduced equation (\ref{reduce-1}).
A point in $\mathcal{C}_{\mu,0}$ with $\partial f/\partial x\neq 0$ is called a  regular point.
Otherwise it is called a contact point.
The set $\mathcal{C}_{\mu,0}$ is called  the critical set
and is called the  critical manifold if it is a submanifold of $\mathbb{R}^{2}$.
This set is useful in investigating the dynamics of the slow-fast system (\ref{fast-1}).
More specifically,
by the Fenichel theory \cite{Fenichel-79},
a normally hyperbolic manifold $\mathcal{M}_{\mu,0}$,
which is a compact submanifold $\mathcal{C}_{\mu,0}$ formed by regular points of a critical set $\mathcal{C}_{\mu,0}$,
is  perturbed to a slow manifold $\mathcal{M}_{\mu,\varepsilon}$  of slow-fast system (\ref{fast-1}) with $0<\varepsilon\ll 1$.
The stable and unstable manifolds of $\mathcal{M}_{\mu,0}$ are also persistent for a sufficiently small $\varepsilon$.

\begin{figure}[!htbp]
\centering
\subfigure[]{
\begin{minipage}[t]{0.22\linewidth}
\centering
\includegraphics[width=1.4in]{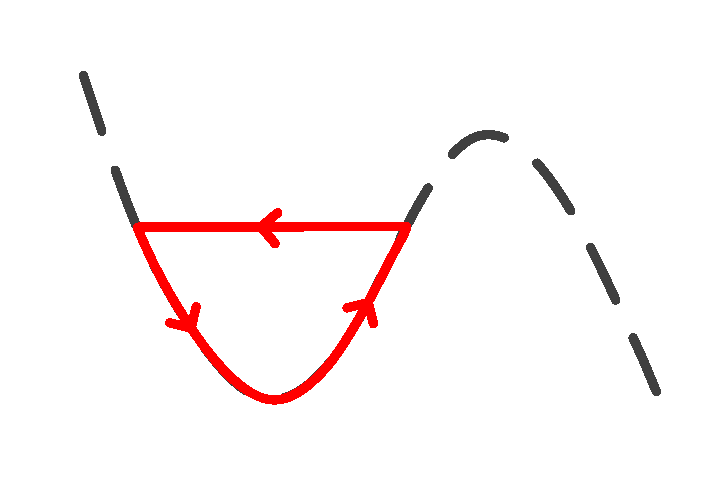}
\label{fig-canard-without-01}
\end{minipage}%
}%
\subfigure[]{
\begin{minipage}[t]{0.22\linewidth}
\centering
\includegraphics[width=1.1in]{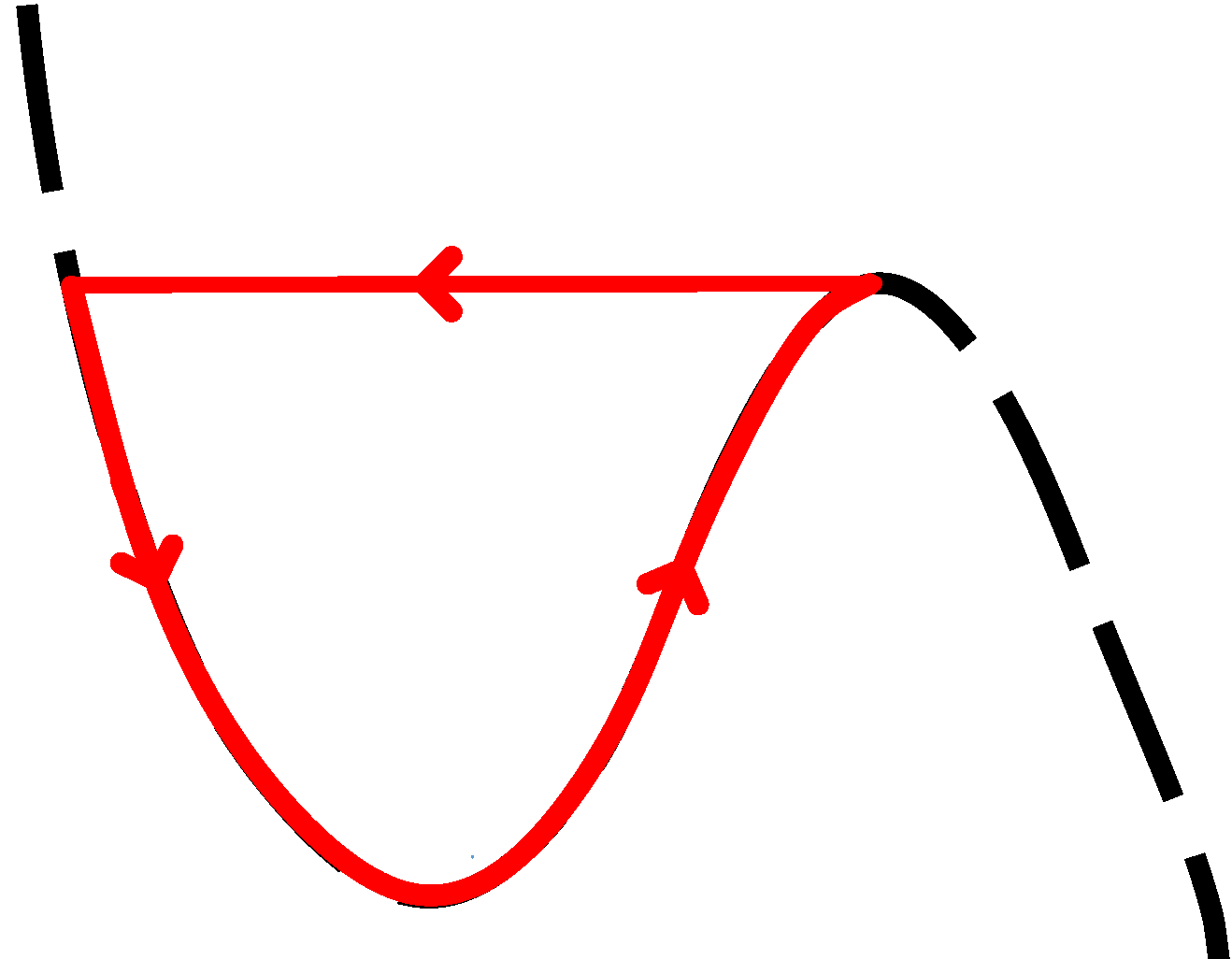}
\label{fig-canard-trans-01}
\end{minipage}%
}%
\subfigure[]{
\begin{minipage}[t]{0.22\linewidth}
\centering
\includegraphics[width=1.3in]{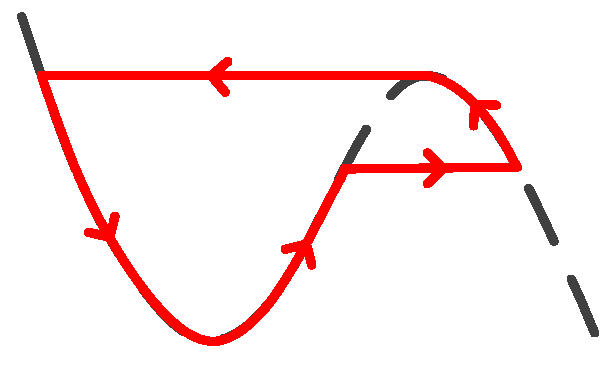}
\label{fig-canard-with-01}
\end{minipage}
}%
\subfigure[]{
\begin{minipage}[t]{0.22\linewidth}
\centering
\includegraphics[width=1.3in]{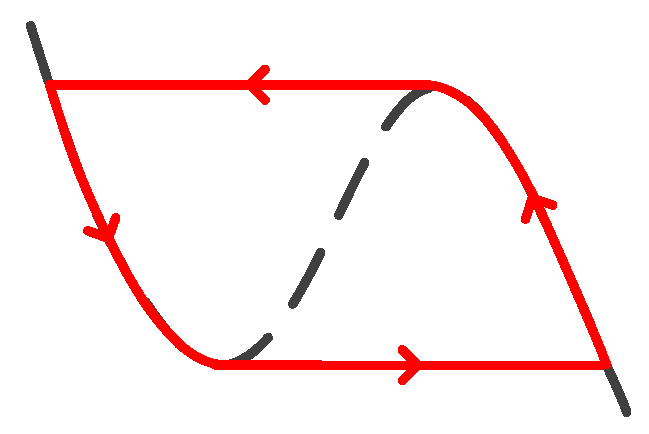}
\label{fig-relaxation-01}
\end{minipage}
}%
\centering
\caption{ \ref{fig-canard-without-01} Canard slow-fast cycle without head.
     \ref{fig-canard-trans-01} Transitory canard.
     \ref{fig-canard-with-01} Canard slow-fast cycle with head.
     \ref{fig-relaxation-01}  Singular relaxation cycle.
    }
        \label{F-singular-cycle}
\end{figure}

The preceding  results show the dynamics near the normally hyperbolic invariant manifolds.
However, non-hyperbolic points at which $\partial f/ \partial x=0$
widely appear in applications, such as the well-known van der Pol equation.
A contact point arising  in a critical manifold is one of the most common forms for the breakdown of normal hyperbolicity.
We analyze two different contact points in planar slow-fast systems,
that is the so-called jump point and canard point
\cite{Benoitetal-81,Dumortieretal-Roussarie-96,Krupa-Szmolyan-01SIMA},
which can induce relaxation oscillation and canard cycle, respectively.
Roughly speaking,
the reduced flow (\ref{reduce-1}) directs towards  a jump point
and passes through a canard point.
Relaxation oscillations and canard cycles can be seen as the  perturbations of slow-fast cycles
formed by gluing the orbits of the reduced system and the layer equations.
Four classical slow-fast cycles shown in Figure \ref{F-singular-cycle}.

Relaxation oscillations,
which perturb from their singular counterparts (see Fig. \ref{fig-relaxation-01}),
are periodic solutions which spend a long time along the slow manifold towards a jump point, jumps from this contact point,
spends a short time parallel to the unstable fibers towards another stable branch of the critical manifold,
follows the slow motion again until another jump point is reached,
and finally forms a closed loop via several similarly successive motions \cite{Grasman,Krupa-Szmolyan-01JDE}.
Canard cycle appearing near a canard point is a periodic solution
which is contained in the intersection of an attracting slow manifold and a repelling slow manifold
\cite{Benoitetal-81,Dumortieretal-Roussarie-96,Krupa-Szmolyan-01JDE}.
This  phenomenon is closely related to canard explosion \cite{Benoitetal-81,Krupa-Szmolyan-01JDE},
which is a transition from a small limit cycle of Hopf type via a family of canard cycles  to a relaxation oscillation.

\section{Model reduction and analysis of equilibria}
\label{sec-singuar-point}
\setcounter{equation}{0}
\setcounter{lemma}{0}
\setcounter{theorem}{0}
\setcounter{remark}{0}

In order to simplify calculations,
we first transform the THTN model into an equivalent system,
and then consider the properties of the equilibria in this system.
Letting
\begin{eqnarray*}
\left(M, P, t\right) \ \to \ \left(\frac{k_3}{8K\nu_{p}}y, \frac{1}{8K}x, \frac{1}{k_{3}}t\right),
\end{eqnarray*}
the THTN model is transformed into
\begin{eqnarray}
\label{2D-model-12}
\begin{aligned}
\frac{d x}{d t} &=x'=  y-\psi_{1}(x),
\\
\frac{d y}{d t} &=y'= \varepsilon\left(\psi_{2}(x)-y\right),
\end{aligned}
\end{eqnarray}
where
\begin{eqnarray} \label{df-h1-h2}
\psi_{1}(x)=\frac{b_{1}\phi(x)+b_2x}{a+x}+x,\ \ \psi_{2}(x)=\frac{v}{c+(x-\phi(x))^{2}},\ \
\phi(x)=2(\sqrt{1+x}-1), \ \ \ \ \ x\geq 0,
\end{eqnarray}
and the positive parameters $a, b_{1}, b_{2}, c, \varepsilon, v$ are given by
\begin{eqnarray*}
a=8J_{P}K,\ \
b_{1}=\frac{8(k_{1}-k_{2})K}{k_3},\ \
b_{2}=\frac{8k_{2}K}{k_3},\ \
c=256K^{2}P_{c}^{2}, \ \
\varepsilon=\frac{k_{m}}{k_{3}},\ \
v=\frac{2048\nu_{m}\nu_{p}P_{c}^{2}K^{3}}{k_3k_m}.
\end{eqnarray*}
Our goal is to study the effects of the rates of mRNA degradation and synthesis on
the periodic behaviors in the THTN model.
For this reason, throughout this paper  we vary the parameters $k_{m}$ and $v_{m}$,
and fix the remaining parameters in the THTN model.
Additionally, we also assume that the rate of mRNA degradation
is proportional to that of mRNA synthesis.
Then the parameters $v$ and $\varepsilon$ are independent of each other and vary,
and other parameters in system (\ref{2D-model-12}) are fixed.

Define
\begin{eqnarray}\label{df-h1-h2-0}
\psi(x):=\psi_{1}(x)-\psi_{2}(x)\ \ \mbox{ for }\ \ x\geq 0.
\end{eqnarray}
Concerning $\psi_{i}$ and $\psi$,
we have the following two lemmas.

\begin{lemma}
\label{lm-psi-1-prpty}
Let  $\psi_{1}$ be defined by (\ref{df-h1-h2}).
Then the second derivative $\psi_{1}''$ of $\psi_{1}$ has a unique positive zero $x_{+}=u_{+}^{2}+2u_{+}$,
where $u_{+}$ is the unique positive zero of the function $\phi_{1}$ defined as
\begin{eqnarray}\label{df-phi-1}
\phi_1(u)=3b_1(u+1)^{4}-(8b_1+4ab_2)(u+1)^{3}+6b_1(1-a)(u+1)^{2}-b_1(a-1)^{2},
\end{eqnarray}
and the following statements hold:

\item{\bf (i)}
$\psi_{1}(0)=0$, $\psi_{1}(x)>0$ for $x>0$ and $\psi_{1}(x)/x\to 1$ as $x\to +\infty$.

\item{\bf (ii)}
$\psi_{1}^{'}(0)=(b_1+b_2)/a+1$, $\psi_{1}^{'}(x)\to 1$ as $x\to +\infty$
and  $\psi_{1}^{'}$ admits the following trichotomies:
\begin{enumerate}
\item[{\bf (ii.1)}]
if $\psi_{1}^{'}(x_{+})>0$, then $\psi_{1}^{'}(x)>0$ for $x\geq 0$.

\item[{\bf (ii.2)}]
if $\psi_{1}^{'}(x_{+})=0$,
then $\psi_{1}^{'}(x)\geq 0$ for $x\geq 0$, and $x_{+}$ is the unique positive zero of $\psi_{1}^{'}$.

\item[{\bf (ii.3)}]
if $\psi_{1}^{'}(x_{+})<0$,
then $\psi_{1}^{'}$ has exactly two zeros $x_{m}$ and $x_{M}$ with $0<x_{m}<x_{+}<x_{M}$,
and  $\psi_{1}^{'}$ satisfies that $\psi_{1}^{'}(x)>0$ for $0<x<x_{m}$ and $x>x_{M}$,
$\psi_{1}^{'}(x)<0$ for $x_{m}<x<x_{M}$.
\end{enumerate}

\item{\bf (iii)}
$\psi_{1}^{''}(x)<0$ for $x\in [0, x_{+})$ and $\psi_{1}^{''}(x)>0$ for $x\in (x_{+},+\infty)$.
\end{lemma}
{\bf Proof.}
Set $u=\sqrt{1+x}-1$ for $x\geq 0$. Then $x=u^{2}+2u$  for $u\geq 0$.
By a direct computation, we have that
$$
2(u+1)^{3}(u^{2}+2u+a)^{3}\psi_{1}^{''}(x(u))=\phi_1(u),
$$
where $\phi_1$ is defined by (\ref{df-phi-1}).
Then by a standard analysis, we obtain this lemma.
\hfill$\Box$

In {\bf (i)} of Theorem \ref{prop-dynamics-1}
we will see that the dynamics of \eqref{2D-model-12} with $\psi_{1}^{'}(x_{+})\geq 0$ is simple.
Consequently, with no confusion,
we always assume that $\psi_{1}^{'}(x_{+})<0$. So the graph of $\psi_{1}$ is $S$-shaped.

\begin{lemma}
\label{lm-psi-0-2-prpty}
Let the functions $\psi_{2}$ and $\psi$ be defined by (\ref{df-h1-h2}) and (\ref{df-h1-h2-0}), respectively.
Then the function $\psi_{2}$ has the following properties:
\item{\bf (i)}
$\psi_{2}(0) =v/c$, $0<\psi_{2}(x)\leq v/c$ for $x\geq 0$, and $\psi_{2}(x)\to 0$ as $x\to +\infty$.

\item{\bf (ii)}
$\psi_{2}'(0) =0$, $-v/(c\sqrt{c})\leq \psi_{2}'(x)<0$ for $x>0$, and $\psi_{2}'(x)\to 0$ as $x\to +\infty$.

\item{\bf (iii)}
the second derivative $\psi_{2}^{''}$ of $\psi_2$ has exactly one zero $x_{1}\in(0,+\infty)$,
which is the unique positive  root of equation
$6(\sqrt{x+1}-1)^{5}+5(\sqrt{x+1}-1)^{4}-2c\sqrt{x+1}-c=0$,
and $\psi_{2}^{''}(x)<0$ for $0<x<x_{1}$ and $\psi_{2}^{''}(x)>0$ for $x>x_{1}$.\\
And the function $\psi$ has the following properties:
\item{\bf (iv)} for each positive parameters $a$, $b_{1}$, $b_{2}$, $c$, $\varepsilon$ and $v$,
the function $\psi$ has at least one positive zero and at most three positive zeros.

\item{\bf (v)}
if the function $\psi$ has precisely two positive zeros $x=\widetilde{x}_{0}$ and $x=\widetilde{x}_{1}$ with $\widetilde{x}_{0}<\widetilde{x}_{1}$,
then  either $\omega=\widetilde{x}_{0}$ or $\omega=\widetilde{x}_{1}$ satisfies that
$\psi(\omega)=\psi^{'}(\omega)=0$ and $\psi^{''}(\omega)\neq 0$.
\end{lemma}
{\bf Proof.}
By a standard analysis, the properties of $\psi_{2}$ can be obtained,
thus the  proof is omitted.

To obtain the properties on  $\psi$,
let $u=\sqrt{1+x}-1$ for $x\geq 0$. Then we have
\begin{eqnarray*}
\begin{aligned}
&(u^{2}+2u+a)(u^{4}+a)\psi(x(u))\\
&\ \ \ =(u^{4}+4u^{3}+(a+b_{2}+4)u^{2}+2(a+b_{1}+b_{2})u)(u^{4}+c)
-v(u^{2}+2u+a):=\phi_{2}(u).
\end{aligned}
\end{eqnarray*}
Since $\phi_{2}(0)=-av<0$ and $\phi_{2}(u)\to +\infty$ as $u\to +\infty$,
then by continuity there exists at least one positive zero for the function $\psi$.
Since the third derivative of $\phi_{2}$ is in the form
\begin{eqnarray*}
\phi_{2}^{(3)}(u)=336u^{5}+840u^{4}+120(a+b_{2}+4)u^{3}+120(a+b_{1}+b_{2})u^{2}+24c u+24c,
\end{eqnarray*}
and $\phi_{2}^{(3)}(u)>0$ for $u\geq 0$,
then $\phi$ has at most three positive zeros.
Thus {\bf (iv)} is proved.
By studying the properties of $\phi_{2}$,
we can obtain  {\bf (v)}.
Therefore, the proof is now complete.
\hfill$\Box$

Under the assumption that $\psi_{1}^{'}(x_{+})<0$,
we observe that the graph of the function $\psi_{1}$ is $S$-shaped.
To consider the properties of the equilibria in (\ref{2D-model-12}),
let $L=L^{0}\cup L^{1}$, $R=R^{0}\cup R^{1}$ and $M=\{(x,y): y=\psi_{1}(x), x_{m}< x < x_{M}\}$,
where the sets
\begin{eqnarray*}
L^{0} \!\!\!& =\{(x_{m},\psi_{1}(x_{m}))\}, \ \ \  L^{1} \!\!\!\!& =\{(x,y): y  =\psi_{1}(x),\  0\leq x < x_{m}\}, \\
R^{0} \!\!\!&=\{(x_{M},\psi_{1}(x_{M}))\},  \ \ \  R^{1} \!\!\!\!&=\{(x,y): y =\psi_{1}(x),\  x> x_{M}\}.
\end{eqnarray*}
We now define symbolic sequences to indicate the numbers
and relative positions of the equilibria on the graph of $\psi_1$.
We use, for example, the symbolic sequence $LMR$ to represent that
$\psi_{2}$ intersects $\psi_{1}$ at points in the sets $L$, $M$ and $R$ in order as the independent variable $x$ increases,
other symbolic sequences are similarly defined.
These symbolic sequences are referred to as the  intersection point sequences.

We next consider all possible intersection point sequences in the case  $\psi_{1}^{'}(x_{+})<0$,
which is useful in the proof for the main results in the  low degradation rate case.

\begin{lemma}
\label{distribution-inter}
Suppose that the function $\psi_{1}$  satisfies $\psi_{1}^{'}(x_{+})<0$,
where the function $\psi_{1}$ and the constant $x_{+}$ are defined as  in Lemma \ref{lm-psi-1-prpty}.
Then the intersection point sequences have the following different types (see Figure \ref{fg-all-cases}):
\begin{figure}[!htbp]
\centering
\subfigure[$L^{1}$.]{
\begin{minipage}[t]{0.23\linewidth}
\centering
\includegraphics[width=1.2in]{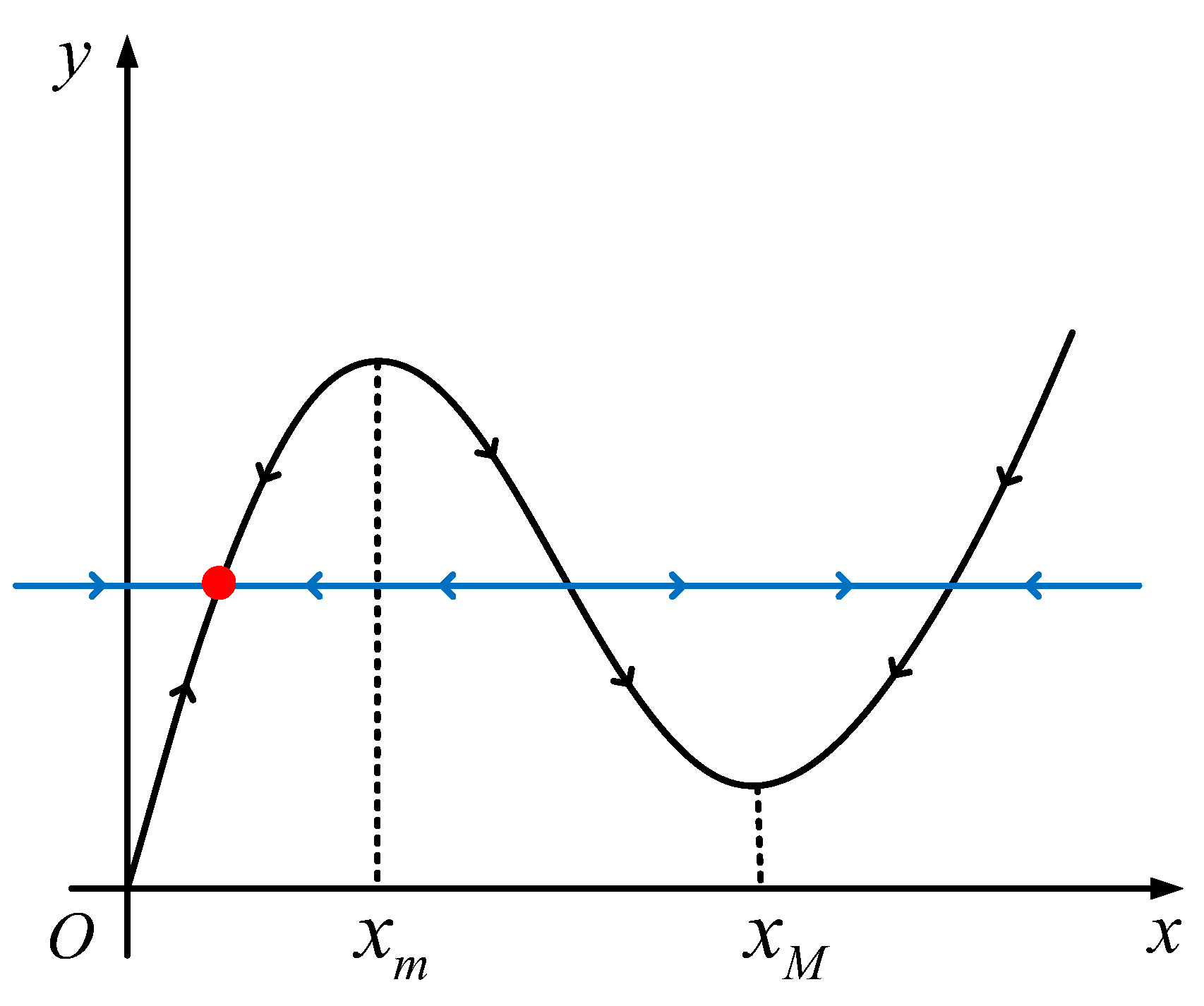}
\label{f-1-1}
\end{minipage}%
}%
\subfigure[$L^{0}$.]{
\begin{minipage}[t]{0.23\linewidth}
\centering
\includegraphics[width=1.2in]{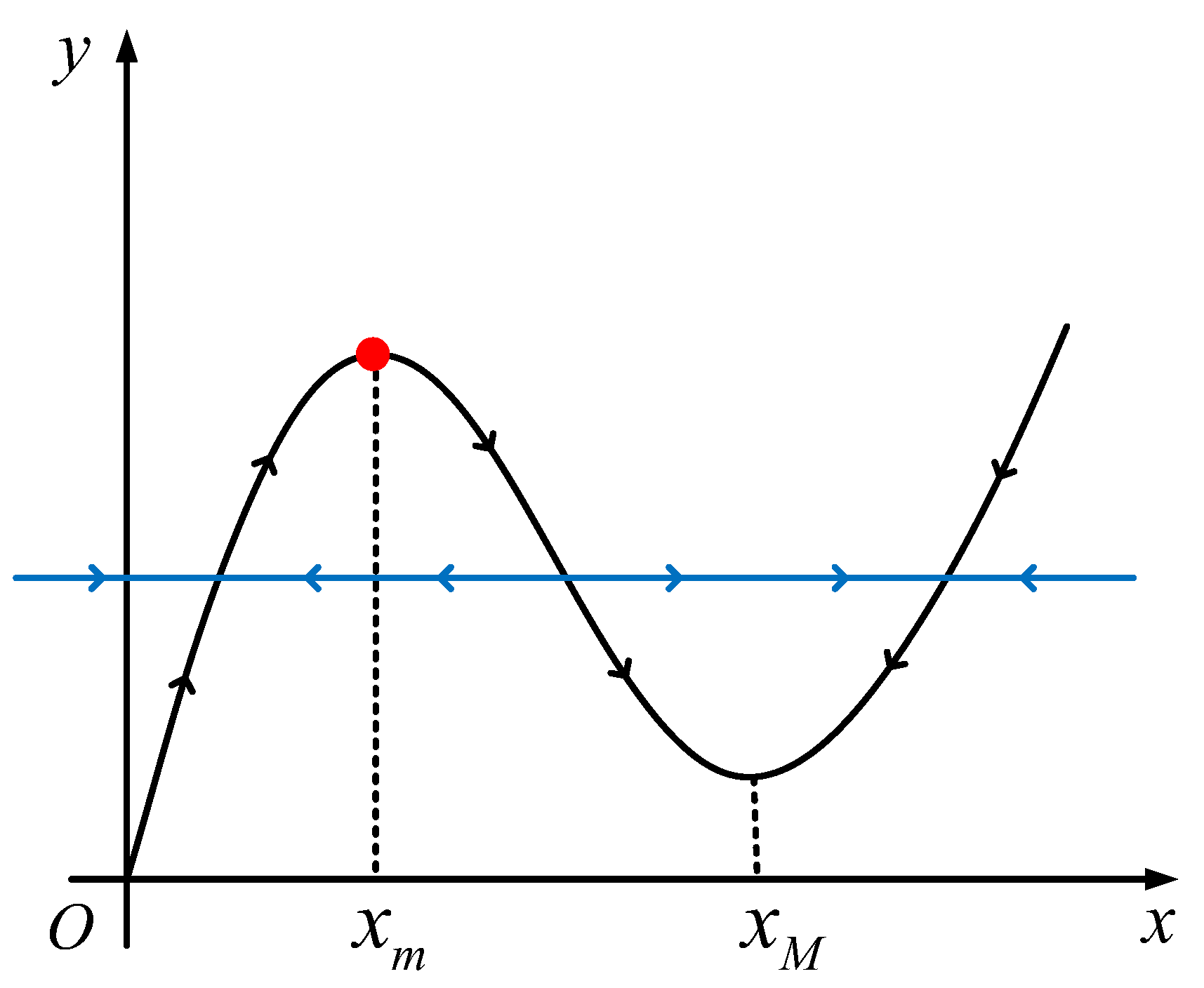}
\label{f-1-2}
\end{minipage}%
}%
\subfigure[$M$.]{
\begin{minipage}[t]{0.23\linewidth}
\centering
\includegraphics[width=1.2in]{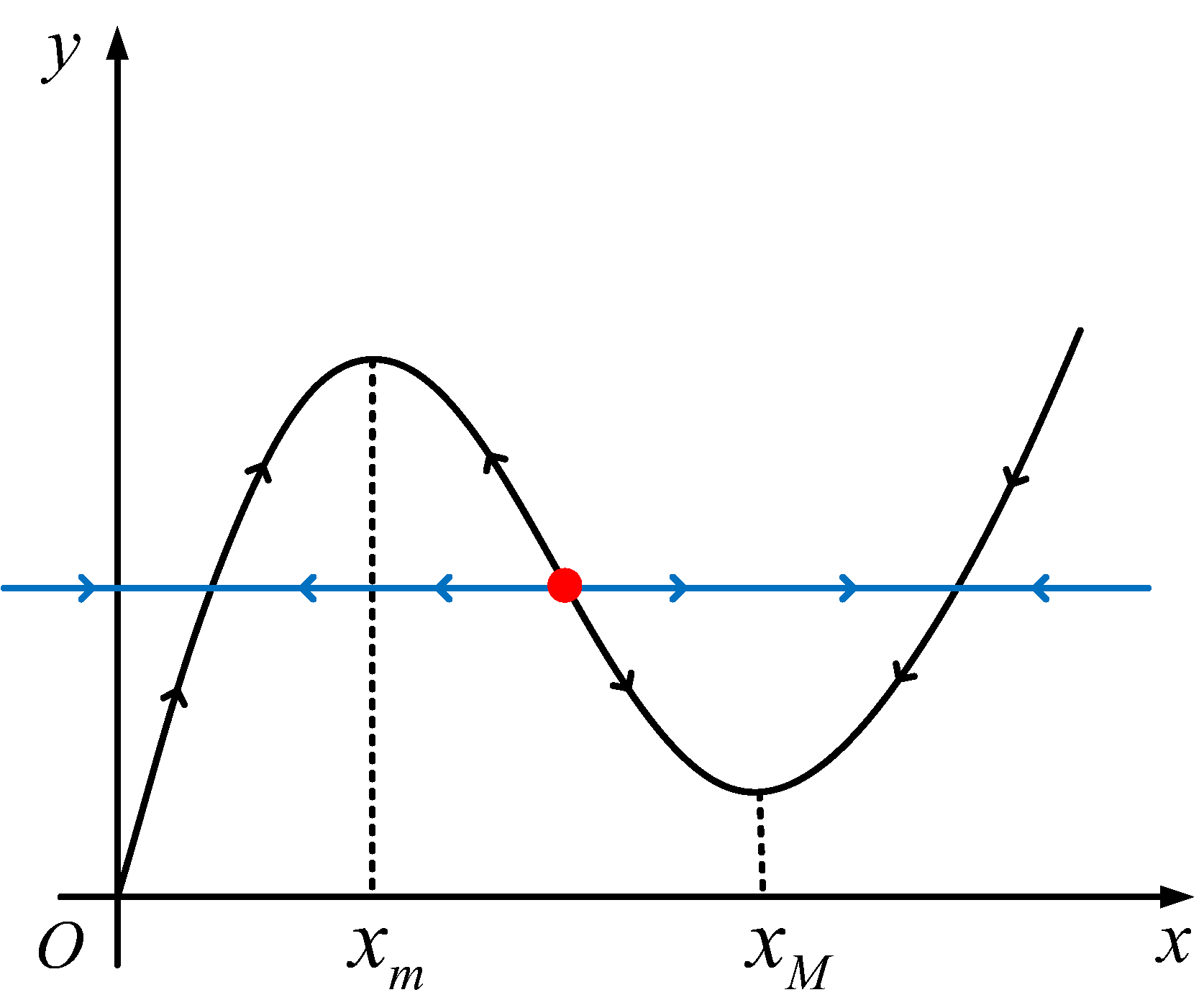}
\label{f-1-3}
\end{minipage}
}%
\subfigure[$R^{0}$.]{
\begin{minipage}[t]{0.23\linewidth}
\centering
\includegraphics[width=1.2in]{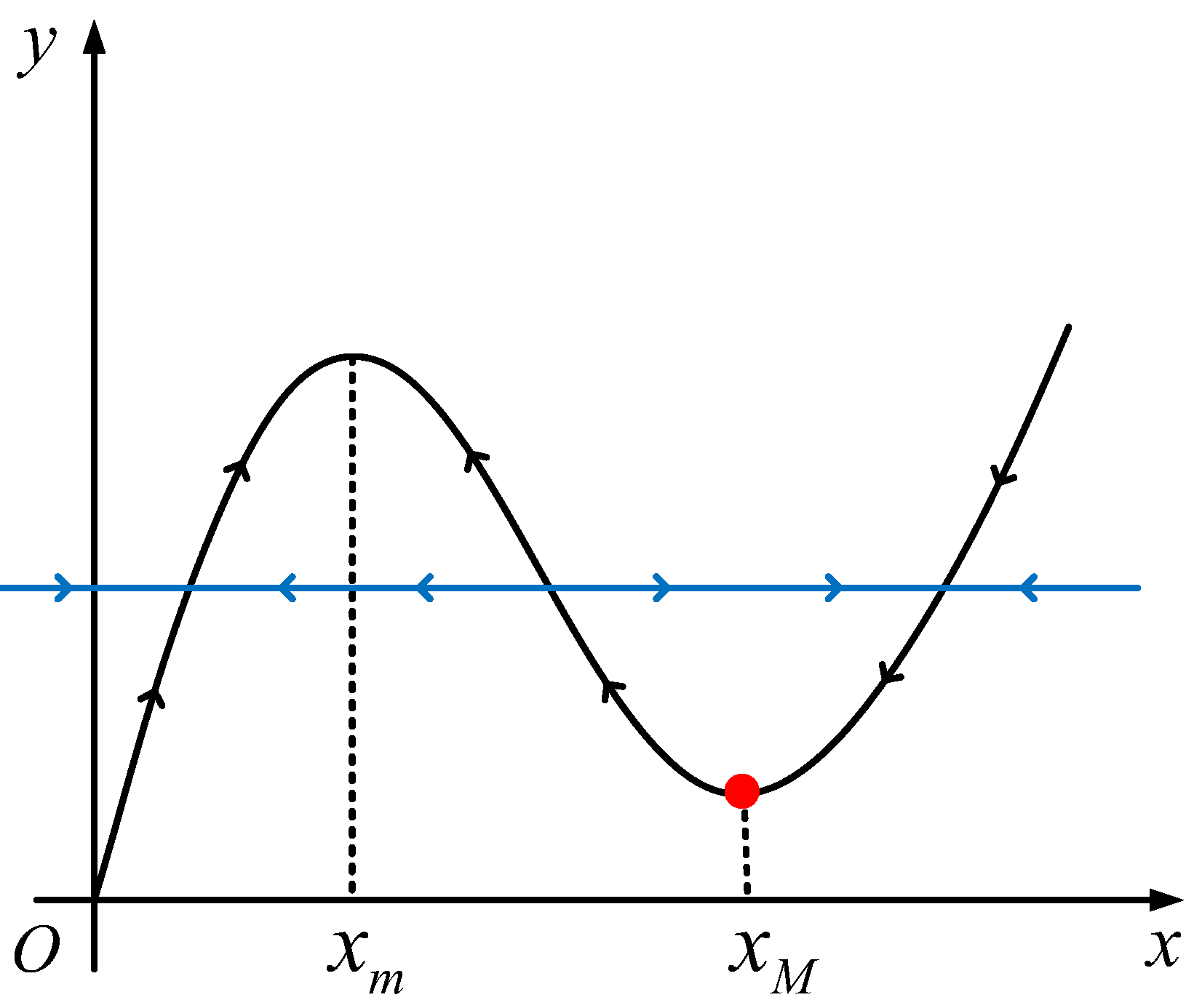}
\label{f-1-4}
\end{minipage}
}%
\\
\subfigure[$R^{1}$.]{
\begin{minipage}[t]{0.23\linewidth}
\centering
\includegraphics[width=1.2in]{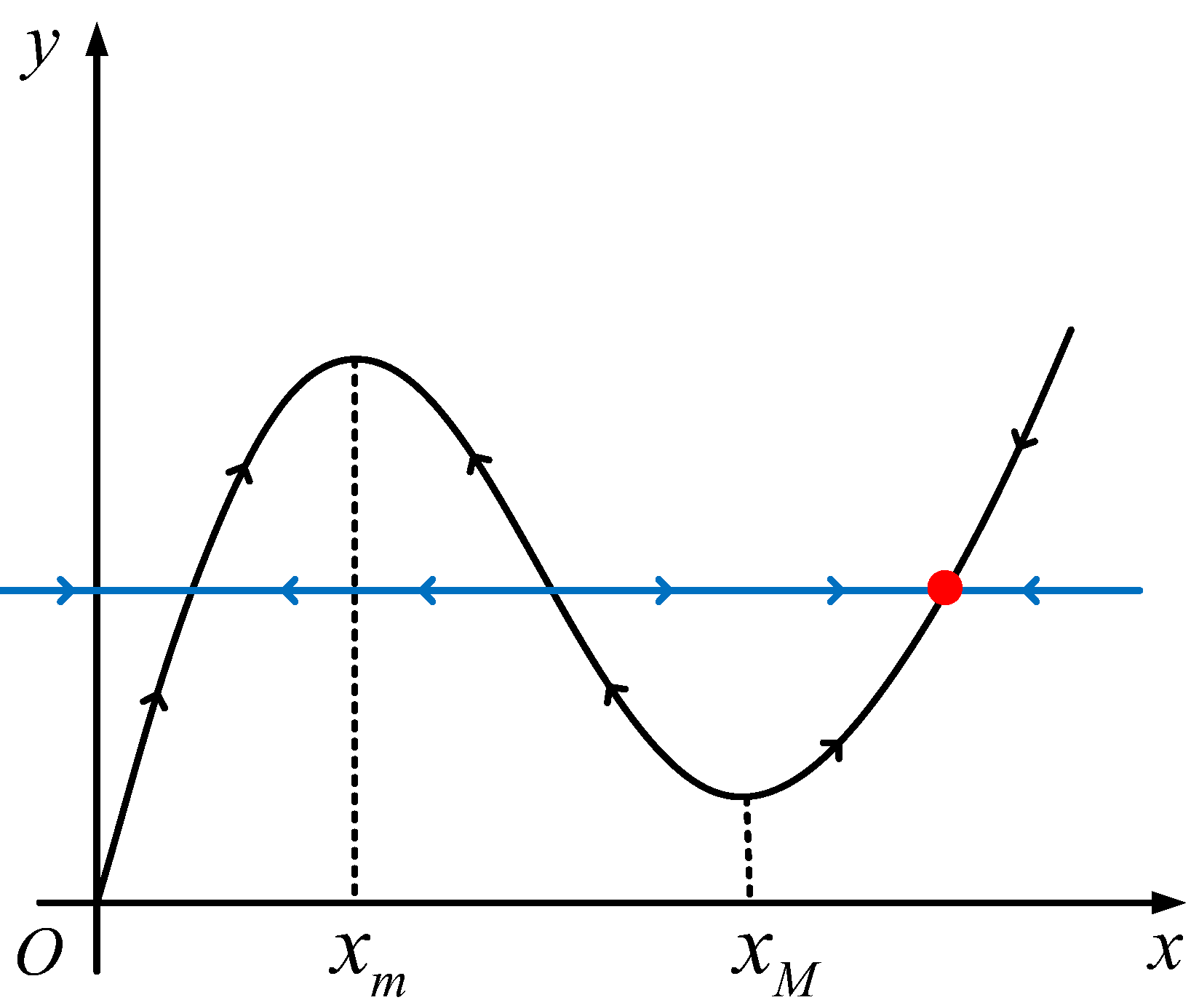}
\label{f-1-5}
\end{minipage}
}%
\subfigure[$L^{0}M$.]{
\begin{minipage}[t]{0.23\linewidth}
\centering
\includegraphics[width=1.2in]{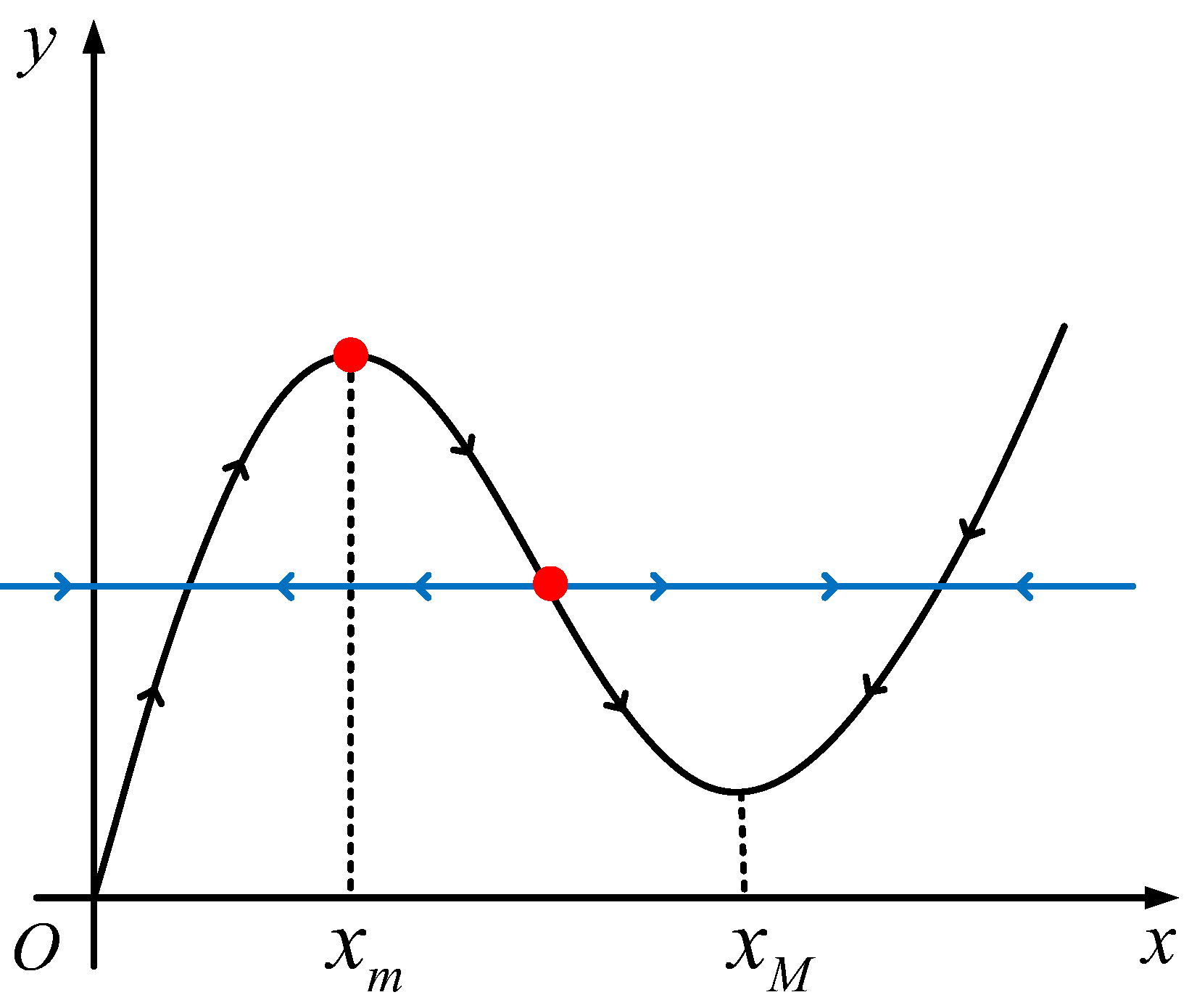}
\label{f-2-1}
\end{minipage}
}%
\subfigure[$L^{1}M$.]{
\begin{minipage}[t]{0.23\linewidth}
\centering
\includegraphics[width=1.2in]{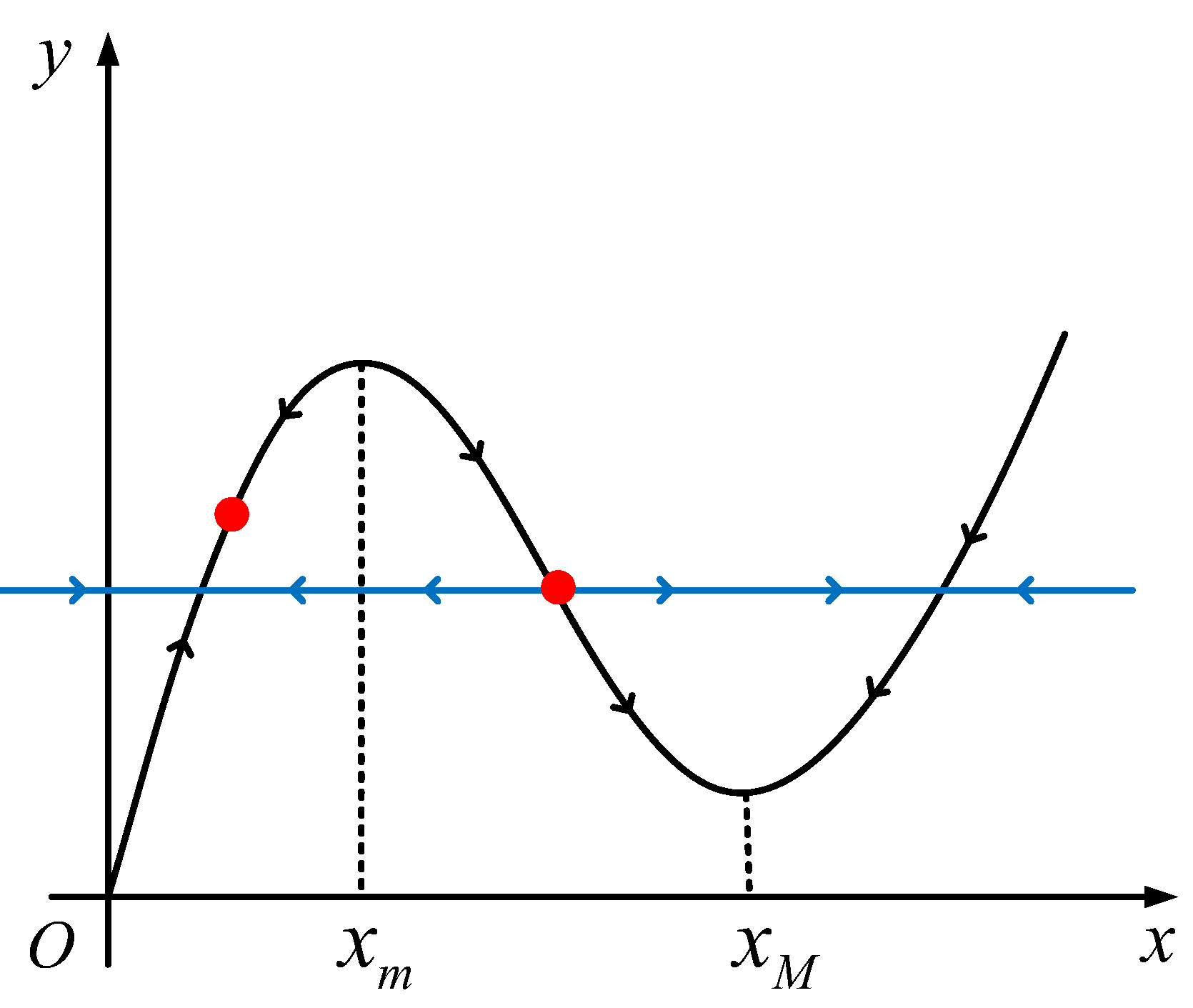}
\label{f-2-2}
\end{minipage}
}%
\subfigure[$MM$.]{
\begin{minipage}[t]{0.23\linewidth}
\centering
\includegraphics[width=1.2in]{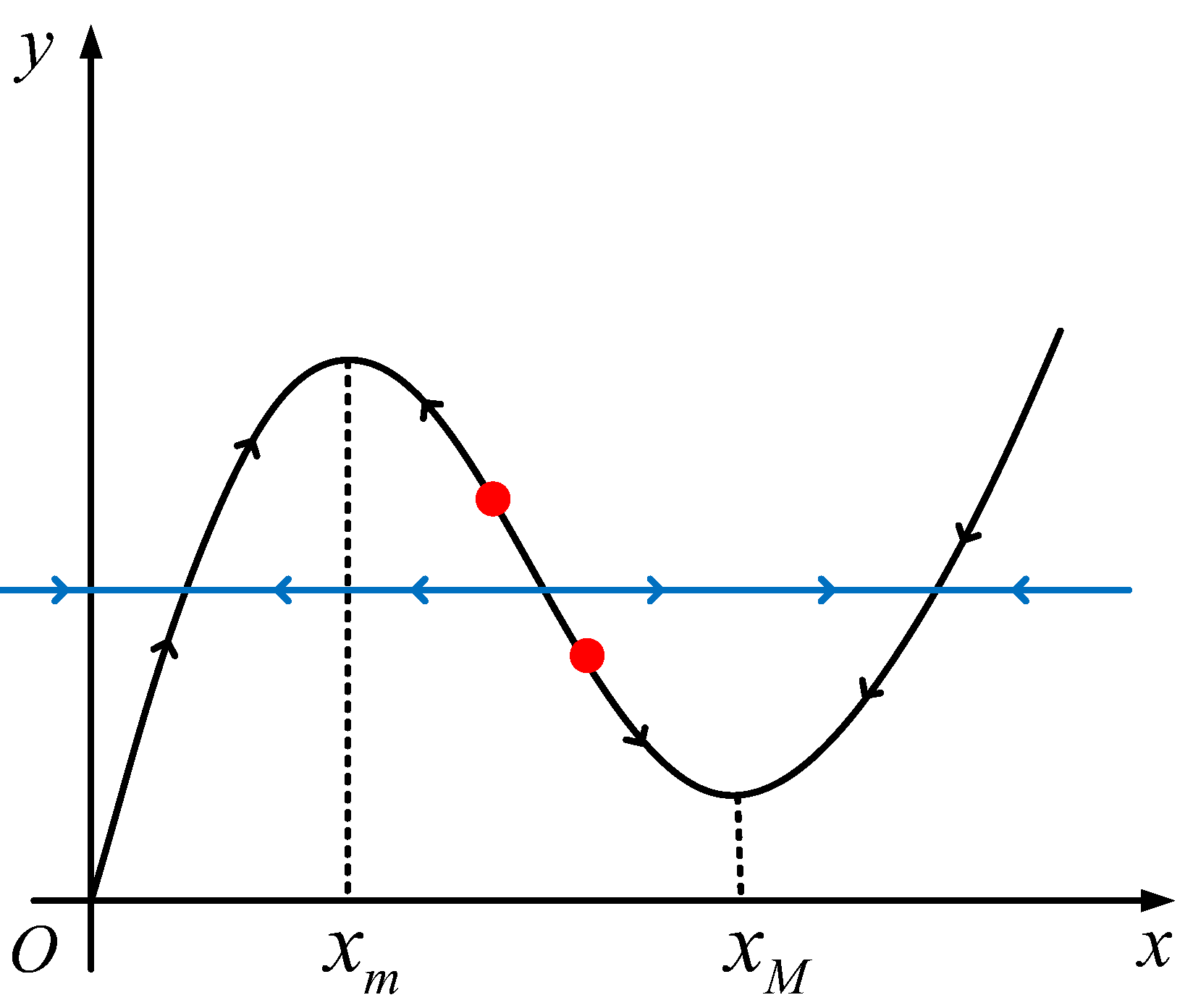}
\label{f-2-3}
\end{minipage}
}%
\\
\subfigure[$MR^{0}$.]{
\begin{minipage}[t]{0.23\linewidth}
\centering
\includegraphics[width=1.2in]{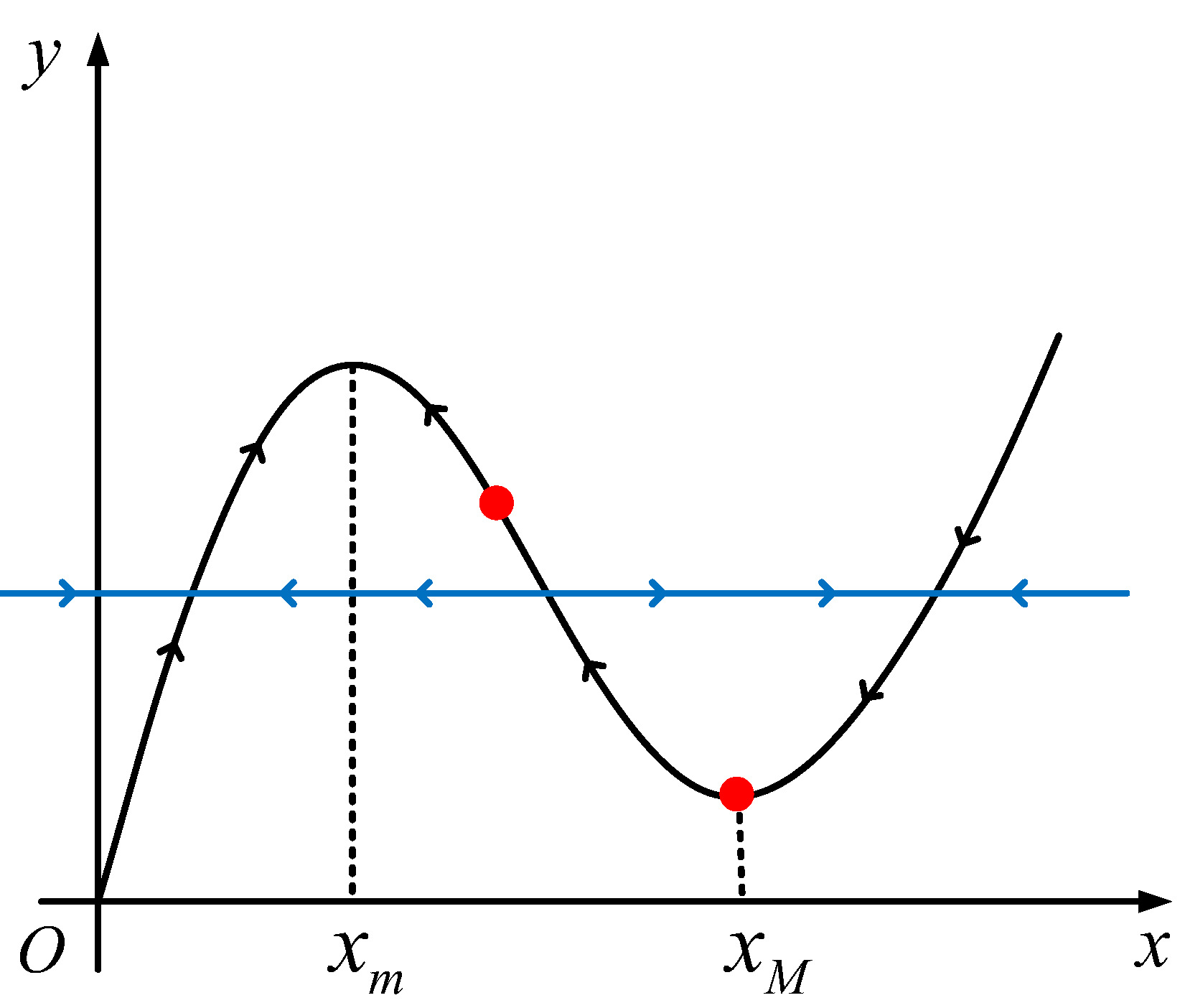}
\label{f-2-4}
\end{minipage}
}%
\subfigure[$MR^{1}$.]{
\begin{minipage}[t]{0.23\linewidth}
\centering
\includegraphics[width=1.2in]{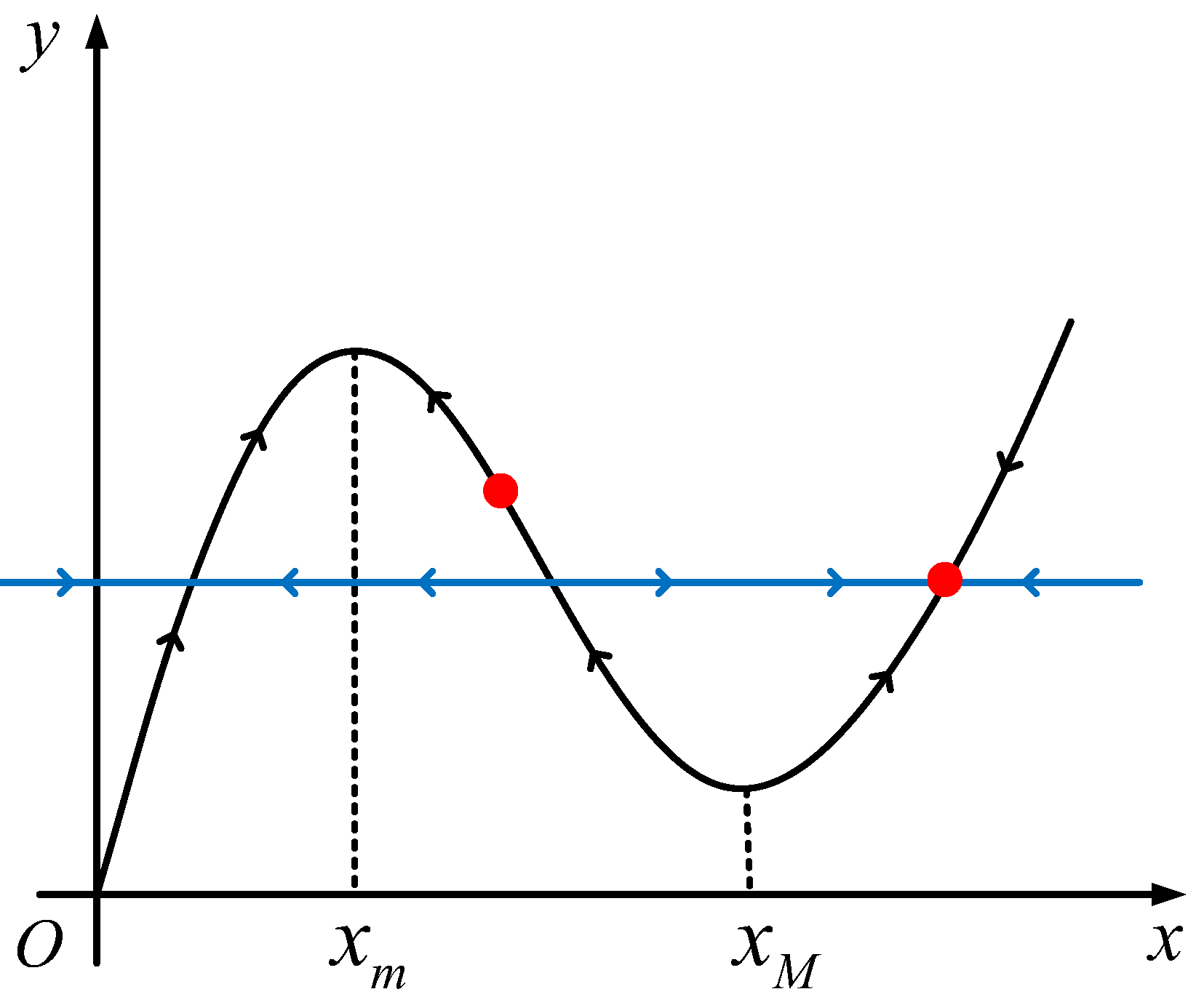}
\label{f-2-5}
\end{minipage}
}%
\subfigure[$L^{0}MR^{0}$.]{
\begin{minipage}[t]{0.23\linewidth}
\centering
\includegraphics[width=1.2in]{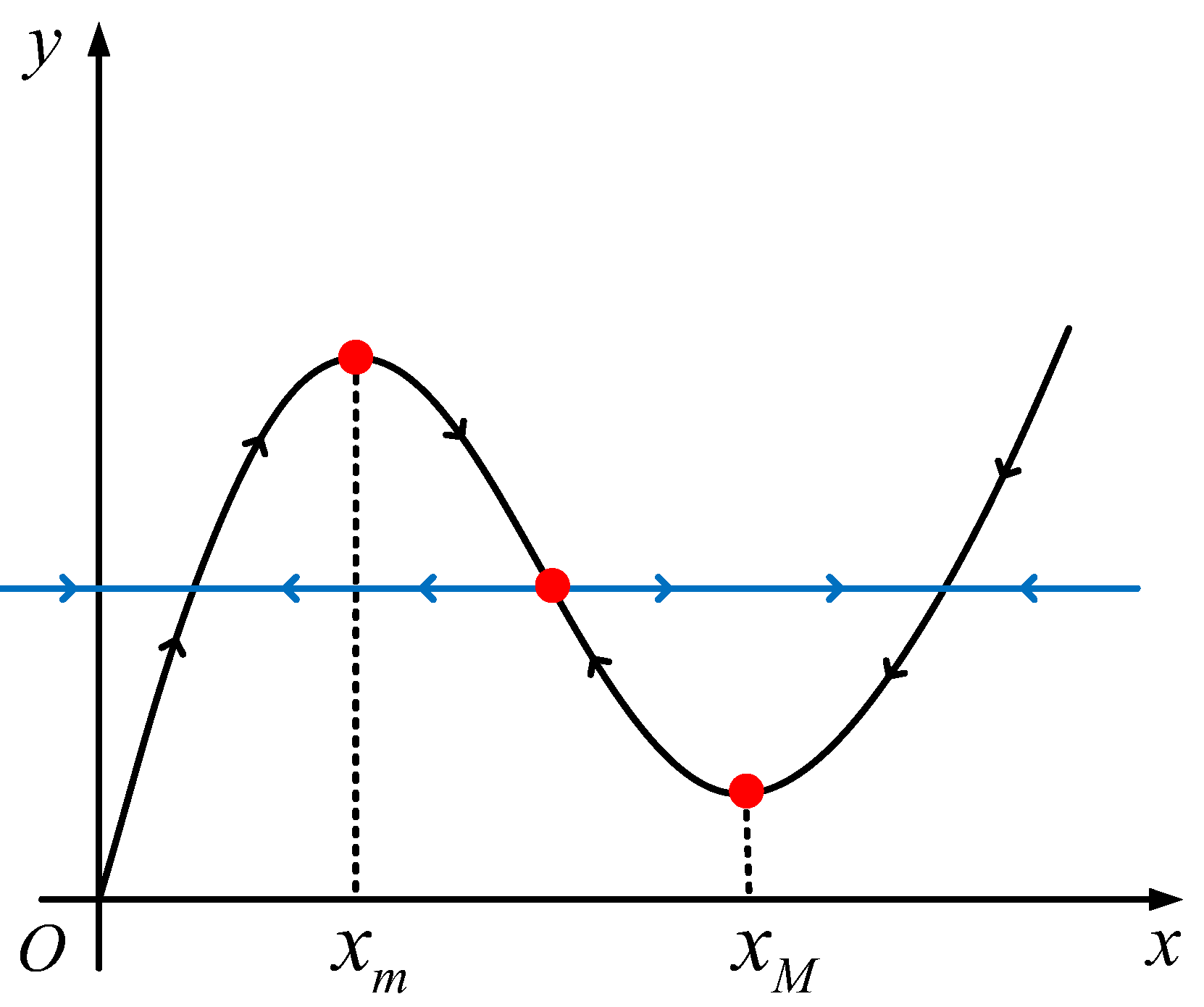}
\label{f-3-1}
\end{minipage}
}%
\subfigure[$L^{0}MR^{1}$.]{
\begin{minipage}[t]{0.23\linewidth}
\centering
\includegraphics[width=1.2in]{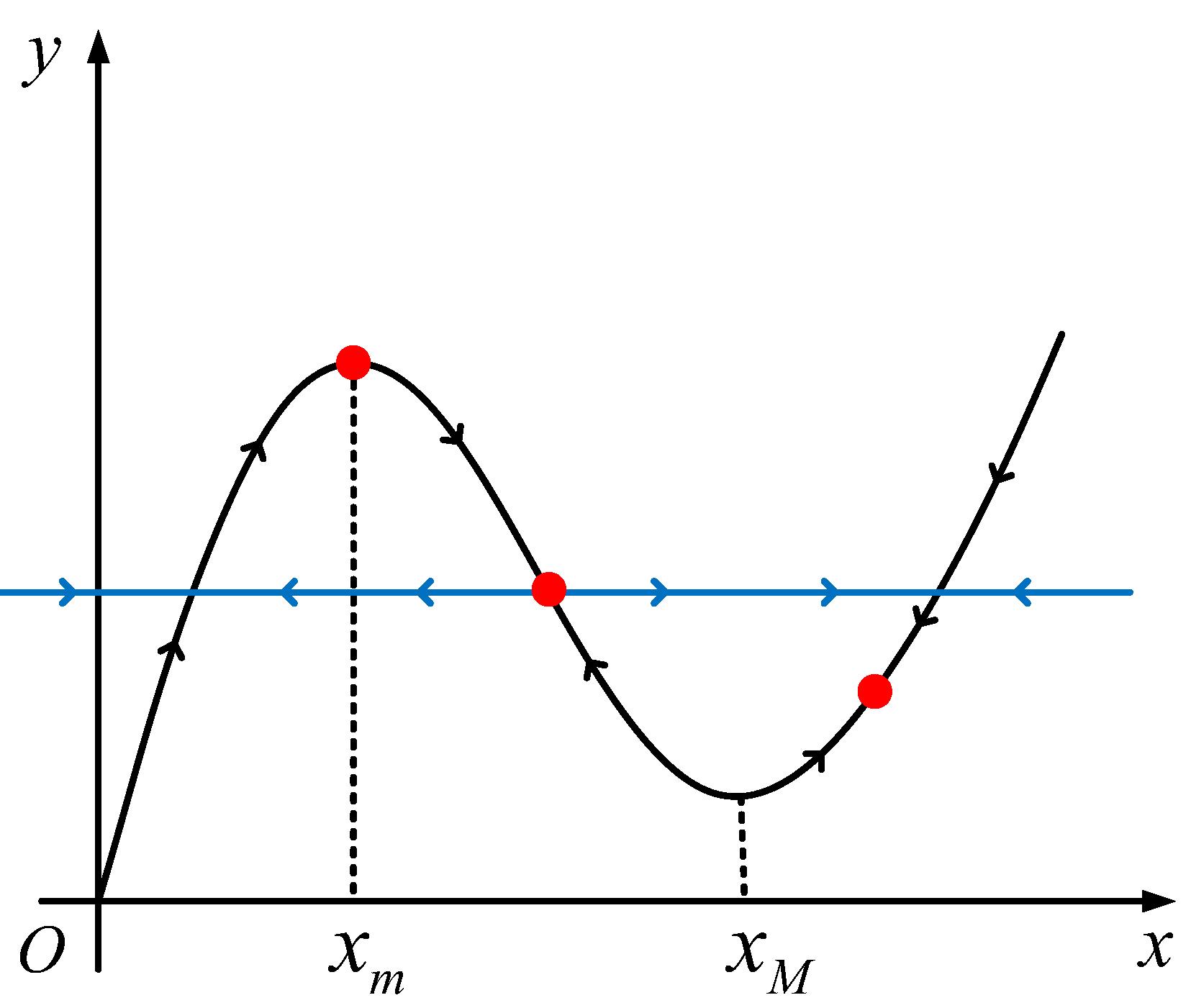}
\label{f-3-2}
\end{minipage}
}%
\\
\subfigure[$L^{1}MR^{0}$.]{
\begin{minipage}[t]{0.23\linewidth}
\centering
\includegraphics[width=1.2in]{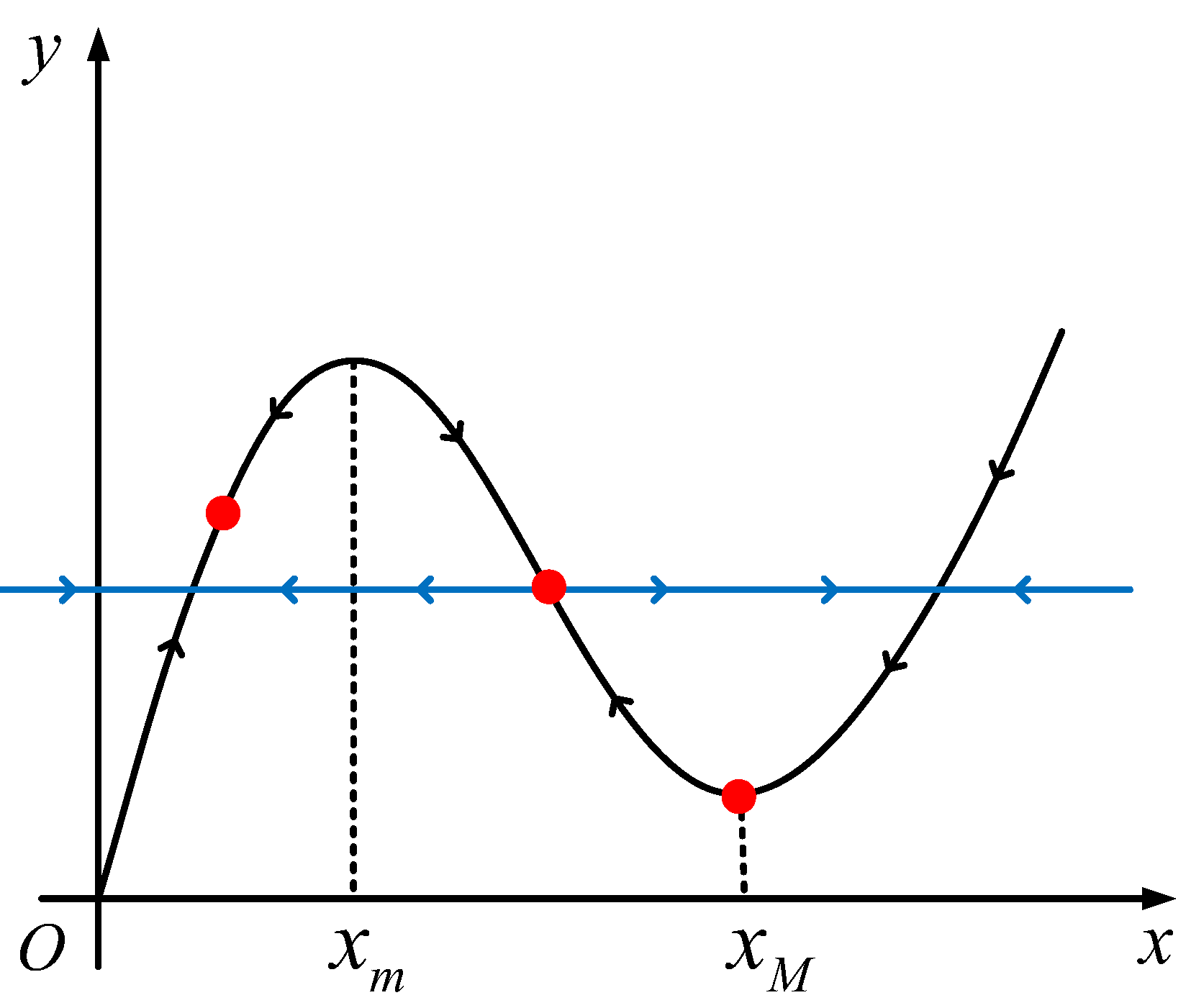}
\label{f-3-3}
\end{minipage}
}%
\subfigure[$L^{1}MR^{1}$.]{
\begin{minipage}[t]{0.23\linewidth}
\centering
\includegraphics[width=1.2in]{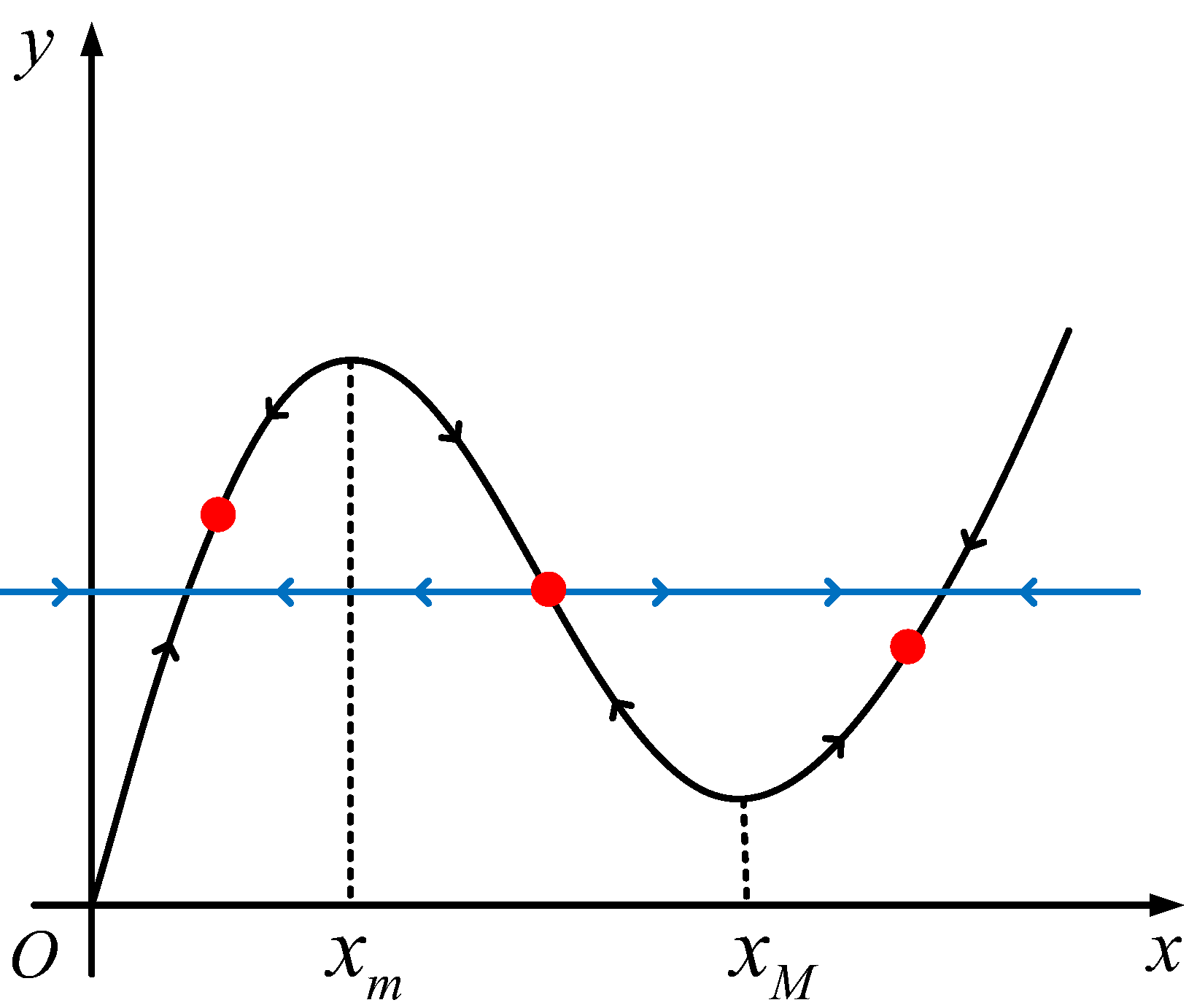}
\label{f-3-4}
\end{minipage}
}%
\subfigure[$L^{0}MM$.]{
\begin{minipage}[t]{0.23\linewidth}
\centering
\includegraphics[width=1.2in]{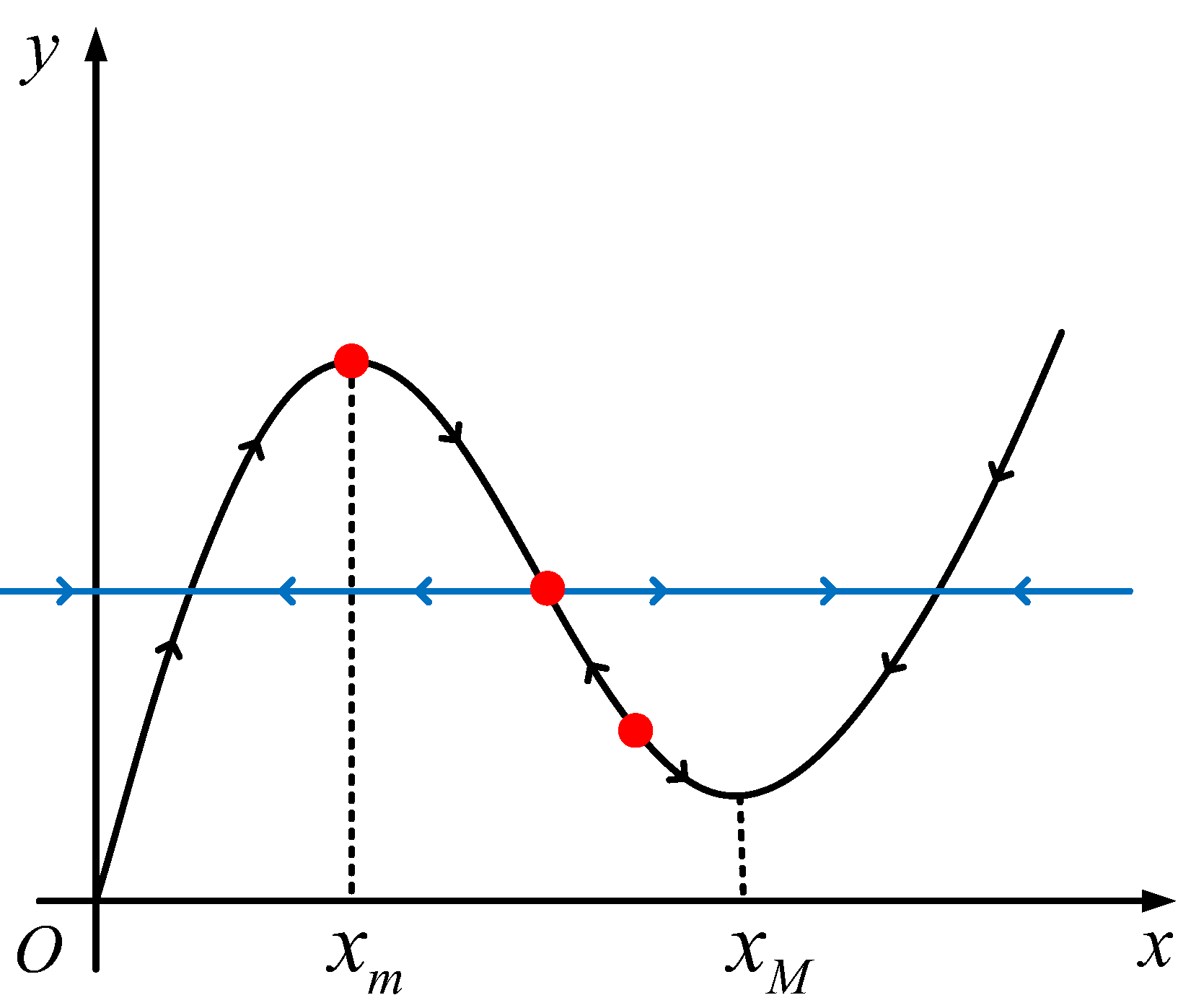}
\label{f-3-5}
\end{minipage}
}%
\subfigure[$L^{1}MM$.]{
\begin{minipage}[t]{0.23\linewidth}
\centering
\includegraphics[width=1.2in]{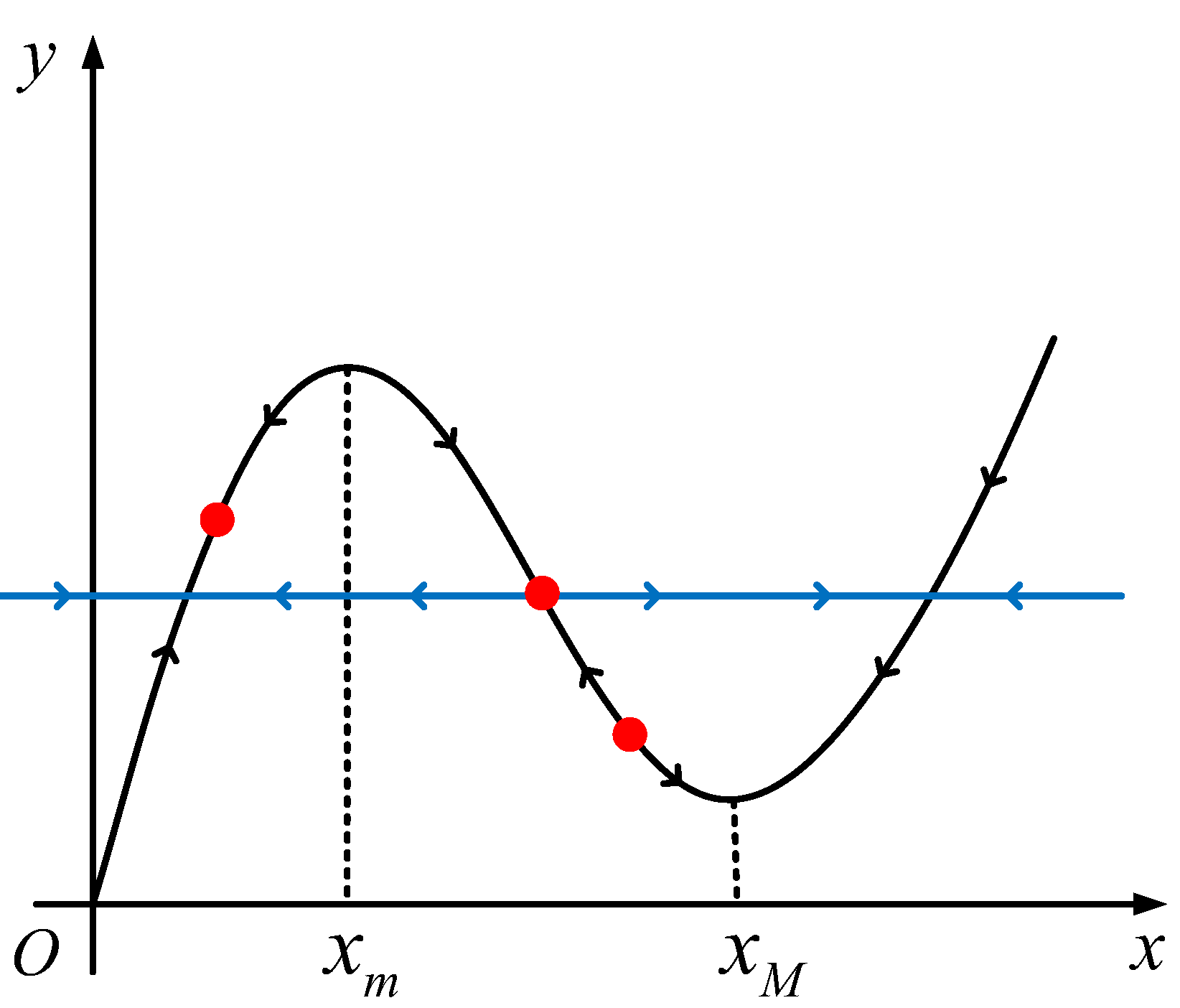}
\label{f-3-6}
\end{minipage}
}%
\\
\subfigure[$MMM$.]{
\begin{minipage}[t]{0.23\linewidth}
\centering
\includegraphics[width=1.2in]{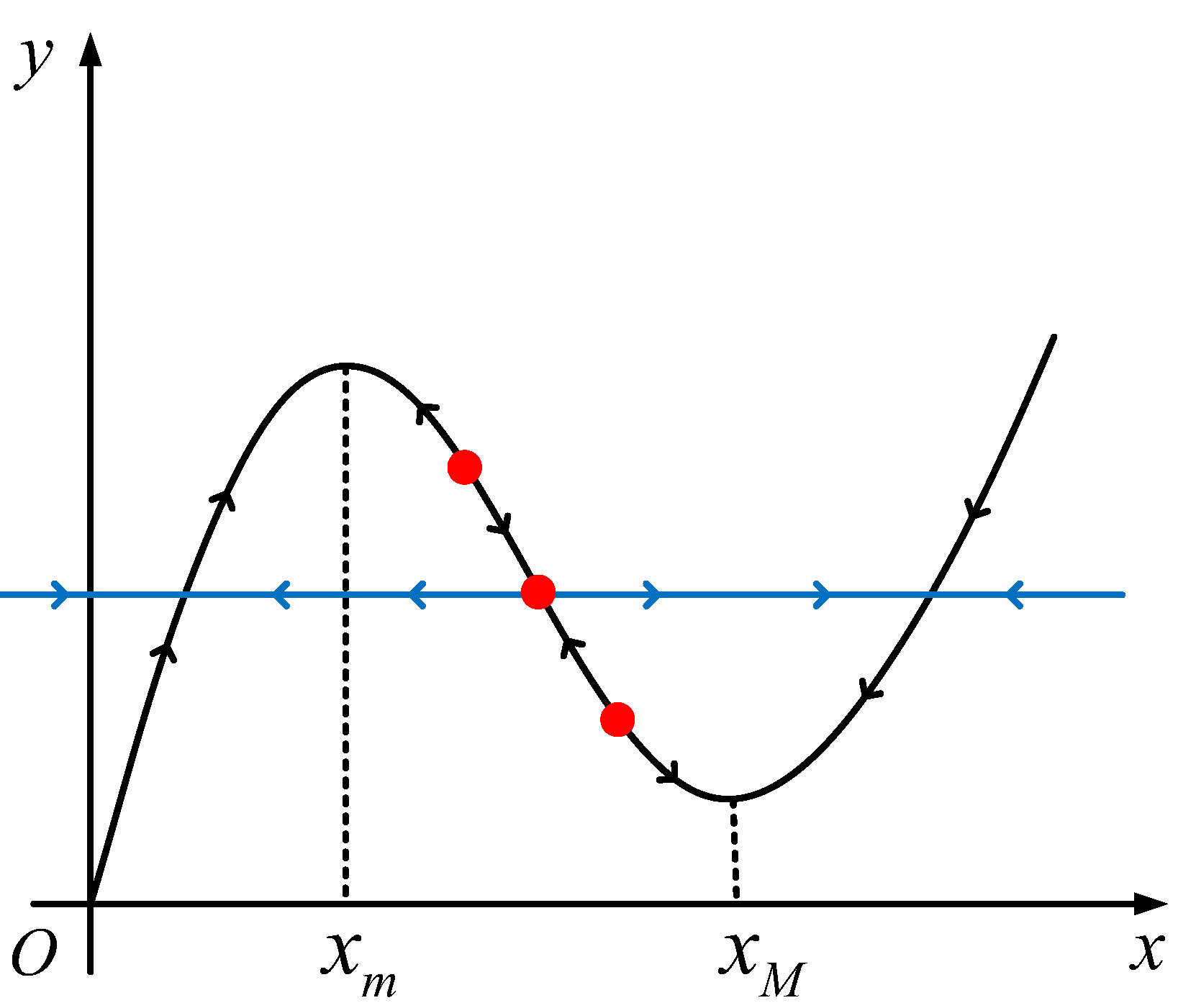}
\label{f-3-7}
\end{minipage}
}%
\subfigure[$MMR^{0}$.]{
\begin{minipage}[t]{0.23\linewidth}
\centering
\includegraphics[width=1.2in]{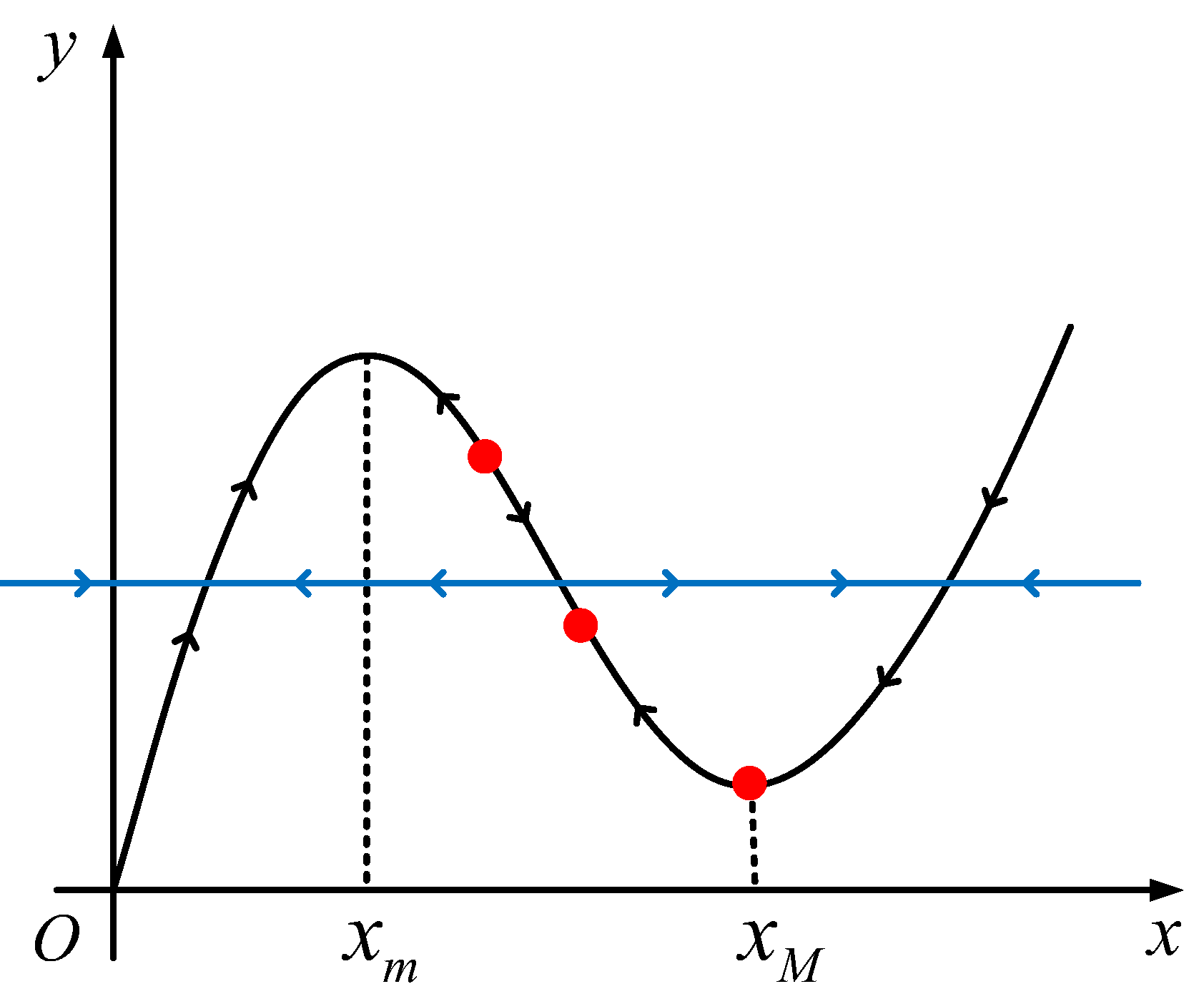}
\label{f-3-8}
\end{minipage}
}%
\subfigure[$MMR^{1}$.]{
\begin{minipage}[t]{0.23\linewidth}
\centering
\includegraphics[width=1.2in]{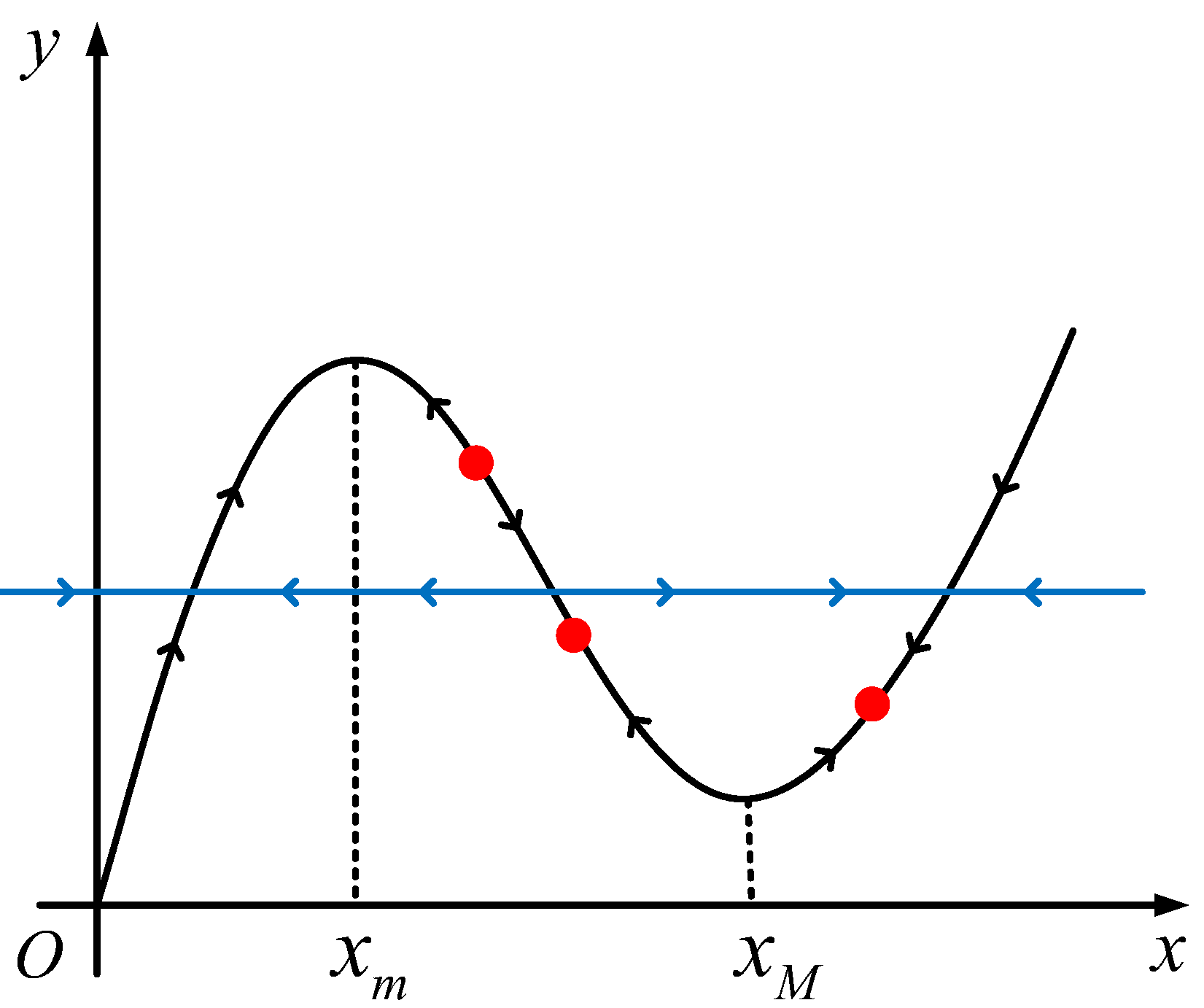}
\label{f-3-9}
\end{minipage}
}%
\centering
\caption{All possible intersection point sequences and the corresponding slow-fast limits.
Red dots are the equilibria lying on the graph of the function $\psi_{1}$ (black curve).
Black arrows indicate the flow of the reduced equation. Blue arrows indicate the flow of the layer equation.}
\label{fg-all-cases}
\end{figure}

\item{\bf (i)}
if the number of the intersection points is one,
then all possible intersection point sequences are  $L^{0}$, $L^{1}$, $M$, $R^{0}$ and $R^{1}$.

\item{\bf (ii)}
if the number of the intersection points is two,
then all possible intersection point sequences are $L^{0}M$, $L^{1}M$, $MM$, $MR^{0}$ and $MR^{1}$.

\item{\bf (iii)}
if the number of the intersection points is three,
then all possible intersection point sequences are
$L^{0}MR^{0}$, $L^{0}MR^{1}$, $L^{1}MR^{0}$, $L^{1}MR^{1}$, $L^{0}MM$, $L^{1}MM$, $MMM$, $MMR^{0}$ and $MMR^{1}$.
\end{lemma}

We give the lengthy proof for this lemma in Appendix A.

\section{Dynamics of the  high degradation rate case}
\label{sec-dynamics-1}
\setcounter{equation}{0}
\setcounter{lemma}{0}
\setcounter{theorem}{0}
\setcounter{remark}{0}

In this section,
we give the detailed study of the dynamics of the THTN model
in the  high degradation rate case,
that is, the rate of mRNA degradation is  high enough. Then $\varepsilon$ is sufficiently large.

\begin{lemma}
\label{thm-attraction}
Let the sets $\mathbb{R}^{2}_{+}$ and $\mathcal{A}$ be respectively defined by
$\mathbb{R}^{2}_{+}=\{(x,y)\in \mathbb{R}^{2}: x\geq 0, y \geq 0\}$ and
\begin{eqnarray*}
\mathcal{A}=\left\{(x,y)\in \mathbb{R}^{2}:
        0\leq x\leq \frac{v}{c}, \ \
        0\leq y \leq \frac{v}{c}\right\}.
\end{eqnarray*}
Then the sets $\mathbb{R}^{2}_{+}$ and $\mathcal{A}$ are both the positive invariant sets of system (\ref{2D-model-12}).
Furthermore, the set $\mathcal{A}$ attracts the set $\mathbb{R}^{2}_{+}$ under the flow of system (\ref{2D-model-12}).
\end{lemma}
{\bf Proof.}
By analyzing the field vector of system (\ref{2D-model-12})
along the boundaries of the sets $\mathbb{R}^{2}_{+}$ and $\mathcal{A}$,
the first statement can be obtained.
For each solution $(x(t),y(t))$ of system (\ref{2D-model-12})
with the initial value $(x(0),y(0))\in \mathbb{R}^{2}_{+}$,
we have that $x(t)\geq 0$ and $y(t)\geq 0$ for $t\geq 0$.
Then by the second equation in system (\ref{2D-model-12}),
we have that $y'(t)\leq -\varepsilon y +v/c$ for $t\geq 0$,
which together with Gronwall's Inequality yields that
\begin{eqnarray}
\label{est-y}
y(t)\leq y(0)e^{-\varepsilon t} +v/c, \  \ \  \ \ t\geq 0.
\end{eqnarray}
Consider the first equation in system (\ref{2D-model-12}) with $0\leq y(t)\leq v/c$.
Similarly, we have that
\begin{eqnarray}
\label{est-x}
x(t)\leq x(0)e^{-t} +v/c, \  \ \  \ \ t\geq 0.
\end{eqnarray}
Then by (\ref{est-y}) and (\ref{est-x}), the second statement holds.
Therefore, the proof is now complete.
\hfill$\Box$

For each finite equilibrium $(x_{0}, y_{0})$ of  system  (\ref{2D-model-12}) with $x_{0}\geq 0$,
in order to obtain the type of equilibrium $(x_{0}, y_{0})$,
we consider the the Jacobian matrix $\mathcal{J}(x_0,y_0)$ of system  (\ref{2D-model-12}) at $(x_{0}, y_{0})$
\begin{eqnarray*}
\mathcal{J}(x_0,y_0)
=\left(
\begin{array}{ll}
-\psi_{1}^{'}(x_0)  & 1
\\
\varepsilon \psi_{2}^{'}(x_0) & -\varepsilon
\end{array}
\right).
\end{eqnarray*}
The determinant and the trace of this Jacobian matrix are respectively given by
\begin{eqnarray}
\label{df-D-T}
\mathcal{D}(x_0,y_0):=\varepsilon (\psi_{1}^{'}(x_0)-\psi_{2}^{'}(x_0)), \ \ \ T(x_0,y_0):=-\varepsilon-\psi_{1}^{'}(x_0).
\end{eqnarray}
To determine the type of this equilibrium,
it is necessary to consider the constant
\begin{eqnarray}
\label{df-Delta}
\Delta (x_{0},y_{0}):=(T(x_0,y_0))^{2}-4\mathcal{D}(x_0,y_0)=(\varepsilon-\psi_{1}^{'}(x_0))^{2} +4\varepsilon \psi_{2}^{'}(x_0).
\end{eqnarray}
By the form of system  (\ref{2D-model-12}),
we observe that the value of $x_{0}$ is independent of the parameter $\varepsilon$
and only relies on the parameters $a$, $b_{i}$, $c$ and $v$.
Based on Bendixson's Theorem (see  \cite[Theorem 7.10, p. 188]{Dumortieretal06}),
we have the following statements.
\begin{theorem}
\label{prop-dynamics-1}
Consider system (\ref{2D-model-12}).
Then the following conclusions hold:

\item{\bf (i)}
if $\psi_{1}$ satisfies $\psi_{1}^{'}(x_{+})\geq 0$,
then there exists a unique equilibrium $(x_{0},y_{0})$ in $\mathbb{R}^{2}_{+}$,
which is a stable focus or node.
Furthermore, system (\ref{2D-model-12}) has no periodic orbits in $\mathbb{R}^{2}_{+}$,
and $(x_0,y_{0})$ attracts the set $\mathbb{R}^{2}_{+}$ under the flow of system (\ref{2D-model-12}).

\item{\bf (ii)}
if $\psi_{1}$ satisfies $-\varepsilon<\psi_{1}^{'}(x_{+})<0$,
then system (\ref{2D-model-12}) has no periodic orbits in $\mathbb{R}^{2}_{+}$,
and at least one equilibrium and at most three equilibria.
Further, the equilibria of system (\ref{2D-model-12}) admit the following trichotomies:
\begin{enumerate}
\item[{\bf (ii.1)}]
if system (\ref{2D-model-12}) has a unique equilibrium $(x_0,y_0)$,
then $(x_0,y_0)$ is a stable focus or node,
and $(x_0,y_{0})$ attracts the set $\mathbb{R}^{2}_{+}$ under the flow of system (\ref{2D-model-12}).

\item[{\bf (ii.2)}]
if system (\ref{2D-model-12}) has two equilibria $(x_0^{1},y_0^{1})$ and $(x_0^{2},y_0^{2})$,
then the point at which $\psi_{1}(x)=\psi_{2}(x)$ holds is a saddle-node,
the other point is a stable focus or node.

\item[{\bf (ii.3)}]
if system (\ref{2D-model-12}) has three equilibria $(x_0^{i},y_0^{i})$, $i=1,2,3$,
satisfying $x_0^{1}<x_0^{2}<x_0^{3}$,
then $(x_0^{1},y_0^{1})$ and $(x_0^{3},y_0^{3})$ are a stable focus or node,
and $(x_0^{2},y_0^{2})$ is a saddle.
\end{enumerate}
\end{theorem}
{\bf Proof.}
Under the condition $\psi_{1}^{'}(x_{+})\geq 0$,
Lemmas \ref{lm-psi-1-prpty} and  \ref{lm-psi-0-2-prpty} yield that
$\psi(0)=-v/c<0$, $\psi^{'}=\psi_{1}^{'}(x)-\psi_{2}^{'}(x)>0$ for $x>0$.
Then there is a unique equilibrium $(x_{0},y_{0})$ for system (\ref{2D-model-12}) in $\mathbb{R}^{2}_{+}$.
Further, this equilibrium satisfies $D(x_0,y_0)>0$ and $T(x_0,y_0)\leq -\varepsilon<0$,
which implies that $(x_0,y_{0})$ is a stable focus for $(\varepsilon-\psi_{1}^{'}(x_0))^{2} +4\varepsilon \psi_{2}^{'}(x_0)< 0$
and is a stable node for $(\varepsilon-\psi_{1}^{'}(x_0))^{2} +4\varepsilon \psi_{2}^{'}(x_0)\geq 0$.
Assume that $\psi_{1}$ satisfies $\psi_{1}^{'}(x_{+})\geq 0$.
Then by Lemma \ref{lm-psi-1-prpty},
\begin{eqnarray}\label{Bdixson-1}
\frac{\partial}{\partial x}(y-\psi_{1}(x))+\frac{\partial}{\partial y}(\varepsilon(\psi_{2}(x)-y)))
=-(\varepsilon+\psi_{1}^{'}(x))\leq -\varepsilon, \ \ \ x\geq 0.
\end{eqnarray}
Hence, Bendixson's Theorem yields that
system (\ref{2D-model-12}) has no periodic orbits in $\mathbb{R}^{2}_{+}$.
Recall that $(x_0,y_{0})$ is a stable focus or node,
then $(x_0,y_{0})$ attracts the set $\mathbb{R}^{2}_{+}$ under the flow of system (\ref{2D-model-12}).
Thus, the statements in {\bf (i)} are proved.

If $\psi_{1}$ satisfies $-\varepsilon<\psi_{1}^{'}(x_{+})<0$,
then by similar method used in the proof for {\bf (i)},
we obtain that system (\ref{2D-model-12}) has no periodic orbits in $\mathbb{R}^{2}_{+}$.
As for the types of equilibria, we only give the proof for the case {\bf (ii.2)}.
Without loss of generality, assume that $\psi_{1}(x_{0}^{1})=\psi_{2}(x_{0}^{1})$ and $x_{0}^{1}>x_{0}^{2}$.
Then by Lemmas \ref{lm-psi-1-prpty} and \ref{lm-psi-0-2-prpty},
we can obtain that $T(x_0^{i},y_0^{i})<0$, $\mathcal{D}(x_0^{1},y_0^{1})=0$, $\mathcal{D}(x_0^{2},y_0^{2})>0$
and $\varepsilon (\psi_{1}^{''}(x_0)-\psi_{2}^{''}(x_0))<0$.
Hence, $(x_0^{2},y_0^{2})$ is a stable focus or node,
and by using \cite[Theorem 7.1, p.114]{ZZF-etal} (see also the proof in Theorem \ref{thm-connection}),
we obtain that $(x_0^{1},y_0^{1})$ is a saddle-node.
Therefore, the proof is now complete.
\hfill$\Box$

\begin{remark}
Whether  an equilibrium is a focus or node,
is determined by the sign of
$\Delta (x_{0},y_{0})=(\varepsilon-\psi_{1}^{'}(x_0))^{2} +4\varepsilon \psi_{2}^{'}(x_0)$ (see \cite{Dumortieretal06,ZZF-etal}).
More precisely,
if $\Delta (x_{0},y_{0})=(\varepsilon-\psi_{1}^{'}(x_0))^{2} +4\varepsilon \psi_{2}^{'}(x_0)<0$ (resp. $\geq 0$),
then it is a focus (resp. node).
We also remark that  for sufficiently large $\varepsilon=k_{m}/k_{3}>|\psi_{1}^{'}(x_{+})|$,
there are no periodic orbits in system (\ref{2D-model-12}).
\end{remark}

\section{Dynamics of the  low degradation rate case}
\label{sec-dynamics-2}
\setcounter{equation}{0}
\setcounter{lemma}{0}
\setcounter{theorem}{0}
\setcounter{remark}{0}

In this section,
we consider the dynamics of the THTN model in the  low degradation rate case,
that is, the rate of mRNA degradation  is  low enough.
Throughout this section, we always assume that $0<\varepsilon\ll1$ and $v$ is independent of $\varepsilon$.

Under the condition that the parameter $\varepsilon $ is sufficiently small,
system (\ref{2D-model-12}) is a standard slow-fast system of the form (\ref{fast-1}).
For convenience, here we write $\psi_{1}(x,\lambda)$ and $\psi_{2}(x,\lambda,v)$,  instead of  $\psi_{1}(x)$ and  $\psi_{2}(x)$,
where $\lambda=(a,b_{1},b_{2},c)$,
then system (\ref{2D-model-12}) can be written as
\begin{eqnarray}
\label{2D-model-14}
\begin{aligned}
\frac{d x}{d t} &=x'=  y-\psi_{1}(x,\lambda):=f(x,y,\lambda),
\\
\frac{d y}{d t} &=y'= \varepsilon\left(\psi_{2}(x,\lambda,v)-y\right):=\varepsilon g(x,y,\lambda,v).
\end{aligned}
\end{eqnarray}
By a time rescaling $s=\varepsilon t$,
the slow system corresponding to system (\ref{2D-model-14}) is in the form
\begin{eqnarray}
\label{2D-slow-1}
\begin{aligned}
\varepsilon \frac{d x}{d s} &= \varepsilon \dot x=y-\psi_{1}(x,\lambda),
\\
\frac{d y}{d s} &=\dot y=\psi_{2}(x,\lambda,v)-y.
\end{aligned}
\end{eqnarray}
Let the set $\mathcal{C}_{0}$ be defined by
$\mathcal{C}_{0}=\left\{(x,y)\in \mathbb{R}\times \mathbb{R}: y=\psi_{1}(x,\lambda)\right\}.$
Throughout this section we always assume that $\psi_{1}$ satisfies $\psi_{1}^{'}(x_{+})<0$ for suitable parameters $\lambda$ and $v$.
Then the set $\mathcal{C}_{0}$ is $S$-shaped.
Due to Lemma \ref{lm-psi-1-prpty}, all points in the set $\mathcal{C}_{0}$,
except $(x_{i},y_{i}):=(x_{i},\psi_{1}(x_{i}))$, $i=m,M$, are normally hyperbolic.
Then by the theory of normally hyperbolic invariant manifolds,
the reduced system on $L^{1}\cup M\cup R^{1}$ is governed by
\begin{eqnarray}
\label{system-redu}
\frac{\partial \psi_{1}}{\partial x}(x,\lambda)\frac{d x}{d s}
   =\psi_{2}(x,\lambda,v)- \psi_{1}(x,\lambda).
\end{eqnarray}
In the following,  we investigate the dynamics of the THTN model in the  low degradation rate
case by employing geometric singular perturbation theory.

\subsection{Local dynamics of canard points}

In this section we study the local dynamics of canard points.
Assume that for $\lambda=\lambda^{0}$ and $v=v^{0}$,
either $(x_{m},y_{m})$ or $(x_{M},y_{M})$ is an equilibrium of the slow-fast system (\ref{2D-model-14}).
Then at this point $(x_{i},y_{i})$, $i=m$ or $M$, we have that
$f(x_{i},y_{i},\lambda^{0})=0$ and $g(x_{i},y_{i},\lambda^{0},v^{0})=0$.
By Lemma \ref{lm-psi-1-prpty}
the function $f$ satisfies
$$
\frac{\partial f}{\partial x}(x_{i},y_{i},\lambda^{0})
=-\frac{\partial \psi_{1}}{\partial x}(x_{i},\lambda^{0})=0,
$$
which yields that the critical manifold $\mathcal{C}_{0}$ loses hyperbolicity at $(x_{i},y_{i})$
and $(x_{i},y_{i})$ is a contact point.
Further,
following Lemmas \ref{lm-psi-1-prpty} and \ref{lm-psi-0-2-prpty},
the slow-fast system (\ref{2D-model-14}) satisfies the  nondegeneracy conditions:
\begin{eqnarray*}
&& \frac{\partial^{2}f}{\partial x^{2}}(x_{i},y_{i},\lambda^{0})
=-\frac{\partial^{2} \psi_{1}}{\partial x^{2}}(x_{i},\lambda^{0})\neq 0,\ \ \ \ \
\frac{\partial f}{\partial y}(x_{i},y_{i},\lambda^{0})=1,\\
&& \frac{\partial g}{\partial x}(x_{i},y_{i},\lambda^{0},v^{0})
=\frac{\partial \psi_{2}}{\partial x}(x_{i},\lambda^{0},v^{0})<0,\ \ \ \ \
\frac{\partial g}{\partial v}(x_{i},y_{i},\lambda^{0},v^{0})=\frac{1}{c^{0}+(x_{i}-\phi(x_{i}))^{2}}>0,
\end{eqnarray*}
where $\frac{\partial^{2} \psi_{1}}{\partial x^{2}}(x_{i},\lambda^{0})<0$ for $i=m$
and $\frac{\partial^{2} \psi_{1}}{\partial x^{2}}(x_{i},\lambda^{0})>0$ for $i=M$.
Then  by (3.2), (3.3) and (3.4) in \cite[p.303]{Krupa-Szmolyan-01SIMA},
the above nondegeneracy conditions insure that the contact point $(x_{i},y_{i})$
is a canard point of the slow-fast system (\ref{2D-model-14}).

We next consider the normal forms of system (\ref{2D-model-14}) near the canard points $(x_{i},y_{i})$, $i=m,M$.
\begin{lemma}
\label{lm-normal-form-1}
Assume that for $\lambda=\lambda^{0}$ and $v=v^{0}$,
either $(x_{m},y_{m})$ or $(x_{M},y_{M})$ is an equilibrium of the slow-fast system (\ref{2D-model-14}).
Then for fixed $\lambda=\lambda^{0}$,
the slow-fast system (\ref{2D-model-14}) near $(x_{m},y_{m})$ and $(x_{M},y_{M})$ can be changed into
\begin{eqnarray}
\label{Normal-form-2D}
\begin{aligned}
x' &=  -y+x^{2}\Phi_{1}(x),
\\
y'&=\varepsilon\left(x\Phi_{2}(x,v)
    -v+\frac{1}{D_{1}\psi_{2}(x_{i},\lambda^{0},v^{0})}y\right),
\end{aligned}
\end{eqnarray}
where $\Phi_{j}$ are defined by
\begin{eqnarray*}
\Phi_{1}(x)\!\!\!&=&\!\!\! 1+\frac{2}{\varphi_{1}^{''}(0)}\widehat{\Phi}_{1}(-\frac{2}{\varphi_{1}^{''}(0)}x),\\
\Phi_{2}(x,v)\!\!\!&=&\!\!\! 1+\frac{1}{D_{1}\varphi_{2}(0,0)}
    \widehat{\Phi}_{2}
    \left(-\frac{2}{\varphi_{1}^{''}(0)}x,\frac{2D_{1}\varphi_{2}(0,0)(c^{0}+(x_{i}-\phi(x_{i}))^{2})}{\varphi_{1}^{''}(0)}v\right),
\end{eqnarray*}
and  the functions $\varphi_{j}$ and $ \widehat{\Phi}_{j}$ are in the form
\begin{eqnarray}
 \varphi_{1}(x)\!\!\!&=&\!\!\!\psi_{1}(x+x_{i},\lambda^{0})-y_{i},
     \ \ \ \varphi_{2}(x,v)=\psi_{2}(x+x_{i},\lambda^{0},v+v^{0})-y_{i},\label{df-varphi-i}\\
\widehat{\Phi}_{1}(x)\!\!\!&=&\!\!\! \int_{0}^{1} \int_{0}^{1}\alpha\varphi_{1}^{''}(\alpha\beta x) d\alpha d\beta-\frac{1}{2}\varphi_{1}^{''}(0),
 \nonumber\\
\widehat{\Phi}_{2}(x,v)\!\!\!&=&\!\!\! x\int_{0}^{1} \int_{0}^{1}\alpha D_{11}\varphi_{2}(\alpha\beta x,0) d\alpha d\beta
                +v\int_{0}^{1} \int_{0}^{1} D_{12}\varphi_{2}(\alpha x,\beta v) d\alpha d\beta\nonumber.
\end{eqnarray}
Here, $D_{ij}=D_{j}\circ D_{i}$ and the operator $D_{j}$ denotes the partial derivative with respect to the {\it j-}th variable.
\end{lemma}
{\bf Proof.}
Assume that $(x_{i},y_{i})$, $i=m$ or $M$, is an equilibrium of
system (\ref{2D-model-14}) with $\lambda=\lambda^{0}$ and $v=v^{0}$.
Let  $\lambda=\lambda^{0}$ be fixed.
Then by a translation transformation $\mathcal{T}_{1}$
of the form
\begin{eqnarray}
\label{df-T-1}
\mathcal{T}_{1}: (x,y,v)\to (x+x_{i},y+y_{i},v+v^{0}),
\end{eqnarray}
system (\ref{2D-model-14}) is transformed into the form
\begin{eqnarray}
\label{2D-model-15}
\begin{aligned}
x' &=y-\varphi_{1}(x),
\\
y'&=\varepsilon\left(\varphi_{2}(x,v)-y\right),
\end{aligned}
\end{eqnarray}
where $\varphi_{i}$ are defined by (\ref{df-varphi-i}) satisfying $\varphi_{1}(0)=0$, $\varphi_{1}^{'}(0)=0$ and $\varphi_{2}(0,0)=0$.
Thus the function $\varphi_{1}$ can be written as the form
\begin{eqnarray*}
\varphi_{1}(x)
              =x\int_{0}^{1} \varphi_{1}^{'}(\alpha x) d\alpha
              =x^{2}\int_{0}^{1} \int_{0}^{1}\alpha\varphi_{1}^{''}(\alpha\beta x) d\alpha d\beta,
\end{eqnarray*}
which implies
\begin{eqnarray*}
\varphi_{1}(x)=x^{2}
\left(\frac{1}{2}\varphi_{1}^{''}(0)+\widehat{\Phi}_{1}(x)\right).
\end{eqnarray*}
Similarly, we have
\begin{eqnarray*}
\varphi_{2}(x,v)\!\!\!&=&\!\!\! \varphi_{2}(x,v)-\varphi_{2}(0,v)+\varphi_{2}(0,v)\\
                \!\!\!&=&\!\!\!x\int_{0}^{1} D_{1}\varphi_{2}(\alpha x,v) d\alpha +\frac{v}{c^{0}+(x_{i}-\phi(x_{i}))^{2}}\\
                \!\!\!&=&\!\!\!
                x\left(D_{1}\varphi_{2}(0,0)+x\int_{0}^{1} \int_{0}^{1}\alpha D_{11}\varphi_{2}(\alpha\beta x,0) d\alpha d\beta
                +v\int_{0}^{1} \int_{0}^{1} D_{12}\varphi_{2}(\alpha x,\beta v) d\alpha d\beta\right)\\
                \!\!\!& &\!\!\!
                +\frac{v}{c^{0}+(x_{i}-\phi(x_{i}))^{2}}\\
                \!\!\!&=&\!\!\!
                x\left(D_{1}\varphi_{2}(0,0)+\widehat{\Phi}_{2}(x,v)\right)+\frac{v}{c^{0}+(x_{i}-\phi(x_{i}))^{2}}.
\end{eqnarray*}
By taking a coordinate transformation $\mathcal{T}_{2}$ of the form
\begin{eqnarray}
\label{df-T-2}
\mathcal{T}_{2}: (x,y,v,\varepsilon)
       \to \left(-\frac{2}{\varphi_{1}^{''}(0)}x,\
              \frac{2}{\varphi_{1}^{''}(0)}y,\
              \frac{2D_{1}\varphi_{2}(0,0)(c^{0}+(x_{i}-\phi(x_{i}))^{2})}{\varphi_{1}^{''}(0)}v,\
              -\frac{1}{D_{1}\varphi_{2}(0,0)}\varepsilon
               \right),
\end{eqnarray}
system (\ref{2D-model-15}) is changed into the form (\ref{Normal-form-2D}).
Therefore, the proof is now complete.
\hfill$\Box$

Next we define several constants,
which play important roles in the analysis of the dynamics near the canard points.
Similarly to the formulae (3.12) and (3.13) in \cite{Krupa-Szmolyan-01JDE},
let
\begin{eqnarray*}
\kappa_{i,1}=\frac{d\Phi_{1}}{dx}(0), \ \  \
\kappa_{i,2}=\frac{\partial \Phi_{2}}{\partial x}(0,0),\ \ \
\kappa_{i,3}=\frac{1}{D_{1}\psi_{2}(x_{i},\lambda^{0},v^{0})},  \ \ \ i=m,\,M,
\end{eqnarray*}
and define $A_{i}$ by
\begin{eqnarray*}
A_{i}=3\kappa_{i,1}-2\kappa_{i,2}-2\kappa_{i,3},\ \ \ i=m,\,M.
\end{eqnarray*}
Here the key constants $A_{i}$ determine the nondegeneracy conditions
for the Hopf bifurcations near the canard points $(x_{i},y_{i})$
and are greatly important for the analysis of canard explosions
(See \cite{Krupa-Szmolyan-01JDE,Kuehn-15}).
By a direct computation we obtain
\begin{eqnarray}
&&\kappa_{i,1}=-\frac{2D_{111}\psi_{1}(x_{i},\lambda^{0})}{3(D_{11}\psi_{1}(x_{i},\lambda^{0}))^{2}},\ \
\kappa_{i,2}=-\frac{D_{11}\psi_{2}(x_{i},\lambda^{0},v^{0})}{D_{11}\psi_{1}(x_{i},\lambda^{0})D_{1}\psi_{2}(x_{i},\lambda^{0},v^{0})},\ \
\kappa_{i,3}=\frac{1}{D_{1}\psi_{2}(x_{i},\lambda^{0},v^{0})},\nonumber\\
&&A_{i}=-\frac{2 D_{111}\psi_{1}(x_{i},\lambda^{0})}{(D_{11}\psi_{1}(x_{i},\lambda^{0}))^{2}}
    +\frac{2D_{11}\psi_{2}(x_{i},\lambda^{0},v^{0})}{D_{11}\psi_{1}(x_{i},\lambda^{0})D_{1}\psi_{2}(x_{i},\lambda^{0},v^{0})}
    -\frac{2}{D_{1}\psi_{2}(x_{i},\lambda^{0},v^{0})}.\label{df-A-i}
\end{eqnarray}
Compared the above notations to the corresponding ones  in \cite{Krupa-Szmolyan-01SIMA},
the functions $h_{j}$ in \cite[system (3.6), p.304]{Krupa-Szmolyan-01SIMA}
are in the form
\begin{eqnarray*}
h_{1}=1, \ \ \ h_{2}=\Phi_{1},\ \ \ h_{3}=0,
\ \ \ h_{4}=\Phi_{2}, \ \ \ h_{5}=1,
\ \ \ h_{6}=\frac{1}{D_{1}\psi_{2}(x_{i},\lambda^{0},v^{0})},
\end{eqnarray*}
and the constants $a_{j}$ introduced in \cite[p.305]{Krupa-Szmolyan-01SIMA}
are in the form
\begin{eqnarray*}
a_{1}=a_{2}=0,\ \ \
a_{3}=\kappa_{i,1},\ \ \
a_{4}=\kappa_{i,2},\ \ \
a_{5}=\kappa_{i,3}.
\end{eqnarray*}
Since the constants $A_{i}$ satisfy
\begin{eqnarray*}
\lefteqn{-\frac{1}{2}(D_{11}\psi_{1}(x_{i},\lambda^{0}))^{2}\cdot D_{1}\psi_{2}(x_{i},\lambda^{0},v^{0})\cdot A_{i}}\\
\!\!\!&=&\!\!\!
  D_{111}\psi_{1}(x_{i},\lambda^{0})\cdot D_{1}\psi_{2}(x_{i},\lambda^{0},v^{0})
       -D_{11}\psi_{1}(x_{i},\lambda^{0}) \cdot D_{11}\psi_{2}(x_{i},\lambda^{0},v^{0})
       +(D_{11}\psi_{1}(x_{i},\lambda^{0}))^{2},
\end{eqnarray*}
then by a direct computation,
three different cases  $A_{i}<0$, $A_{i}>0$ and $A_{i}=0$
can appear under some suitable conditions.
We follow \cite{Krupa-Szmolyan-01JDE} and analyze the canard explosion in (\ref{2D-model-14}).
Thus, we assume that $A_{i}\neq 0$ for $i=m,M$.
This implies that two Hopf bifurcations near $(x_{i},y_{i})$ are both nondegerate
(see {\bf (iv)} of Theorem \ref{thm-one}).
The cases $A_{i}=0$ will be studied in the future.

\begin{figure}[!htbp]
\centering
\subfigure[]{
\begin{minipage}[t]{0.45\linewidth}
\centering
\includegraphics[width=2.2in]{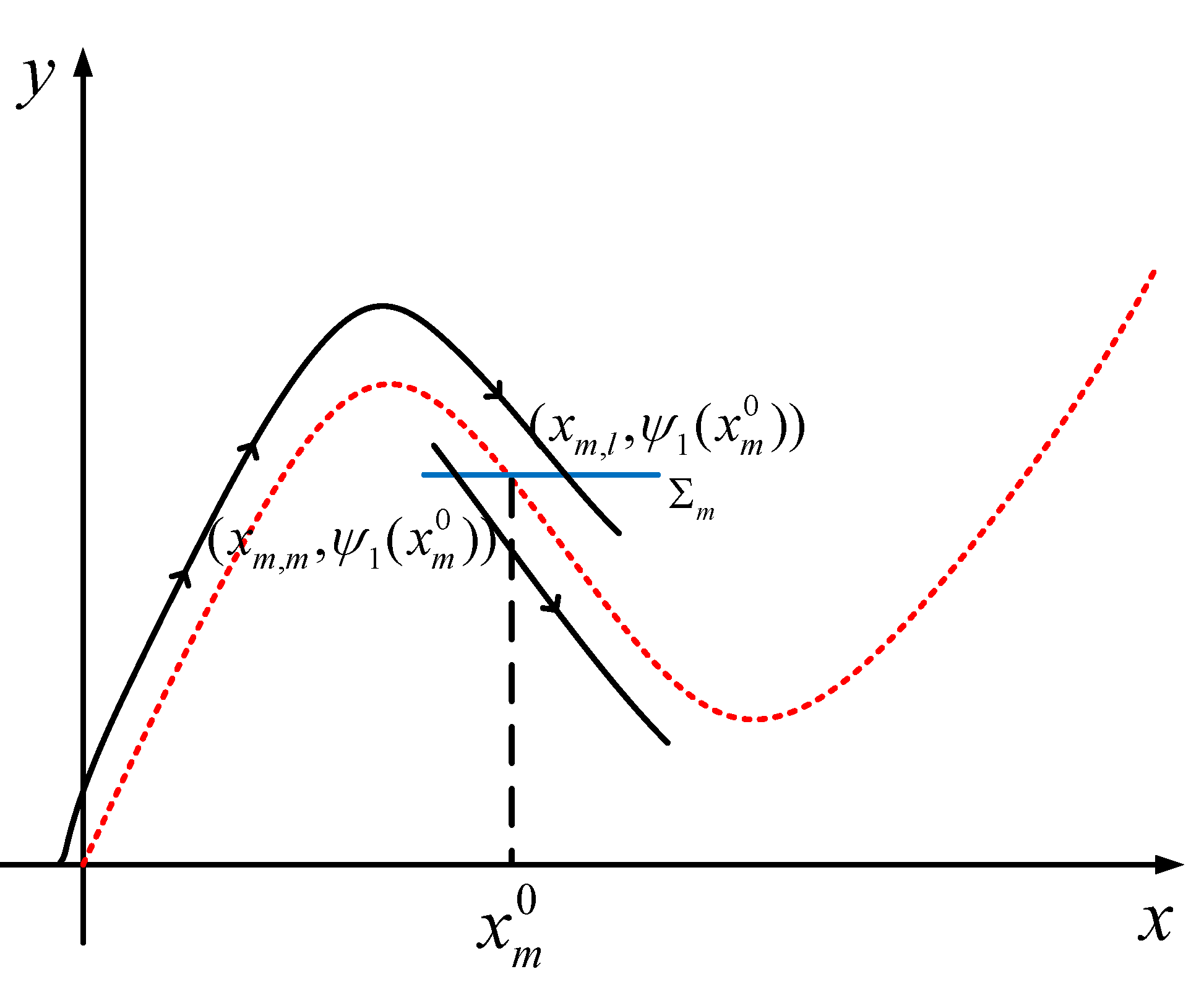}
\label{fig-canard-m}
\end{minipage}%
}%
\subfigure[]{
\begin{minipage}[t]{0.45\linewidth}
\centering
\includegraphics[width=2.2in]{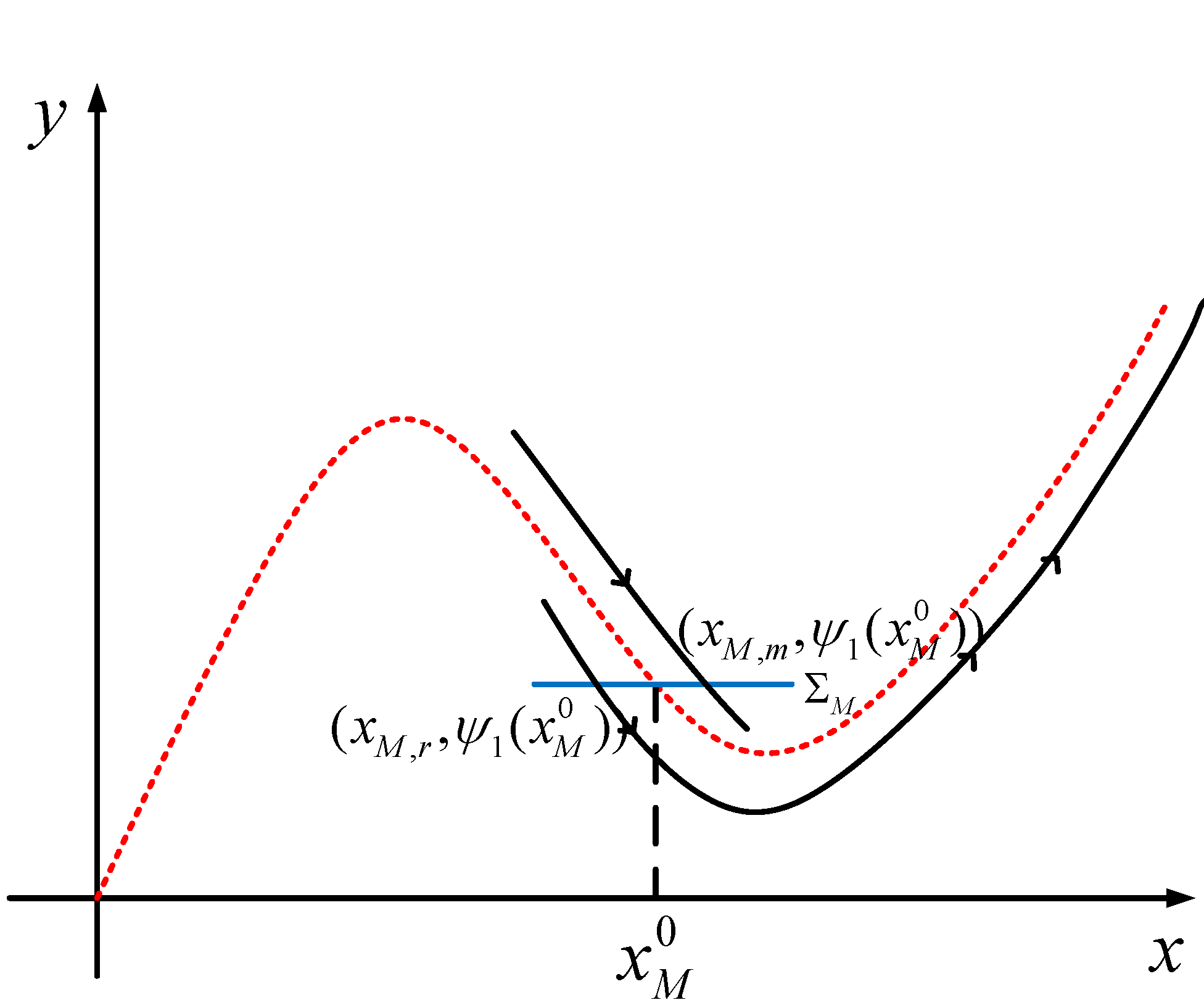}
\label{fig-canard-M}
\end{minipage}
}%
\centering
\caption{Dynamics of the slow-fast system (\ref{2D-model-14}) near the canard points $(x_{m},y_{m})$ and $(x_{M},y_{M})$.
The black curves are the orbits of system (\ref{2D-model-14}).
The dashed red curve is the graph of function $\psi_{1}$.}
         \label{Fg-near-canard}
\end{figure}

For sufficiently small $\varepsilon>0$,
one can see that the manifold $L^{1}$, $M$ and $R^{1}$ perturb smoothly to locally invariant manifolds $L^{1}_{\varepsilon}$,
$M_{\varepsilon}$ and $R^{1}_{\varepsilon}$, respectively.
Assume that $(x_{m},y_{m})$ (resp.  $(x_{M},y_{M})$) is a canard point.
Let $\Sigma_{m}$ (resp. $\Sigma_{M}$) be the cross-section of the curve $M$
at the point $(x_{m}^{0},\psi_{1}(x_{m}^{0}))$ (resp. $(x_{M}^{0},\psi_{1}(x_{M}^{0}))$ along the $x$-direction,
where $x_{m}^{0}$ (resp. $x_{M}^{0}$) satisfies that $x_{m}^{0}-x_{m}$ (resp. $x_{M}-x_{M}^{0}$)
is positive and sufficiently small.
Let the manifold $L^{1}_{\varepsilon}$ (resp. $R^{1}_{\varepsilon}$) and $M_{\varepsilon}$   extend in the neighborhood of  this canard point.
Assume that they respectively intersect with the section $\Sigma_{m}$ (resp. $\Sigma_{M}$) at points $(x_{m,l},\psi_{1}(x_{m}^{0}))$
and $(x_{m,m},\psi_{1}(x_{m}^{0}))$
(resp. $(x_{M,m}, \psi_{1}(x_{M}^{0}))$ and $(x_{M,r}, \psi_{1}(x_{M}^{0}))$).
See Figure \ref{Fg-near-canard}.
%
We have the following.

\begin{lemma}
\label{thm-connection}
Assume that  for $\lambda=\lambda^{0}$ and $v=v^{0}$,
the slow-fast system (\ref{2D-model-14}) has an equilibrium at either $(x_{m},y_{m})$  or $(x_{M},y_{M})$ for $x\geq 0$.
Then for sufficiently small $\varepsilon>0$,
there exist two smooth functions $v^{c}_{i}$, $i=m,M$, defined by
\begin{eqnarray}
\label{df-canard-curve}
v^{c}_{i}(\varepsilon)=v^{0}+
\mathcal{K}_{i}\varepsilon+O(\varepsilon^{3/2}), \ \ \ i=m,M,
\end{eqnarray}
such that the slow-fast system (\ref{2D-model-14}) with $\lambda=\lambda^{0}$ has $x_{m,l}=x_{m,m}$ for $i=m$ and $x_{M,m}=x_{M,r}$ for $i=M$
if and only if  $v=v^{c}_{i}(\varepsilon)$,
where the constants $\mathcal{K}_{i}$ are defined by
\begin{eqnarray}
\mathcal{K}_{i}=(\kappa_{i,3}+\frac{A_{i}}{4})\cdot
\frac{(D_{1}\psi_{2}(x_{i},\lambda^{0},v^{0}))^{2}(c^{0}+(x_{i}-\phi(x_{i}))^{2})}{D_{11}\psi_{1}(x_{i},\lambda^{0})}, \ \ \  i=m,M.
\label{df-K-i}
\end{eqnarray}
Furthermore, if $(x_{m},y_{m})$ is a canard point,
then $x_{m,l}>x_{m,m}$ for $0<v-v^{c}_{m}(\varepsilon)\ll 1$
and $x_{m,l}<x_{m,m}$ for $0<v^{c}_{m}(\varepsilon)-v\ll 1$.
If $(x_{M},y_{M})$ is a canard point,
then  $x_{M,m}>x_{M,r}$ for $0<v^{c}_{M}(\varepsilon)-v\ll 1$ and $x_{M,m}<x_{M,r}$ for $0<v-v^{c}_{M}(\varepsilon)\ll 1$.
\end{lemma}
{\bf Proof.}
We only give the proof for the case $(x_{m},y_{m})$.
Under the transformation $\mathcal{T}_{2}\circ \mathcal{T}_{1}$,
we assume that the points $(x_{m,l},\psi_{1}(x_{m}^{0}))$ and $(x_{m,m},\psi_{1}(x_{m}^{0}))$
are changed to the points $(w_{m,l}, z_{m})$ and $(w_{m,m}, z_{m})$, respectively.
Recall that the transformations $\mathcal{T}_{j}$, $j=1,2$, are given by (\ref{df-T-1}) and (\ref{df-T-2}),
and $\varphi_{1}^{''}(0)=D_{11}\psi_{1}(x_{m},\lambda^{0},v^{0})<0$,
then $x_{m,l}-x_{m,m}$ and $w_{m,l}-w_{m,m}$ have the same sign.
To finish the proof for this lemma,
we consider the normal form (\ref{Normal-form-2D}) of system (\ref{2D-model-14}) near $(x_{m},y_{m})$.
By \cite[Theorem 3.1]{Krupa-Szmolyan-01SIMA}
there exists a smooth function $\widehat{v}^{c}_{m}(\cdot)$ defined by
\begin{eqnarray*}
\widehat{v}^{c}_{m}(\varepsilon)=-\frac{4\kappa_{m,3}+A_{m}}{8}\varepsilon+O(\varepsilon^{3/2})
\end{eqnarray*}
such that system (\ref{Normal-form-2D}) has $w_{m,l}=w_{m,m}$ if and only if  $v=\widehat{v}^{c}_{m}(\varepsilon)$ .
Thus, by taking the variable transformation $\mathcal{T}_{1}^{-1}\circ \mathcal{T}_{2}^{-1}$
we obtain that (\ref{df-canard-curve}) holds for $i=m$.
Since the constant $d_{\lambda_{2}}$ in \cite[formula (3.23)]{Krupa-Szmolyan-01SIMA} is negative,
then the remaining statements hold.
Thus, the proof is finished.
\hfill$\Box$

\subsection{Global dynamics of the slow-fast system (\ref{2D-model-14})}
In this section,
we study the global dynamics of the slow-fast system (\ref{2D-model-14}).
The discussion is divided into three different parts according to the number of equilibria.

\subsubsection{One equilibrium}

Assume that the slow-fast system (\ref{2D-model-14}) with $\lambda=\lambda^{0}$ and $v=v^{0}$
has exactly one equilibrium $(x_{0},y_{0})$ in the set $x\geq 0$.
Then all types of the intersection point sequences are
$L^{1}$, $L^{0}$, $M$, $R^{0}$ and $R^{1}$.
See Figures \ref{f-1-1}, \ref{f-1-2}, \ref{f-1-3}, \ref{f-1-4} and \ref{f-1-5}.

\begin{figure}[htbp]
\centering
\subfigure[]{
\label{fig-relaxtion}
\includegraphics[height=4cm,width=4.8cm]{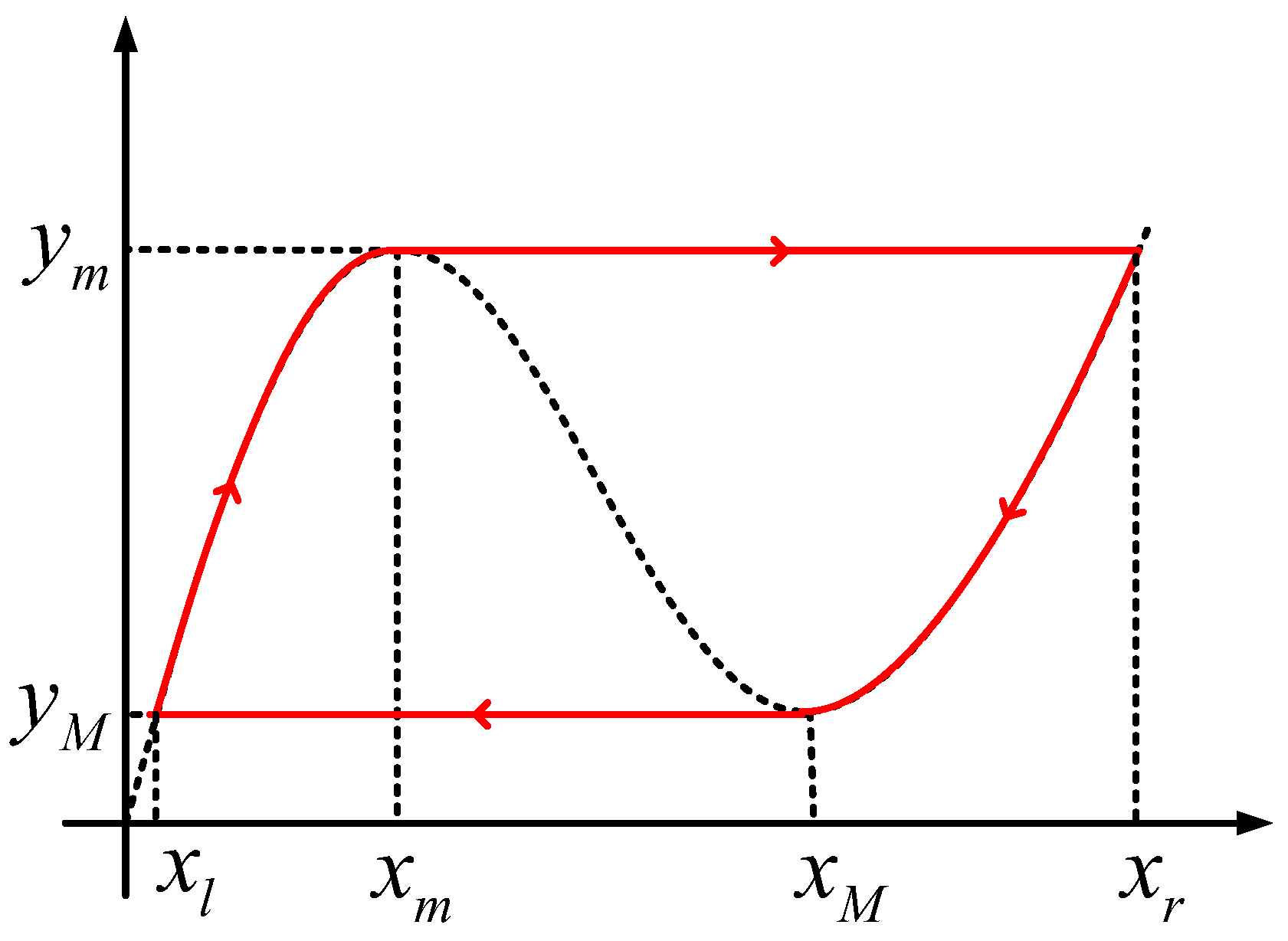}}
\subfigure[]{
\label{fig-canard-nohead}
\includegraphics[height=4cm,width=4.8cm]{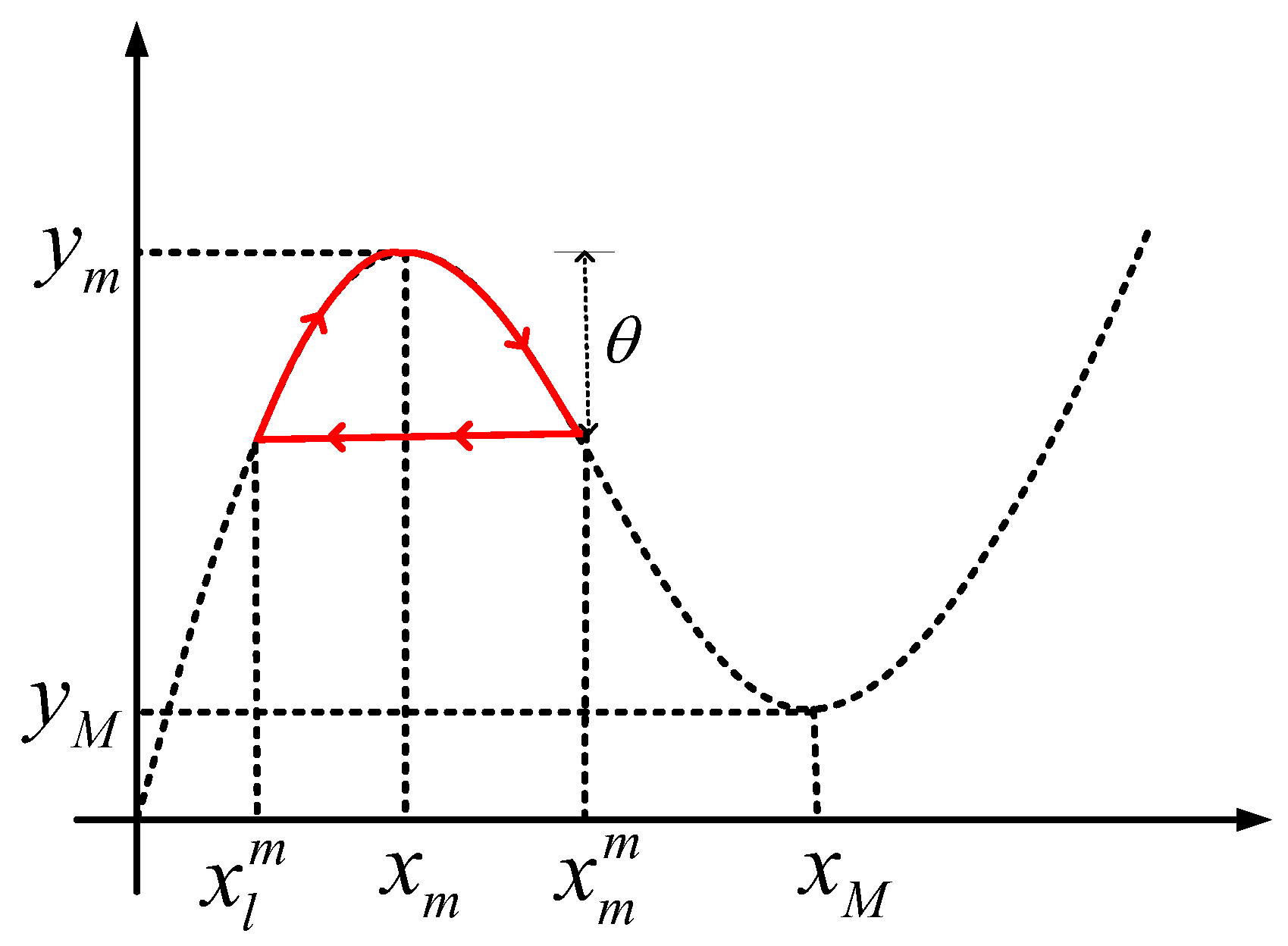}}
\subfigure[]{
\label{fig-canard-head}
\includegraphics[height=4cm,width=4.8cm]{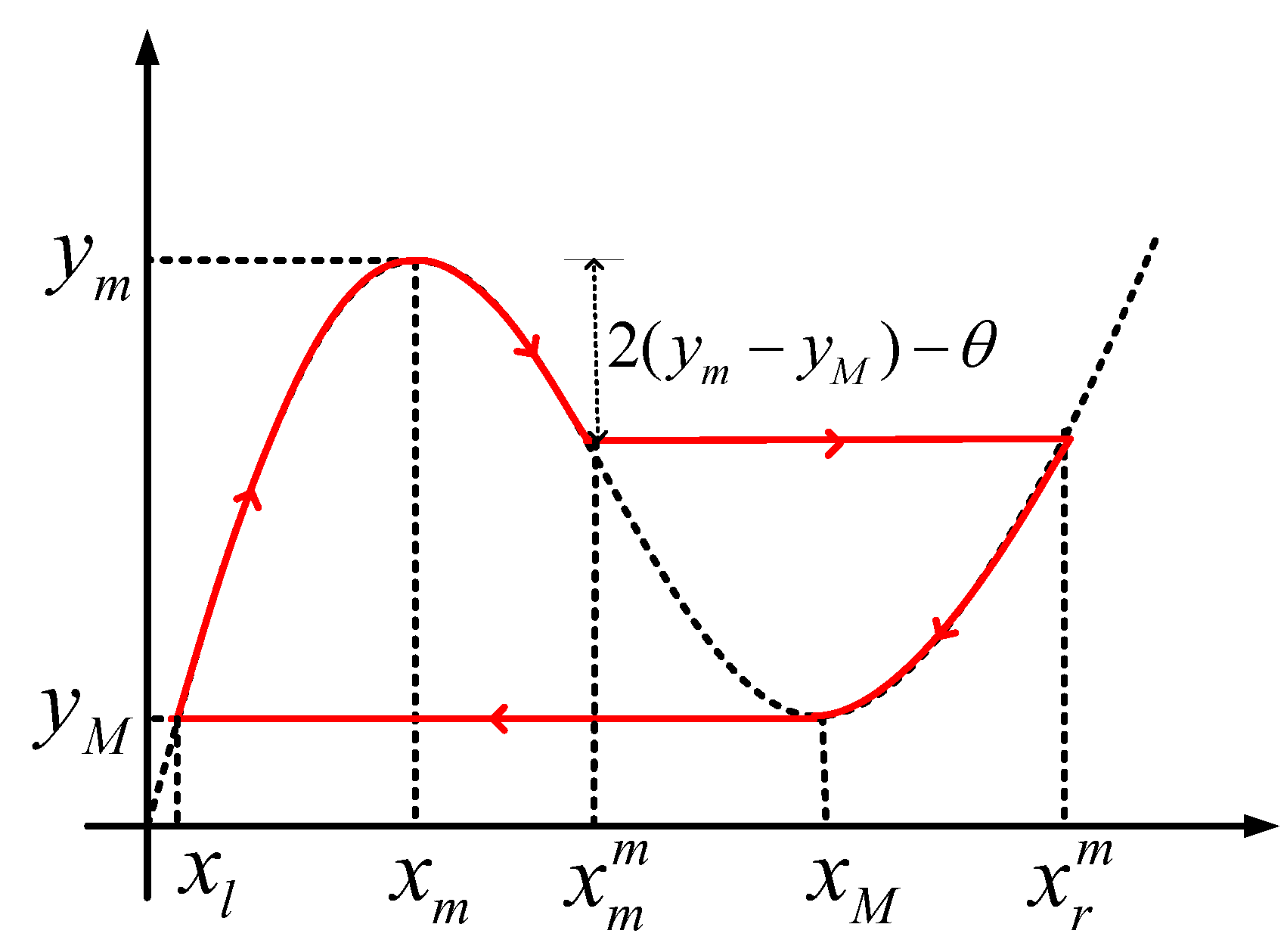}}
\caption{Slow-fast cycles (the red curves) are constructed:
        \ref{fig-relaxtion} Singular relaxation cycle. \
        \ref{fig-canard-nohead} Canard slow-fast cycle without head.\
         \ref{fig-canard-head} Canard slow-fast cycle with head.
         }
         \label{fig-singular-cycle}
\end{figure}

If the unique equilibrium $(x_{0},y_{0})$ is of type $M$,
then $(x_{i},y_{i})$ are both jump points.
Let $x_{l}$ (resp. $x_{r}$) be the value such that $\psi_{1}(x_{l},\lambda^{0})=y_{M}$
(resp. $\psi_{1}(x_{r},\lambda^{0})=y_{m}$)
and $(x_{l},y_{M})\in L$ (resp. $(x_{r},y_{m})\in R$).
We define a singular relaxation cycle $\Gamma_{r}$.
See Figure \ref{fig-relaxtion}.
This cycle $\Gamma_{r}$ consists of four branches,
among which two branches are the critical fibers of the layer equation
joining $(x_{m},y_{m})$ to $(x_{r},y_{m})$ and $(x_{M},y_{M})$ to $(x_{l},y_{M})$,
another two branches are the parts of the critical manifolds
joining $(x_{l},y_{M})$ to $(x_{m},y_{m})$ and $(x_{r},y_{m})$ to $(x_{M},y_{M})$.

If  the unique equilibrium $(x_{0},y_{0})$ is of type $L^{0}$ or type $R^{0}$,
then $(x_{m},y_{m})$ or $(x_{M},y_{M})$  is a canard point.
As a preparation, we next begin with the construction of canard slow-fast  cycles.
See Figures \ref{fig-canard-nohead} and \ref{fig-canard-head}.
For a positive constant $\theta$ with $0<\theta<y_{m}-y_{M}$,
let the constants $x^{m}_{j}$, $j=l,m,r$, with $0<x^{m}_{l}(\theta)<x_{m}<x^{m}_{m}(\theta)<x_{M}<x^{m}_{r}(\theta)$,
denote the roots of equation $\psi_{1}(x,\lambda^{0})=y_{m}-\theta$.
We define the canard slow-fast  cycles $\Gamma_{m}(\theta)$, $0\leq \theta \leq 2(y_{m}-y_{M})$,
for the canard point $(x_{m},y_{m})$ as follows.
For $0\leq \theta\leq y_{m}-y_{M}$,
$$\Gamma_{m}(\theta):=\left\{(x,\psi_{1}(x,\lambda^{0})): x\in [x_{l}^{m}(\theta), x_{m}^{m}(\theta)]\right\}
       \cup \left\{(x,y_{m}-\theta): x\in [x_{l}^{m}(\theta), x_{m}^{m}(\theta)]\right\},$$
and for $y_{m}-y_{M}\leq \theta\leq 2(y_{m}-y_{M})$,
\begin{eqnarray*}
\Gamma_{m}(\theta)\!\!\!&:=&\!\!\!
        \left\{(x,\psi_{1}(x,\lambda^{0})): x\in [x_{l}, x_{m}^{m}(2(y_{m}-y_{M})-\theta)]\right\}\\
       \!\!\!& &\!\!\! \cup \left\{(x,2y_{M}+\theta-y_{m}): x\in [x_{m}^{m}(2(y_{m}-y_{M})-\theta),x_{r}^{m}(2(y_{m}-y_{M})-\theta)]\right\}\\
        \!\!\!& &\!\!\! \cup \left\{(x,\psi_{1}(x,\lambda^{0})): x\in [x_{M},x_{r}^{m}(2(y_{m}-y_{M})-\theta)]\right\}\\
         \!\!\!& &\!\!\! \cup \left\{(x,y_{M}): x\in [x_{l}, x_{M}]\right\}.
\end{eqnarray*}
Similarly, we can define the family of slow-fast cycles $\Gamma_{M}(\cdot)$ for the canard point $(x_{M},y_{M})$,
the detail is omitted.
Then we have the following statements.
\begin{theorem}
\label{thm-one}
Assume that  for $\lambda=\lambda^{0}$ and $v=v^{0}$,
the slow-fast system (\ref{2D-model-14}) has a unique equilibrium $(x_{0},y_{0})$ in the set $x\geq 0$.
Then for  $\lambda=\lambda^{0}$, $v=v^{0}$ and sufficiently small $\varepsilon>0$,
the following statements hold:

\item{\bf (i)}
if the equilibrium $(x_{0},y_{0})$ is in the set $L^{1}$ (resp. $R^{1}$),
then system (\ref{2D-model-14}) has no periodic orbits in the set $\mathbb{R}^{2}_{+}$,
and $(x_{0},y_{0})$ is a stable node
and attracts the set $\mathbb{R}^{2}_{+}$ under the flow of system (\ref{2D-model-14}).

\item{\bf (ii)}
if the equilibrium $(x_{0},y_{0})$ is in the set $M$,
then for sufficiently small $\varepsilon>0$,
the equilibrium $(x_{0},y_{0})$ is an unstable node,
and there exists a unique limit cycle $\Gamma_{r,\varepsilon}$
in a small neighborhood of the slow-fast cycle $\Gamma_{r}$.
Furthermore, the limit cycle $\Gamma_{r,\varepsilon}$ is locally asymptotically stable
with the Floquet exponent bounded above by $-C/\varepsilon$ for some $C>0$,
and $\Gamma_{r,\varepsilon}\to \Gamma_{r}$ as $\varepsilon\to 0$ in the sense of Hausdorff distance.

\item{\bf (iii)}
if the equilibrium $(x_{0},y_{0})$ is in the set $L^{0}$ (resp. $R^{0}$),
then $(x_{0},y_{0})$ is a stable focus.

Further,
for the intersection point sequences $L^{0}$ and $R^{0}$,
let $\lambda=\lambda^{0}$ be fixed and the parameter $v$ vary.
Then for sufficiently small $\varepsilon>0$, the following assertions hold:

\item{\bf (iv)}
there exists a $V_{0}>0$ such that for each $v$ with $|v-v^{0}|<V_{0}$,
system (\ref{2D-model-14}) possesses a unique equilibrium near $(x_{m},y_{m})$ (resp. $(x_{M},y_{M})$) in the set $x\geq 0$,
which converges to $(x_{m},y_{m})$ (resp. $(x_{M},y_{M})$) as $(v,\varepsilon)\to (v^{0},0)$.
Moreover, there exist two Hopf bifurcation curves $v^{H}_{i}$ defined by
\begin{eqnarray}
\label{df-Hopf-bif}
v^{H}_{i}(\varepsilon)=v^{0}+
\frac{\kappa_{i,3}(D_{1}\psi_{2}(x_{i},\lambda^{0},v^{0}))^{2}(c^{0}+(x_{i}-\phi(x_{i}))^{2})}
{D_{11}\psi_{1}(x_{i},\lambda^{0})}\varepsilon+O(\varepsilon^{3/2}), \ \ \ i=m,M,
\end{eqnarray}
such that this equilibrium is stable for $v<v^{H}_{m}(\varepsilon)$ (resp. $v>v^{H}_{M}(\varepsilon)$)
and is unstable for $v>v^{H}_{m}(\varepsilon)$ (resp. $v<v^{H}_{M}(\varepsilon)$).
These Hopf bifurcations are nondegenerate if the constants $A_{i}$ given by (\ref{df-A-i}) satisfy $A_{i}\neq 0$, $i=m,M$,
and are supercritical for $A_{m}<0$ (resp. $A_{M}>0$) and are subcritical for $A_{m}>0$ (resp. $A_{M}<0$).

\item{\bf (v)}
fix some $\gamma \in (0,1)$  and assume that $A_{i}$ defined by  (\ref{df-A-i}) satisfy $A_{i}\neq 0$.
Then for each $i=m,M$,
there exists a smooth family of periodic orbits
\begin{eqnarray*}
(\theta, \varepsilon) \to (v_{i}(\theta,\varepsilon),\, \Gamma_{i}(\theta,\varepsilon)),
\ \ \ \varepsilon\in (0,\varepsilon_{0}),\ \ \theta\in (0,2(y_{m}-y_{M})),
\end{eqnarray*}
such that $\Gamma_{i}(\theta,\varepsilon)\to \Gamma_{i}(\theta)$ as $\varepsilon\to0$.
More precisely,
the periodic orbit $\Gamma_{i}(\theta,\varepsilon)$ is $O(\varepsilon^{\gamma})$-close to the canard point $(x_{i},y_{i})$
for each $\theta\in \left(0, \left(-D_{1}\psi_{2}(x_{i},\lambda^{0},v^{0})\varepsilon\right)^{\gamma}\right)$,
a relaxation oscillation for
each $\theta\in \left(2y_{m}-\left(-D_{1}\psi_{2}(x_{i},\lambda^{0},v^{0})\varepsilon\right)^{\gamma},2y_{m}\right)$,
and  a canard cycle for $v=v_{i}(\theta,\varepsilon)$ and
each $\theta\in \left[\left(-D_{1}\psi_{2}(x_{i},\lambda^{0},v^{0})\varepsilon\right)^{\gamma},
2y_{m}-\left(-D_{1}\psi_{2}(x_{i},\lambda^{0},v^{0})\varepsilon\right)^{\gamma}\right]$,
here $v_{i}(\theta,\varepsilon)$ satisfies
\begin{eqnarray}
|v_{i}(\theta,\varepsilon)-v^{c}_{i}(\varepsilon)|\leq
  \frac{D_{11}\psi_{1}(x_{i},\lambda^{0})}
  {2D_{1}\psi_{2}(x_{i},\lambda_{0},v^{0})(c^{0}+(x_{i}-\phi(x_{i}))^{2})}e^{-\left(-D_{1}\psi_{2}(x_{i},\lambda^{0},v^{0})\varepsilon\right)^{\gamma-1}},
\label{canard-vc}
\end{eqnarray}
where $v^{c}_{i}$ is in the form (\ref{df-canard-curve}).

\item{\bf (vi)}
if $(x_{0},y_{0})=(x_{m},y_{m})$ is a canard point,
then for $A_{m}>0$ and some $v$ with $v^{c}_{m}(\varepsilon)<v<v^{H}_{m}(\varepsilon)$,
there are two coexistent periodic orbits surrounding the equilibrium $(x_{m},y_{m})$,
where the inner one is unstable and the outer one is stable.
If $(x_{0},y_{0})=(x_{M},y_{M})$ is a canard point,
then for $A_{M}<0$ and some $v$ with $v^{c}_{M}(\varepsilon)<v<v^{H}_{M}(\varepsilon)$,
there are two coexistent periodic orbits surrounding the equilibrium $(x_{M},y_{M})$,
where the inner one is stable and the outer one is unstable.
\end{theorem}
{\bf Proof.}
We omitted the proofs for the types of the equilibria,
which can be obtained by a standard analysis.
The dynamics of the layer equations and the reduced systems are shown in Figure \ref{fg-all-cases}.

To prove {\bf (i)},
we only consider the case $(x_{0},y_{0})\in L^{1}$,
as the other one can be similarly proved.
Since the manifold $L^{1}$ is normally hyperbolic
and transversally intersects with $x$-axis,
then by \cite[Theorem 9.1]{Fenichel-79}
the manifold $L^{1}$ perturbs smoothly to locally invariant manifolds $L^{1}_{\varepsilon}$
which connects $(x_{0},y_{0})$ to a point at $x$-axis and transversally intersects with $x$-axis.
Then no periodic orbits surround $(x_{0},y_{0})$,
together with Theorem \ref{thm-attraction},
yields the attraction of $(x_{0},y_{0})$.
Thus, {\bf (i)} is obtained.

To prove {\bf (ii)},
assume that for $\lambda=\lambda^{0}$ and $v=v^{0}$ type $M$ appears.
By Lemmas \ref{lm-psi-1-prpty} and  \ref{lm-psi-0-2-prpty},
system (\ref{system-redu}) satisfies
$\dot y>0$ for $0<x<x_{m}$ and $\dot x<0$ for $x>x_{M}$,
and the stability of the critical manifold $\mathcal{C}_{0}$ changes at points $(x_{i},y_{i})$ for the layer equation.
The statements on the limit cycle $\Gamma_{r,\varepsilon}$ can be proved
by applying \cite[Theorem 2.1, p.318]{Krupa-Szmolyan-01JDE} and \cite[Theorem 9.1]{Fenichel-79}.
Thus, {\bf (ii)} is obtained.

To prove {\bf (iv)},
we recall that the existence and location of equilibria for the slow-fast system (\ref{2D-model-14})
are independent of $\varepsilon$,
then we can check that the first statement holds.
By \cite[formula (3.15), p.326]{Krupa-Szmolyan-01JDE}, for each $i=m,M$,
the Hopf bifurcation curve $\widehat{V}^{H}_{i}$ for the normal form (\ref{Normal-form-2D}) is in the form
\begin{eqnarray*}
\widehat{V}^{H}_{i}(\varepsilon)=-\frac{\kappa_{i3}}{2}\varepsilon+O(\varepsilon^{3/2})
         =-\frac{1}{2D_{1}\psi_{2}(x_{i},\lambda^{0},v^{0})}\varepsilon+O(\varepsilon^{3/2}).
\end{eqnarray*}
Thus by the transformation $\mathcal{T}_{1}^{-1}\circ \mathcal{T}_{2}^{-1}$,
we obtain the Hopf bifurcation curve given by (\ref{df-Hopf-bif}).
For canard point $(x_{m},y_{m})$ (resp. $(x_{M},y_{M})$),
the transformation $\mathcal{T}_{2}$ does not change (resp. changes) the sign of  $v$,
then from \cite[Theorem 3.1]{Krupa-Szmolyan-01JDE} it follows that
the remaining statements in {\bf (iv)} hold.

To prove {\bf (v)},
we first consider the normal form (\ref{Normal-form-2D}) of the slow-fast system (\ref{2D-model-14}) near the canard points $(x_{i},y_{i})$,
then by applying Theorems 3.3 and 3.5 in \cite{Krupa-Szmolyan-01JDE},
we can prove {\bf (v)} by similar method used in the proof for {\bf (iv)}.

To prove {\bf (vi)}, we only consider the case $(x_{0},y_{0})=(x_{m},y_{m})$,
as the other one can be similarly proved.
Assume that $A_{m}>0$.
Then by $D_{11}\psi_{1}(x_{m},\lambda^{0})<0$, (\ref{df-canard-curve}) and (\ref{df-Hopf-bif}),
we have that
\begin{eqnarray*}
v^{c}_{m}(\varepsilon)<v^{H}_{m}(\varepsilon), \ \ \ \ \
v^{H}_{m}(\varepsilon)-v^{c}_{m}(\varepsilon)=O(\varepsilon)
\end{eqnarray*}
for sufficiently small $\varepsilon>0$,
where $v^{c}_{m}(\varepsilon)$ and $v^{H}_{m}(\varepsilon)$
control the Hopf bifurcation and the intersection of slow manifolds near $(x_{m},y_{m})$, respectively.
Let a sufficiently small $\varepsilon>0$ be fixed
and vary $v$ from $v^{c}_{m}(\varepsilon)$ to  $v^{H}_{m}(\varepsilon)$.
When $v$ is in an exponentially small neighborhood of $v^{c}_{m}(\varepsilon)$
and satisfies $v^{c}_{m}(\varepsilon)<v<v^{H}_{m}(\varepsilon)$,
by Lemma \ref{thm-connection} and {\bf (v)} in this theorem
we have $x_{m,l}>x_{m,m}$ and a canard cycle with head appears.
By the bifurcation diagram in \cite[Figure 7 (b), p. 328]{Krupa-Szmolyan-01JDE},
the amplitude of this limit cycle increases as $v$ increases
and this persistent limit cycle is a relaxation oscillation or  a stable canard cycle with head
for each $v$ in a small neighborhood of  $v^{H}_{m}(\varepsilon)$.
Then we obtain the outer limit cycle.
By {\bf (iv)} in this theorem,
the Hopf bifurcation is subcritical for $A_{m}>0$.
Then there exists  a sufficiently small $\tilde{V}_{0}>0$
such that for each $v$ with $0<v^{H}_{m}(\varepsilon)-v<\tilde{V}_{0}$,
an unstable limit cycle arises from  the subcritical Hopf bifurcation
and coexists with the obtained large amplitude  limit cycle.
Thus, two coexistent periodic orbits are obtained and {\bf (vi)} is proved.
This finishes the proof.
\hfill$\Box$

\subsubsection{Two equilibria}

Assume that the slow-fast system (\ref{2D-model-14}) has precisely two equilibria in the set $x\geq 0$
for some $\lambda=\lambda^{0}$ and $v=v^{0}$.
Then
all possible intersection point sequences are as follows:
$L^{0}M$, $L^{1}M$, $MM$,  $MR^{0}$ and $MR^{1}$.
See Figures \ref{f-2-1}, \ref{f-2-2}, \ref{f-2-3}, \ref{f-2-4} and \ref{f-2-5}.

We first show that one of equilibria in  $M$ is a saddle-node
and  the slow-fast system (\ref{2D-model-14}) undergoes saddle-node bifurcation \cite[Section 3.4]{Guck-Holmes-83}
as the parameter $v$ varies.

\begin{theorem}
\label{thm-SN}
Assume that for $\lambda=\lambda^{0}$ and $v=v^{0}$,
the slow-fast system (\ref{2D-model-14}) has precisely two equilibria in the half plane $x\geq 0$.
Then the following statements hold:

\item{\bf (i)}
for sufficiently small $\varepsilon>0$,
system (\ref{2D-model-14}) has a saddle-node point $(x_{0},y_{0})\in M$,
at which  system (\ref{2D-model-14}) satisfies $D_{1}\psi_{1}(x_{0},\lambda^{0})=D_{1}\psi_{2}(x_{0},\lambda^{0},v^{0})$
and $D_{11}\psi_{1}(x_{0},\lambda^{0})\neq D_{11}\psi_{2}(x_{0},\lambda^{0},v^{0})$.

\item{\bf (ii)}
let $\lambda=\lambda^{0}$ be fixed and  the parameter $v$ vary.
Then system (\ref{2D-model-14}) undergoes a saddle-node bifurcation,
more precisely,
if $\psi_{1}(x,\lambda^{0})\leq \psi_{2}(x,\lambda^{0},v^{0})$
(resp. $\psi_{1}(x,\lambda^{0})\geq \psi_{2}(x,\lambda^{0},v^{0})$) near $x=x_{0}$,
then for small $|v-v^{0}|$,
system (\ref{2D-model-14}) has no equilibria near $(x_{0},y_{0})$ for $v>v^{0}$ (resp. $v<v^{0}$),
and system (\ref{2D-model-14}) has two equilibria $(x_{0}^{1},y_{0}^{1})$ and $(x_{0}^{2},y_{0}^{2})$
satisfying $x_{0}^{1}<x_{0}^{2}$
near $(x_{0},y_{0})$ for $v<v^{0}$ (resp. $v>v^{0}$),
where $(x_{0}^{1},y_{0}^{1})$ is an unstable node (resp. a saddle)
and $(x_{0}^{1},y_{0}^{1})$ is a saddle (resp. an unstable node).
\end{theorem}
{\bf Proof.}
Assume that system (\ref{2D-model-14}) has precisely two equilibria in the set $x\geq 0$
for $\lambda=\lambda^{0}$ and $v=v^{0}$,
then by Lemmas \ref{lm-psi-0-2-prpty} and \ref{distribution-inter},
there exists precisely one equilibrium $(x_{0},y_{0})$ in $M$,
which is a tangent point between functions $\psi_{1}$ and $\psi_{2}$,
that is, $D_{1}\psi_{1}(x_{0},\lambda^{0})=D_{1}\psi_{2}(x_{0},\lambda^{0},v^{0})$.
Then for sufficiently small $\varepsilon>0$,
the functions  $\mathcal{D}(\cdot,\,\cdot)$, $T(\cdot,\,\cdot)$ and $\Delta(\cdot,\,\cdot)$
defined by (\ref{df-D-T}) and (\ref{df-Delta}) satisfy
\begin{eqnarray*}
\mathcal{D}(x_0,y_0)=0, \ \  \ T(x_0,y_0)>0, \  \ \ \Delta(x_0,y_0)>0,
\end{eqnarray*}
and the eigenvalues of the Jacobian matrix $\mathcal{J}(x_{0},y_{0})$ are
$\mu_{1}=-\varepsilon-D_{1}\psi_{1}(x_{0},\lambda^{0})>0$ and $\mu_{2}=0$.
By a change
\begin{eqnarray*}
(x,y)\to (\bar{x}+\bar{y}+x_{0},D_{1}\psi_{1}(x_{0},\lambda^{0})\bar{x}-\varepsilon \bar{y}+y_{0}),
\end{eqnarray*}
and then dropping the bars over the variables,
we can change system (\ref{2D-model-14}) into
\begin{eqnarray}
\label{SN-normal-form}
\begin{aligned}
\frac{dx}{dt}&= X_{2}(x+y),
\label{SN-normal-form-1}
\\
\frac{dy}{dt}&= \mu_{1}y+Y_{2}(x+y),
\label{SN-normal-form-2}
\end{aligned}
\end{eqnarray}
where $X_{2}$ and $Y_{2}$ are given by
\begin{eqnarray*}
X_{2}(x)\!\!\!&=&\!\!\! \frac{\varepsilon}{D_{1}\psi_{1}(x_{0},\lambda^{0})+\varepsilon}
        \left(\psi_{2}(x+x_{0},\lambda^{0},v^{0})-\psi_{1}(x+x_{0},\lambda^{0})\right),\\
Y_{2}(x) \!\!\!&=&\!\!\!
       - \frac{1}{D_{1}\psi_{1}(x_{0},\lambda^{0})+\varepsilon}
       \left(D_{1}\psi_{1}(x_{0},\lambda^{0})\psi_{1}(x+x_{0},\lambda^{0})+\varepsilon\psi_{2}(x+x_{0},\lambda^{0},v^{0})\right)\\
        \!\!\!& &\!\!\! + D_{1}\psi_{1}(x_{0},\lambda^{0})x+y_{0}.
\end{eqnarray*}
Clearly, $X_{2}(0)=Y_{2}(0)=X_{2}^{'}(0)=Y_{2}^{'}(0)=0$.
Then by the  Implicit Function Theorem,
there exists a smooth function $y=y(x)$ with $y(0)=y'(0)=0$ such that $\mu_{1}y(x)+Y_{2}(x,y(x))=0$
in a neighbourhood of $(0,0)$.
By a direct computation,  for small $|x|$ the  function $X_{2}(\cdot+y(\cdot))$ can be expanded as the form
$$X_{2}(x+y(x))= K_{2} x^{2}+O(x^{3}),$$
where  the coefficient $K_{2}$ is in the form
\begin{eqnarray*}
K_{2}=\frac{\varepsilon \left(D_{11}\psi_{2}(x_{0},\lambda^{0},v^{0})-D_{11}\psi_{1}(x_{0},\lambda^{0})\right)}
      {D_{1}\psi_{1}(x_{0},\lambda^{0})+\varepsilon}.
\end{eqnarray*}
By Lemma \ref{lm-psi-0-2-prpty} we have $K_{2}\neq 0$.
Thus, \cite[Theorem 2.19, p.74]{Dumortieretal06} yields that
the equilibrium $(x_{0},y_{0})\in M$ is a saddle-node.
Then {\bf (i)} holds.

\begin{figure}[!htbp]
\centering
\subfigure[$v<0$]{
\begin{minipage}[t]{0.32\linewidth}
\centering
\includegraphics[width=1in]{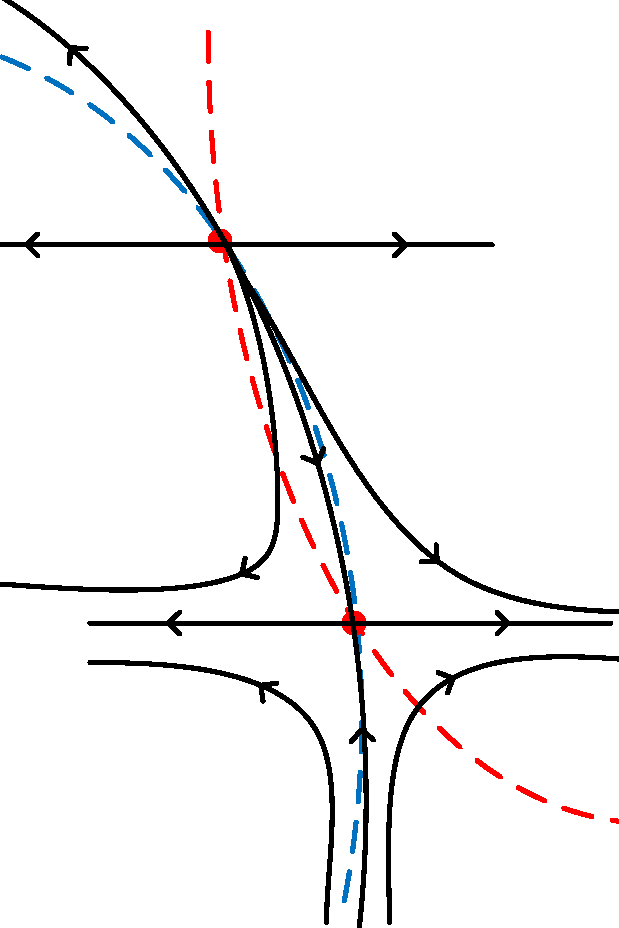}
\end{minipage}%
}%
\subfigure[$v=0$]{
\begin{minipage}[t]{0.32\linewidth}
\centering
\includegraphics[width=0.8in]{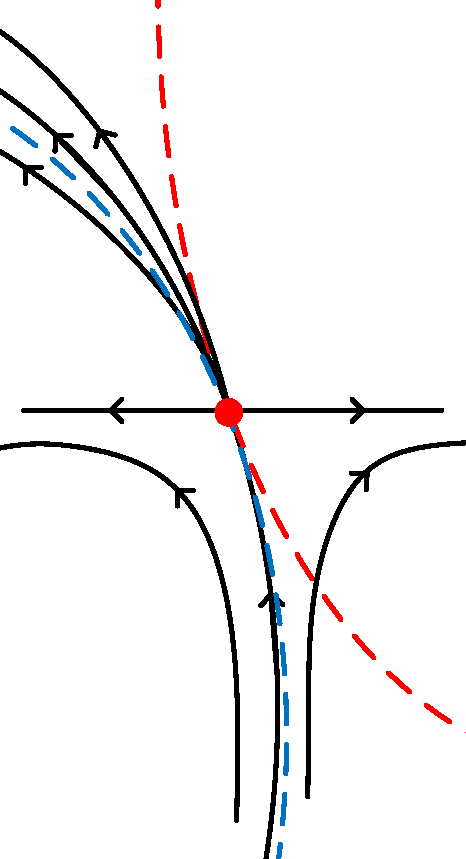}
\end{minipage}
}%
\subfigure[$v>0$]{
\begin{minipage}[t]{0.32\linewidth}
\centering
\includegraphics[width=1in]{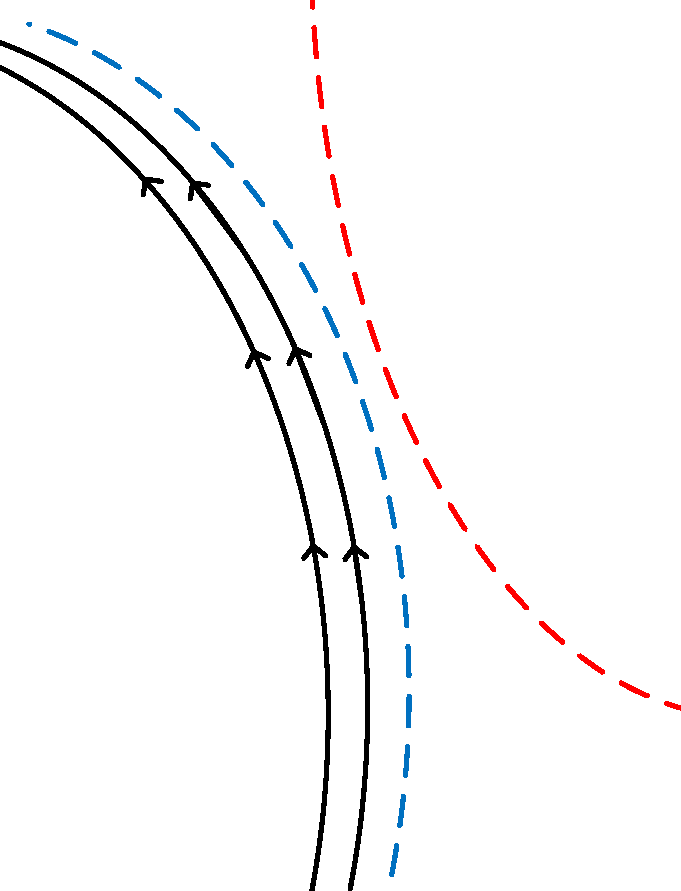}
\end{minipage}
}%
\centering
\caption{Saddle-node bifurcation.}
\label{Fig-SN-Bif}
\end{figure}

To prove {\bf (ii)},
we only consider the case that $\psi_{1}(x,\lambda^{0})\leq \psi_{2}(x,\lambda^{0},v^{0})$
for small $|x-x_{0}|$,
as the other case can be similarly discussed.
Then we have
$$D_{11}\psi_{2}(x_{0},\lambda^{0},v^{0})-D_{11}\psi_{1}(x_{0},\lambda^{0})>0.$$
Consider (\ref{SN-normal-form}) with  $v^{0}$ replaced by $v+v^{0}$.
By the Center Manifold Theory \cite[Section 1.3]{Carr-81},
the flow on the center manifold for an equivalent system of (\ref{SN-normal-form}) is governed by
\begin{eqnarray}
 \label{system-CM}
\begin{aligned}
\frac{dx}{dt}=&\, \frac{\varepsilon}{D_{1}\psi_{1}(x_{0},\lambda^{0})+\varepsilon}
        \left(\left(\frac{1}{c^{0}+(x_{0}-\phi(x_{0}))^{2}}+D_{13}\psi_{2}(x_{0},\lambda^{0},v^{0})x\right)v\right.\\
       &\left.+\frac{1}{2}\left( D_{11}\psi_{2}(x_{0},\lambda^{0},v^{0})-D_{11}\psi_{1}(x_{0},\lambda^{0})\right)x^{2}\right)
        +O(|(x,v)|^{3})\\
\frac{dv}{dt}=&\, 0.
\end{aligned}
\end{eqnarray}
The proof for (\ref{system-CM}) is given in Appendix B.
Since $D_{1}\psi_{1}(x_{0},\lambda^{0})<0$,
then for sufficiently small $\varepsilon$,
system (\ref{system-CM}) has no equilibria near $x=0$ for $v>0$
and has two equilibria $x=x_{1}(v)$  and $x=x_{2}(v)$ with $x_{1}(v)<x_{2}(v)$ near $x=0$ for $v<0$,
where $x=x_{1}(v)$ and $x=x_{2}(v)$ are an unstable node and a stable node, respectively.
See Figure \ref{Fig-SN-Bif}.
Then {\bf (ii)} holds.
Therefore, the proof is now complete.
\hfill$\Box$

By the above theorem,
we observe that the equilibrium of type $M$ in the sequences $L^{1}M$, $MR^{1}$, $L^{0}M$ and $MR^{0}$
is a saddle-node, so is one of the equilibria in the sequence $MM$.
More properties of the slow-fast system (\ref{2D-model-14})  with two equilibria are given in the next results.

\begin{theorem}
\label{thm-two}
Assume that the slow-fast system (\ref{2D-model-14}) has precisely two equilibria in the set $x\geq 0$
for $\lambda=\lambda^{0}$ and $v=v^{0}$.
Then for sufficiently small $\varepsilon>0$,
the following statements hold:

\item{\bf (i)}
if the intersection point sequence is $L^{1}M$ (resp. $MR^{1}$),
then system (\ref{2D-model-14}) has a stable node $(x^{1}_{0},y^{1}_{0})$ in  $L^{1}$ (resp. $R^{1}$),
a saddle-node $(x^{2}_{0},y^{2}_{0})$ in $M$, no periodic orbits in the set $x\geq0$
and infinitely many heteroclinic orbits joining $(x^{2}_{0},y^{2}_{0})$ to $(x^{1}_{0},y^{1}_{0})$.
Further, all orbits starting from the first quadrant including its boundary,
except a unique center manifold of $(x^{2}_{0},y^{2}_{0})$, converge to the stable node $(x^{1}_{0},y^{1}_{0})$
as time goes to infinity.

\item{\bf (ii)}
if the intersection point sequence is $MM$,
then system (\ref{2D-model-14}) has an unstable node $(x^{1}_{0},y^{1}_{0})\in M$,
a saddle-node $(x^{2}_{0},y^{2}_{0})\in M$,
and a unique heteroclinic orbit joining the unstable node to the saddle-node.

\item{\bf (iii)}
if the intersection point sequence is $L^{0}M$ (resp. $MR^{0}$),
then system (\ref{2D-model-14}) has a stable focus $(x^{1}_{0},y^{1}_{0})$ in  $L^{0}$ (resp. $R^{0}$)
and a saddle-node $(x^{2}_{0},y^{2}_{0})$ in $M$.
Let $\lambda=\lambda^{0}$ be fixed and the parameter $v$ satisfy $|v-v^{0}|\ll 1$.
Then system (\ref{2D-model-14}) has a homoclinic orbit, which closes to
either a canard slow-fast  cycle without head or a canard slow-fast cycle with head,
if and only if
$\kappa_{i,3}+A_{i}/4<0$ and $v=v^{c}_{i}(\varepsilon)$,
where the functions $v^{c}_{i}$  are defined by (\ref{df-canard-curve}).
Furthermore, if $\kappa_{i,3}+A_{i}/4<0$ and $0<v-v^{c}_{m}(\varepsilon)\ll 1$ (resp. $0<v^{c}_{M}(\varepsilon)-v\ll 1$),
then either an unstable canard cycle with head or an unstable canard cycle without head bifurcates from this homoclinic orbit.
\end{theorem}

\begin{figure}[!htbp]
\centering
\subfigure[$L^{1}M$]{
\begin{minipage}[t]{0.32\linewidth}
\centering
\includegraphics[width=1.7in]{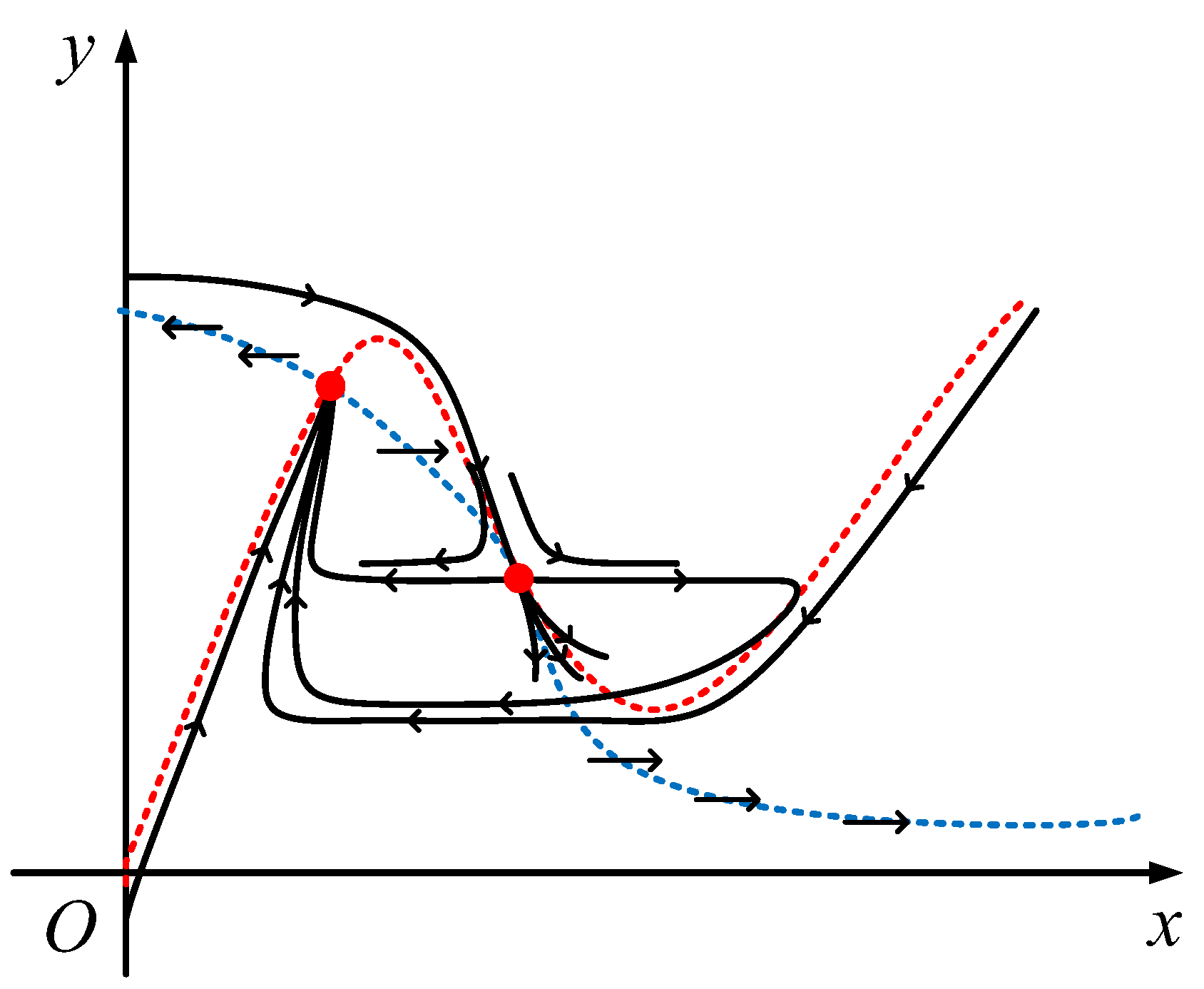}
\label{fig-S2-L1-M}
\end{minipage}%
}%
\subfigure[$MM$]{
\begin{minipage}[t]{0.32\linewidth}
\centering
\includegraphics[width=1.6in]{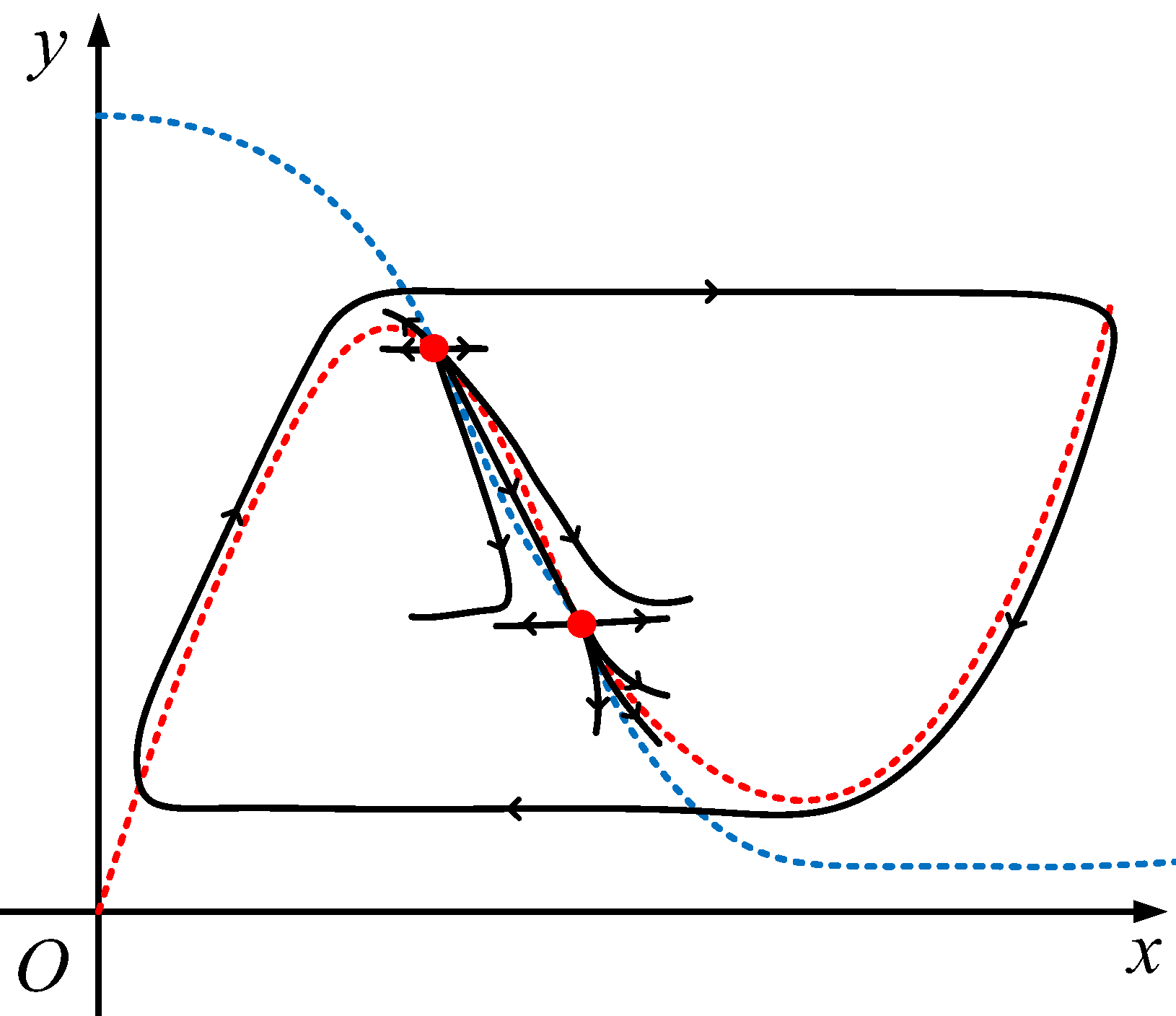}
\label{fig-S2-M-M}
\end{minipage}
}%
\subfigure[$L^{0}M$]{
\begin{minipage}[t]{0.32\linewidth}
\centering
\includegraphics[width=1.6in]{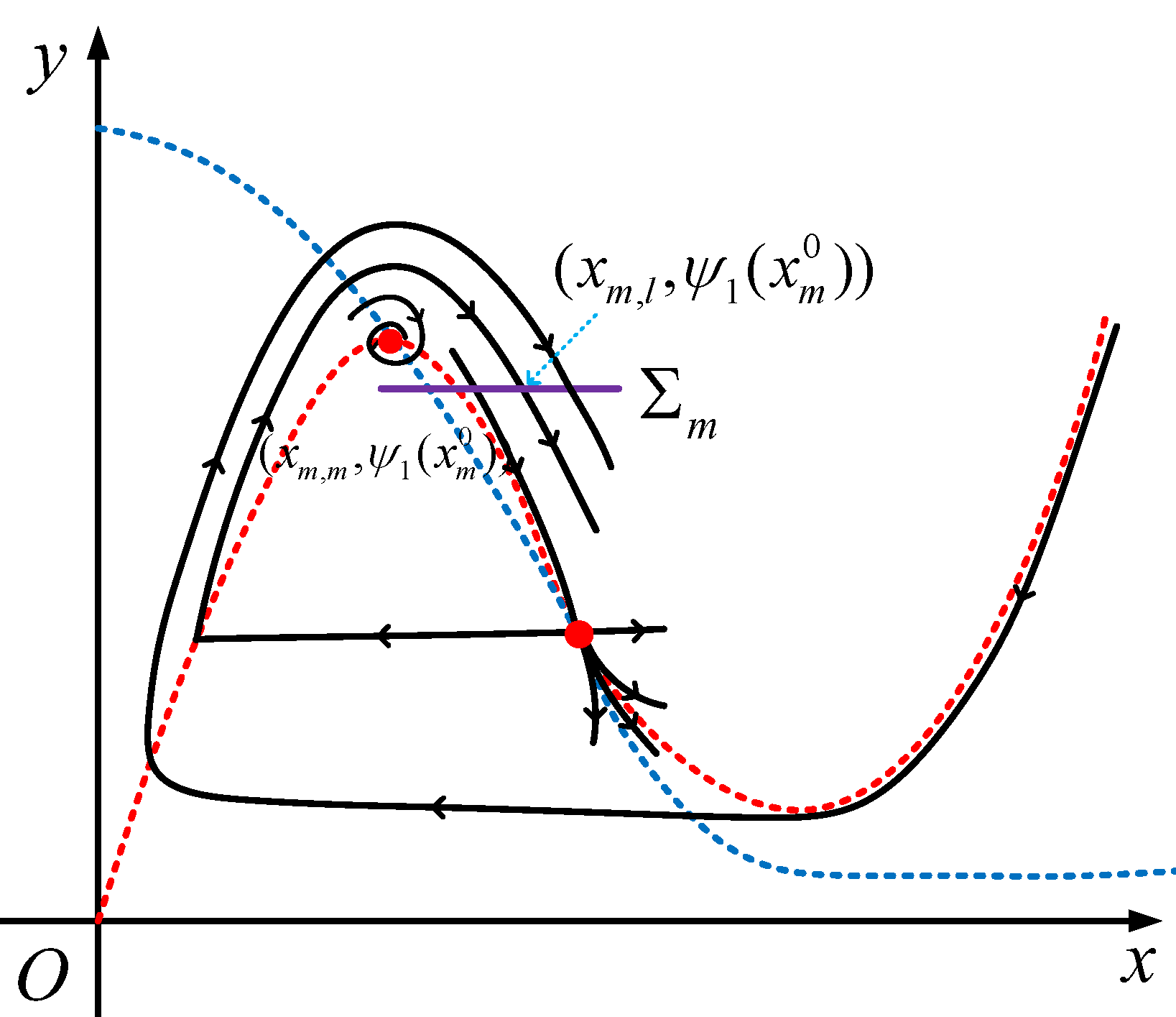}
\label{fig-S2-L0-M}
\end{minipage}
}%
\centering
\caption{Dynamics of the slow-fast system (\ref{2D-model-14}) with two equilibria in the set $x\geq 0$.
The solid black curves are the orbits of system (\ref{2D-model-14}),
the graphes of the functions $\psi_{1}$ and $\psi_{2}$ respectively indicate the dashed red and the dashed blue curves.
}
         \label{fig-S2-Dynamics}
\end{figure}

Throughout the proof for this theorem,
we omit the proofs for the types of the equilibria.
Dynamics of the cases $L^{1}M$, $MM$ and $L^{0}M$  are illustrated by Figure \ref{fig-S2-Dynamics}.\
\\
{\bf Proof.}
To prove {\bf (i)}, we only give the proof for type $L^{1}M$.
Similarly to Theorem \ref{thm-one} {\bf (i)},
system (\ref{2D-model-14}) with sufficiently small $\varepsilon>0$ has no periodic orbits surrounding $(x^{1}_{0},y^{1}_{0})$.
Clearly, along the curve $y=\psi_{2}(x,\lambda^{0})$ for $x>x_{0}^{1}$ and $x\neq x_{0}^{2}$ we have $dx/dt<0$,
which yields that no periodic orbits surround  $(x_0^{2},y_0^{2})$.
Thus, no periodic orbits exist in the first quadrant.
By Theorem \ref{thm-SN} we obtain that the saddle-node $(x_0^{2},y_0^{2})$
possesses a unique center manifold approaching to it and infinitely many center manifolds leaving it.
Hence,  there are  infinitely many orbits, which leave the saddle-node point $(x_{0}^{2},y_{0}^{2})$,
joining $(x_{0}^{2},y_{0}^{2})$ to $(x_{0}^{1},y_{0}^{1})$,
and a unique orbit approaching to $(x_{0}^{2},y_{0}^{2})$.
Thus, the proof for  {\bf (i)} is finished by using Theorem \ref{thm-attraction}.

To prove {\bf (ii)},
assume that  $(x^{1}_{0},y^{1}_{0})$ and $(x^{2}_{0},y^{2}_{0})$
are a transversal point and a tangent point of the functions $\psi_{1}$ and $\psi_{2}$, respectively.
Without loss of generality, assume that $x_{0}^{1}<x_{0}^{2}$ (see Figure \ref{fig-S2-M-M}).
By Theorem \ref{thm-SN},
there are a unique center manifold on which  the orbit approaches to $(x^{2}_{0},y^{2}_{0})$ from the above
and infinitely many orbits leaving $(x^{2}_{0},y^{2}_{0})$.
The existence and uniqueness of heteroclinic orbits are derived from the persistence of normally hyperbolic invariant manifolds
and the uniqueness of the orbits approaching to $(x^{2}_{0},y^{2}_{0})$.
Thus, the proof for {\bf (ii)} is finished.

\begin{figure}[!htbp]
\centering
\subfigure[Small homoclinic orbit.]{
\begin{minipage}[t]{0.4\linewidth}
\centering
\includegraphics[width=2in]{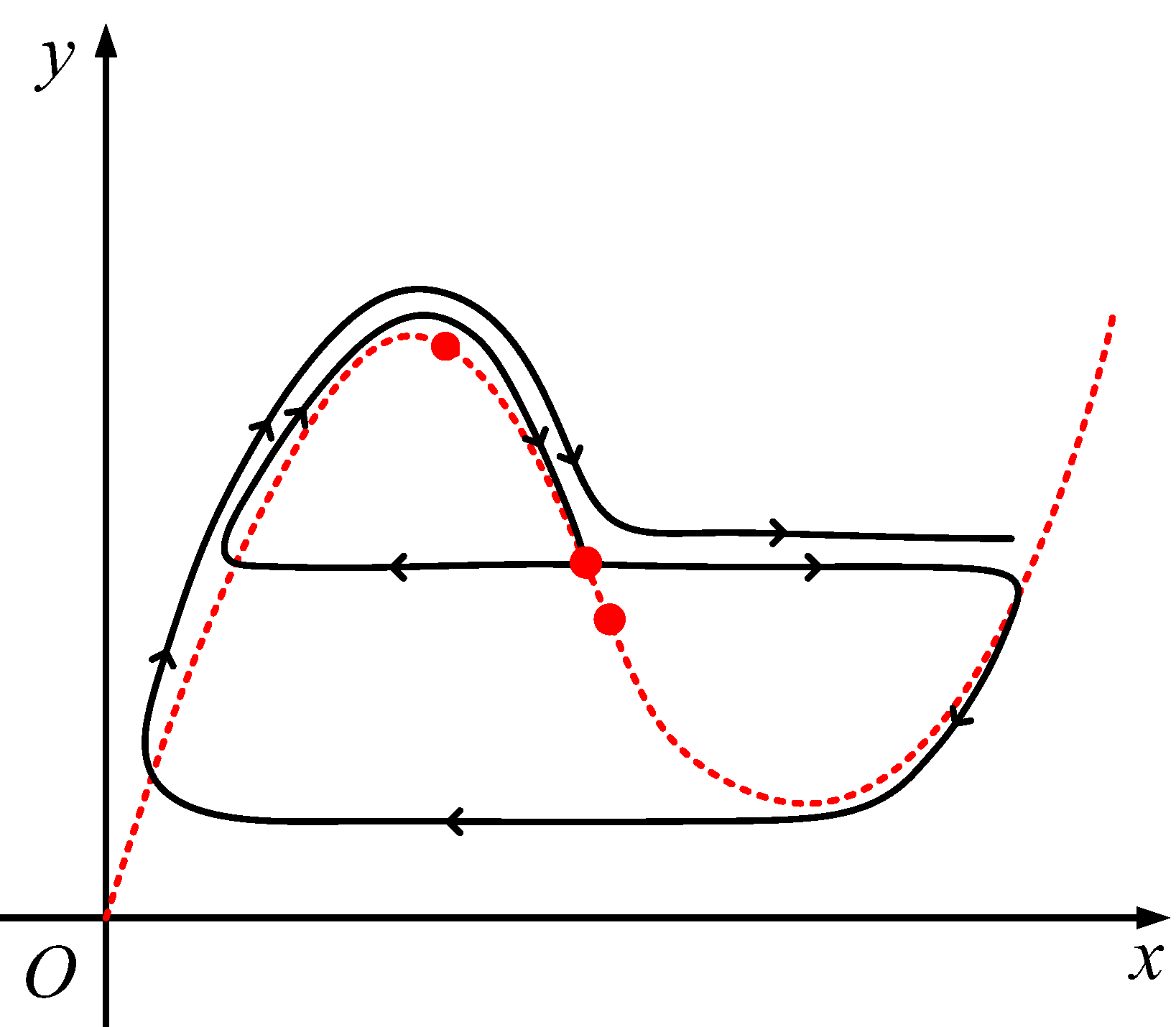}
\label{fig-S2-L0-M-small}
\end{minipage}%
}%
\subfigure[Big homoclinic orbit.]{
\begin{minipage}[t]{0.4\linewidth}
\centering
\includegraphics[width=2in]{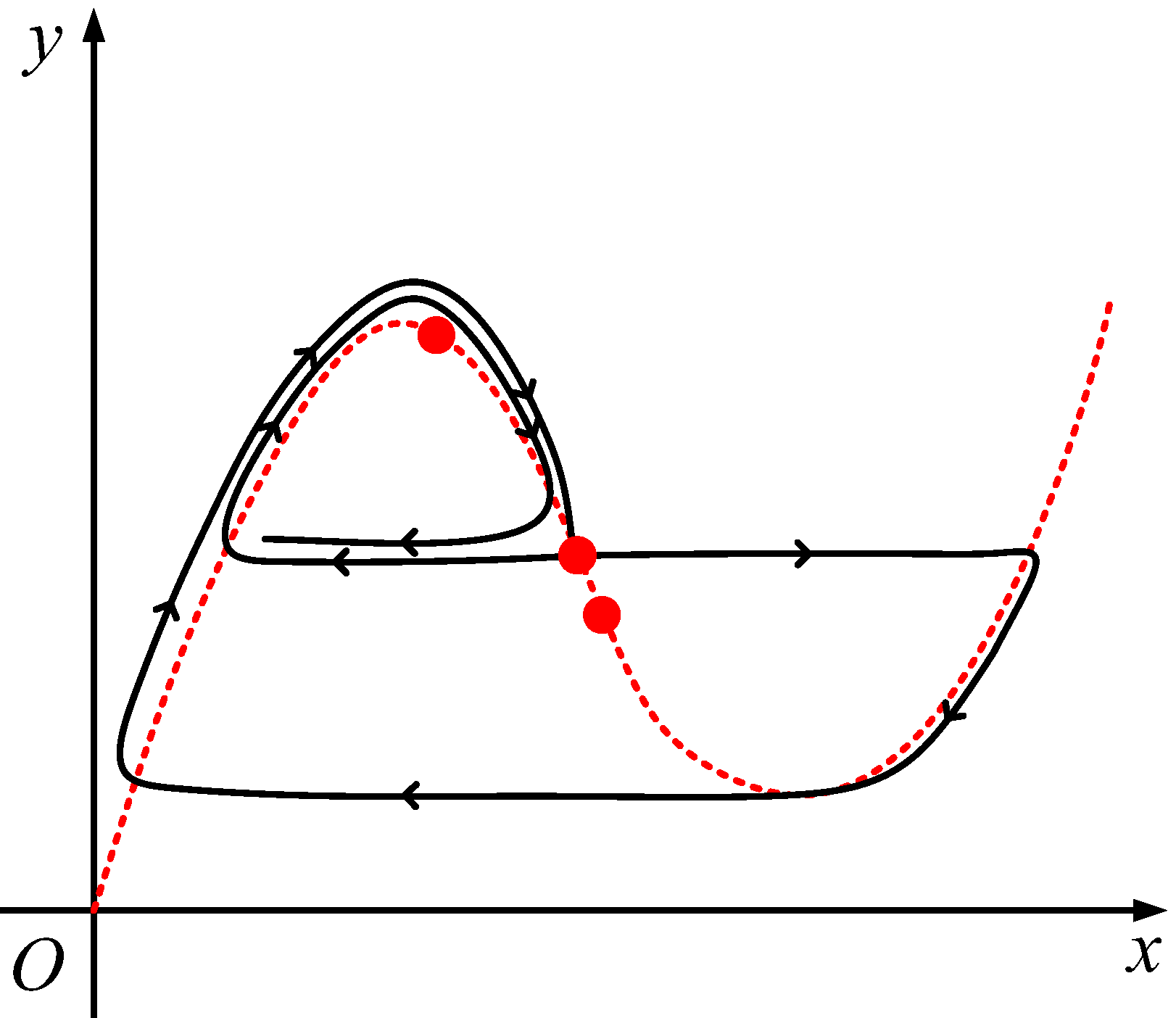}
\label{fig-S2-L0-M-big}
\end{minipage}
}%
\centering
\caption{Two possible homoclinic orbits arise in type $L^{0}M$.
The red dots are equilibria,
the solid black curves are the orbits of system (\ref{2D-model-14}),
and the graphes of the functions $\psi_{1}$ and $\psi_{2}$ respectively indicate the dashed red and the dashed blue curves.}
\end{figure}

To prove {\bf (iii)},
we only consider type $L^{0}M$ (see Figure \ref{fig-S2-L0-M}).
Let the notations be given as in Lemma \ref{thm-connection}.
If $\kappa_{m,3}+A_{m}/4<0$ and $v=v^{c}_{m}(\varepsilon)$,
then by $D_{11}\psi_{1}(x_{m},\lambda^{0})<0$ and (\ref{df-canard-curve}),
we obtain that $v=v^{c}_{m}(\varepsilon)>0$.
This together with Theorem \ref{thm-SN} yields that
there are a saddle $(\widetilde{x}^{2}_{0},\widetilde{y}^{2}_{0})$
and an unstable node $(x^{3}_{0},y^{3}_{0})$,
which bifurcate from the saddle-node $(x^{2}_{0},y^{2}_{0})$
and satisfy $x^{1}_{0}<\widetilde{x}^{2}_{0}<x^{3}_{0}$.
Lemma \ref{thm-connection} yields
the existence of the homoclinic orbit,
which is homoclinic to the saddle $(\widetilde{x}^{2}_{0},\widetilde{y}^{2}_{0})$
and together with this saddle forms
either a small loop near a canard slow-fast  cycle without head (see Figure \ref{fig-S2-L0-M-small})
or a big one  near a canard slow-fast cycle with head (see Figure \ref{fig-S2-L0-M-big}).
If either $\kappa_{m,3}+A_{m}/4>0$ or $v\neq v^{c}_{m}(\varepsilon)$,
then by Lemma \ref{thm-connection}  and Theorem \ref{thm-SN},
no homoclinic orbits exist for system (\ref{2D-model-14}) with sufficiently small $\varepsilon$.
To prove the last statement,
assume that $\kappa_{m,3}+A_{m}/4<0$ and $0<v-v^{c}_{m}(\varepsilon)\ll 1$,
then by Lemma \ref{thm-connection} we obtain that $x_{m,l}>x_{m,m}$.
Since the first order saddle quantity $T(\widetilde{x}^{2}_{0},\widetilde{y}^{2}_{0})$ of the saddle $(\widetilde{x}^{2}_{0},\widetilde{y}^{2}_{0})$ satisfies $T(\widetilde{x}^{2}_{0},\widetilde{y}^{2}_{0})>0$ for sufficiently small $\varepsilon$,
then by \cite[Theorem 3.3, p. 357]{Chow-Hale-82}
an unstable periodic orbit bifurcating from this homoclinic orbit is either a canard cycle without head
if the homoclinic orbit is small or a canard cycle with head if the homoclinic orbit is big.
Thus, we obtain {\bf (iii)}.
Therefore, the proof is now complete.
\hfill$\Box$

\subsubsection{Three equilibria}

Assume that
the slow-fast system (\ref{2D-model-14}) possesses three equilibria
for some $\lambda=\lambda^{0}$ and $v=v^{0}$.
Then all  possible intersection point sequences are
$L^{0}MR^{0}$, $L^{0}MR^{1}$, $L^{1}MR^{0}$, $L^{1}MR^{1}$, $L^{0}MM$, $L^{1}MM$, $MMM$, $MMR^{0}$ and $MMR^{1}$.
See Figures \ref{f-3-1}, \ref{f-3-2}, \ref{f-3-3}, \ref{f-3-4}, \ref{f-3-5}, \ref{f-3-6}, \ref{f-3-7}, \ref{f-3-8} and \ref{f-3-9}.
The main results for this case are summarized as follows.

\begin{theorem}
\label{thm-three}
Assume that the slow-fast system (\ref{2D-model-14}) has precisely three  equilibria $(x^{i}_{0},y^{i}_{0})$,
$i=1,2,3$, in the set $x\geq 0$ for $\lambda=\lambda^{0}$ and $v=v^{0}$,
where $x^{1}_{0}<x^{2}_{0}<x^{3}_{0}$.
Then for sufficiently small $\varepsilon>0$, the following statements hold:

\item{\bf (i)}
if the intersection point sequence is $L^{1}MR^{1}$,
then $(x^{1}_{0},y^{1}_{0})\in L^{1}$ and $(x^{3}_{0},y^{3}_{0})\in R^{1}$ are stable nodes
and $(x^{2}_{0},y^{2}_{0})\in M$ is a saddle,
system  (\ref{2D-model-14}) has no periodic orbits in the set $x\geq 0$,
and two heteroclinic orbits joining $(x^{2}_{0},y^{2}_{0})$ to $(x^{1}_{0},y^{1}_{0})$
and $(x^{2}_{0},y^{2}_{0})$ to $(x^{3}_{0},y^{3}_{0})$, respectively.
Furthermore,
the set $\mathcal{A}$ defined as in Theorem \ref{thm-attraction} is divided into
two disjoint sets $\Omega_{1}$ and $\Omega_{2}$ by the stable manifolds of $(x^{2}_{0},y^{2}_{0})$,
and all orbits starting from the interior of $\Omega_{1}$ (resp. $\Omega_{2}$)
converge to  $(x^{1}_{0},y^{1}_{0})$ (resp. $(x^{3}_{0},y^{3}_{0})$)
as time goes to infinity.

\item{\bf (ii)}
if the intersection point sequence is $MMM$,
then $(x^{1}_{0},y^{1}_{0})$ and $(x^{3}_{0},y^{3}_{0})$ are unstable nodes
and $(x^{2}_{0},y^{2}_{0})$ is a saddle,
and a locally asymptotically stable relaxation oscillation $\Gamma_{r,\varepsilon}$
arising from the  singular relaxation cycle $\Gamma_{r}$
approaches to $\Gamma_{r}$ in the sense of Hausdorff distance as $\varepsilon\to 0$,
where  the singular relaxation cycle $\Gamma_{r}$ is constructed as in Figure \ref{fig-relaxtion}.

\item{\bf (iii)}
if the intersection point sequence is $L^{1}MM$ (resp. $MMR^{1}$),
then $(x^{1}_{0},y^{1}_{0})\in L^{1}$ (resp. $(x^{3}_{0},y^{3}_{0})\in R^{1}$) is a stable node,
$(x^{2}_{0},y^{2}_{0})\in M$ is a saddle
and $(x^{3}_{0},y^{3}_{0})\in M$ (resp. $(x^{1}_{0},y^{1}_{0})\in M$) is an unstable node,
and system  (\ref{2D-model-14}) has no periodic orbits in the first quadrant,
a heteroclinic orbit connecting $(x^{2}_{0},y^{2}_{0})$ to $(x^{3}_{0},y^{3}_{0})$ (resp. $(x^{1}_{0},y^{1}_{0})$),
two heteroclinic orbits connecting $(x^{2}_{0},y^{2}_{0})$ to $(x^{1}_{0},y^{1}_{0})$ (resp. $(x^{3}_{0},y^{3}_{0})$)
and infinitely many heteroclinic orbtis connecting $(x^{3}_{0},y^{3}_{0})$ to $(x^{1}_{0},y^{1}_{0})$.

\item{\bf (iv)}
if the intersection point sequence is $L^{0}MM$ (resp. $MMR^{0}$),
then $(x^{1}_{0},y^{1}_{0})\in L^{0}$ (resp. $(x^{3}_{0},y^{3}_{0})\in R^{0}$) is a stable focus,
$(x^{2}_{0},y^{2}_{0})\in M$ is a saddle
and $(x^{3}_{0},y^{3}_{0})\in M$ (resp. $(x^{1}_{0},y^{1}_{0})\in M$) is an unstable node,
and system  (\ref{2D-model-14}) has
a heteroclinic orbit connecting $(x^{2}_{0},y^{2}_{0})$ to $(x^{3}_{0},y^{3}_{0})$ (resp. $(x^{1}_{0},y^{1}_{0})$).
Further, let $\lambda=\lambda^{0}$ be fixed and the parameter $v$ vary.
Then for $(x^{1}_{0},y^{1}_{0})\in L^{0}$ (resp. $(x^{3}_{0},y^{3}_{0})\in R^{0}$),
system (\ref{2D-model-14}) undergoes Hopf bifurcation and canard explosion
in the ways stated in Theorem \ref{thm-one}  {\bf (iv)} and Theorem \ref{thm-one}  {\bf (v)}, respectively.
\item{\bf (v)}
if the intersection point sequence is $L^{0}MR^{1}$ (resp. $L^{1}MR^{0}$),
then $(x^{1}_{0},y^{1}_{0})\in L^{0}$ (resp. $(x^{3}_{0},y^{3}_{0})\in R^{0}$) is a stable focus,
$(x^{2}_{0},y^{2}_{0})\in M$ is a saddle
and $(x^{3}_{0},y^{3}_{0})\in R^{1}$ (resp. $(x^{1}_{0},y^{1}_{0})\in L^{1}$) is a stable node.
Further, let $\lambda=\lambda^{0}$ be fixed and $v$ vary.
Then the following statements hold:

\begin{enumerate}
\item[{\bf (v.1)}]
system (\ref{2D-model-14}) undergoes Hopf bifurcation according to Theorem \ref{thm-one} {\bf (iv)}.

\item[{\bf (v.2)}]
system (\ref{2D-model-14}) has no relaxation oscillations
 or canard cycles with head as varying $v$ near $v^{0}$.

\item[{\bf (v.3)}]
there are two smooth functions $v^{c}_{i}$, $i=m,M$, having the expansions in (\ref{df-canard-curve})
such that system (\ref{2D-model-14})
possesses a homoclinic orbit, which is homoclinic to a saddle in  $M$ and
lies near a canard slow-fast cycle without head, if and only if $v=v^{c}_{i}(\varepsilon)$.

\item[{\bf (v.4)}]
if $0<v-v^{c}_{m}(\varepsilon)\ll 1$ (resp. $0<v^{c}_{M}(\varepsilon)-v\ll 1$),
then  an unstable canard cycle without head bifurcates from this homoclinic orbit.
If $0<v^{c}_{m}(\varepsilon)-v\ll 1$ (resp. $0<v-v^{c}_{M}(\varepsilon)\ll 1$),
then there exist no periodic orbits bifurcating from this homoclinic orbit.
\end{enumerate}

\item{\bf (vi)}
if the intersection point sequence is $L^{0}MR^{0}$,
then $(x^{1}_{0},y^{1}_{0})\in L^{0}$ is a stable focus,
$(x^{2}_{0},y^{2}_{0})\in M$ is a saddle
and $(x^{3}_{0},y^{3}_{0})\in R^{0}$ is a stable focus.
Further, let $\lambda=\lambda^{0}$  be fixed and $v$ vary.
Then the following statements hold:
\begin{enumerate}
\item[{\bf (vi.1)}]
system (\ref{2D-model-14}) undergoes a Hopf bifurcation near $(x_{m},y_{m})$ or $(x_{M},y_{M})$
according to the way stated in Theorem \ref{thm-one} {\bf (iv)}, but not simultaneously.

\item[{\bf (vi.2)}]
there are two smooth functions $v^{c}_{i}$, $i=m,M$, defined by (\ref{df-canard-curve})
such that system (\ref{2D-model-14}) has a homoclinic orbit, which is homoclinic to a saddle in $M$,
if and only if $v=v^{c}_{i}(\varepsilon)$.

\item[{\bf (vi.3)}]
assume that the constants $\mathcal{K}_{i}$ defined by (\ref{df-K-i}) satisfy $\mathcal{K}_{m}\neq \mathcal{K}_{M}$.
Then for  $v$ satisfying $0<v-v^{c}_{m}(\varepsilon)\ll 1$ (resp. $0<v^{c}_{M}(\varepsilon)-v\ll 1$),
there exists an unstable canard cycle bifurcating from the homoclinic orbit corresponding to $v=v^{c}_{m}(\varepsilon)$
(resp. $v=v^{c}_{M}(\varepsilon)$), and these two canard cycles can not appear simultaneously.
If $v$ satisfies $0<v^{c}_{m}(\varepsilon)-v\ll 1$ (resp. $0<v-v^{c}_{M}(\varepsilon)\ll 1$),
then there are no periodic orbits bifurcating from the homoclinic orbit corresponding to $v=v^{c}_{m}(\varepsilon)$
(resp. $v=v^{c}_{M}(\varepsilon)$).
\end{enumerate}
\end{theorem}
{\bf Proof.}
Here we also omitted the proofs for the types of the equilibria.

To prove {\bf (i)}, we first consider the existence of periodic orbits.
Similarly to Theorem \ref{thm-one} {\bf (i)},
no periodic orbits surround stable nodes $(x^{1}_{0},y^{1}_{0})$ and $(x^{3}_{0},y^{3}_{0})$.
Since $(x_{i}, y_{i})$, $i=m,M$, are jump points,
then by \cite[Theorem 2.1, p.290]{Krupa-Szmolyan-01SIMA}
the stable manifolds of $(x^{2}_{0},y^{2}_{0})$ extend to the boundary of the set $\mathcal{A}$.
Hence, the stable manifolds of $(x^{2}_{0},y^{2}_{0})$ cut $\mathcal{A}$ into two disjoint parts,
and no periodic orbits surround $(x^{2}_{0},y^{2}_{0})$.
Thus, no periodic orbits exist.
The invariant property of $\mathcal{A}$ yields the last statement.
Thus,  {\bf (i)} is proved.

Similarly to Theorem \ref{thm-one} {\bf (ii)},
we can obtain {\bf (ii)} in this theorem.

To prove {\bf (iii)}, we only consider type $L^{1}MM$.
Similarly to Theorem \ref{thm-one} {\bf (i)},
no periodic orbits surround $(x^{1}_{0},y^{1}_{0})$.
Since the manifold $M$ smoothly perturbs to locally invariant manifold $M_{\varepsilon}$,
which connects $(x^{2}_{0},y^{2}_{0})$ to $(x^{3}_{0},y^{3}_{0})$,
then system (\ref{2D-model-14}) with sufficiently small $\varepsilon$ has no periodic orbits in the first quadrant.
Thus,  {\bf (iii)} is obtained.

To prove {\bf (iv)},
for type  $L^{0}MM$ (resp. $MMR^{0}$),
the slow manifold $M_{\varepsilon}$ connects
$(x^{2}_{0},y^{2}_{0})$ to $(x^{3}_{0},y^{3}_{0})$ (resp. $(x^{1}_{0},y^{1}_{0})$).
Then the existence of the heteroclinic orbit is obtained.
The assertions {\bf (iv)} and {\bf (v)} in Theorem \ref{thm-one} yield that the last statement holds.
Thus, {\bf (iv)} is proved.

To prove {\bf (v)}, we only discuss type $L^{0}MR^{1}$.
Similarly to {\bf (iv)} in Theorem \ref{thm-one},  {\bf (v.1)} holds.
Since  $(x_{M},y_{M})$ is a jump point,
then by \cite[Theorem 2.1, p.290]{Krupa-Szmolyan-01SIMA}
the locally invariant manifold $M_{\varepsilon}$, which is a stable manifold of the saddle $(x^{2}_{0},y^{2}_{0})$,
can extend to the boundary of the invariant region $\mathcal{A}$.
Consequently, neither relaxation oscillations nor canard cycles with head appear.
Hence, {\bf (v.2)} holds.
The statements {\bf (v.3)} and {\bf (v.4)} can be similarly proved by the method used in Theorem \ref{thm-two} {\bf (iii)}.
Thus, {\bf (v)} is proved.

To prove {\bf (vi.1)},
by  Theorem \ref{thm-one} {\bf (iv)}, near $(x_{i},y_{i})$
Hopf bifurcations can take place by varying  $v$,
and the corresponding Hopf bifurcation curves $v^{H}_{i}(\cdot)$ are given by (\ref{df-Hopf-bif}).
Since $D_{11}\psi_{1}(x_{m},\lambda^{0})<0$, $D_{11}\psi_{1}(x_{M},\lambda^{0})>0$,
and $D_{1}\psi_{2}(x_{i},\lambda^{0},v^{0})<0$,
then  $v^{H}_{m}(\varepsilon)>0$ and $v^{H}_{M}(\varepsilon)<0$ for sufficiently small $\varepsilon$,
which implies that two Hopf bifurcations does not appear simultaneously.
Thus, {\bf (vi.1)} is proved.
Similarly to {\bf (v.3)} in this theorem, we can obtain {\bf (vi.2)}.
To prove {\bf (vi.3)}, assume that $\mathcal{K}_{m}\neq \mathcal{K}_{M}$.
Then by Lemma \ref{thm-connection}, two homoclinic orbits stated in {\bf (vi.2)} can not appear simultaneously.
By (\ref{canard-vc}) we obtain that canard cycles appear for the parameter $v$ in the exponentially small interval of $v^{c}_{i}(\varepsilon)$,
together with $|v^{c}_{m}(\varepsilon)- v^{c}_{M}(\varepsilon)|=O(\varepsilon)$,
yields that two canard cycles can not appear simultaneously.
The remaining statements can be proved by the  way in {\bf (v.4)}.
Thus, {\bf (vi)} is proved.
Therefore, the proof is  complete.
\hfill$\Box$

\subsection{Numerical examples}
\setcounter{equation}{0}
\setcounter{lemma}{0}
\setcounter{theorem}{0}
\setcounter{remark}{0}

Now we give several concrete numerical examples to illustrate the obtained results as follows.

\begin{example}
Let the parameters $a$, $b_{1}$, $b_{2}$, $c$, $\varepsilon$ and $v$ satisfy
$a=0.01$, $b_{1}=20$, $b_{2}=0.1$, $c=1$, $\varepsilon=0.01$ and $v=37.9$
in system (\ref{2D-model-12}).
A numerical simulation shows that there exists a big limit cycle enclosing a small one.
This indicates the coexistence of two limit cycles.
\begin{figure}[!htp]
  \centering
\includegraphics[width=6.6cm,height=5.4cm]{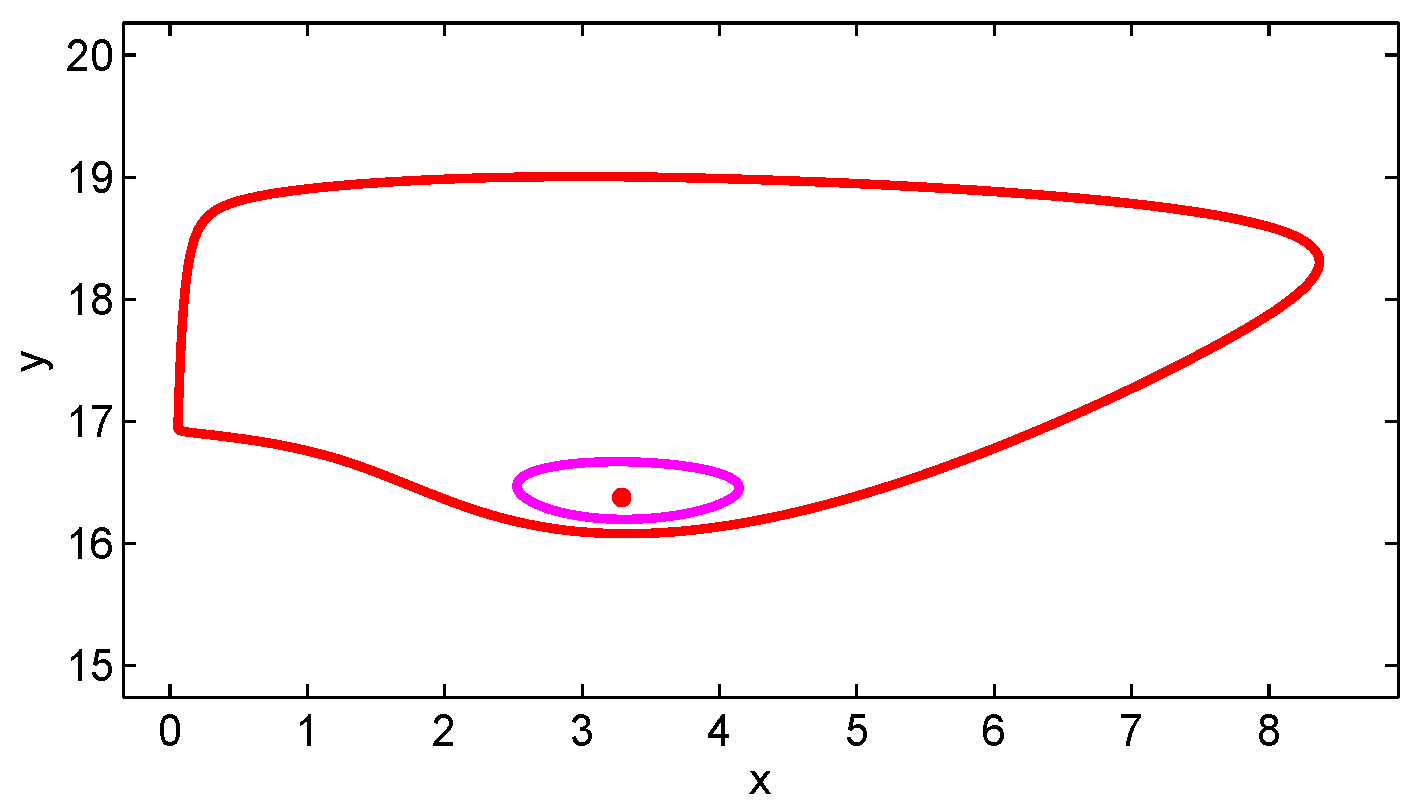}
\caption{A big limit cycle (the red cycle) encloses a small limit cycle (the mauve cycle).
The red point indicates an equilibrium.
}
\label{fig-S1-coexistence}
\end{figure}
\end{example}

\begin{example}
Let the parameters $a$, $b_{1}$, $b_{2}$, $c$ and $\varepsilon$ be
given by $a=0.1$, $b_{1}=30$, $b_{2}=0.6$, $c=1$ and $\varepsilon=0.005$
in system (\ref{2D-model-12}).
Then a canard explosion appears as the parameter $v$ varies.
See Figures \ref{fig-S1-canard-nohead-R0}, \ref{fig-S1-canard-head-R0} and \ref{fig-S1-M-relaxation}.
\begin{figure}[!htbp]
\centering
\subfigure[]{
\begin{minipage}[t]{0.32\linewidth}
\centering
\includegraphics[width=4.5cm,height=4.5cm]{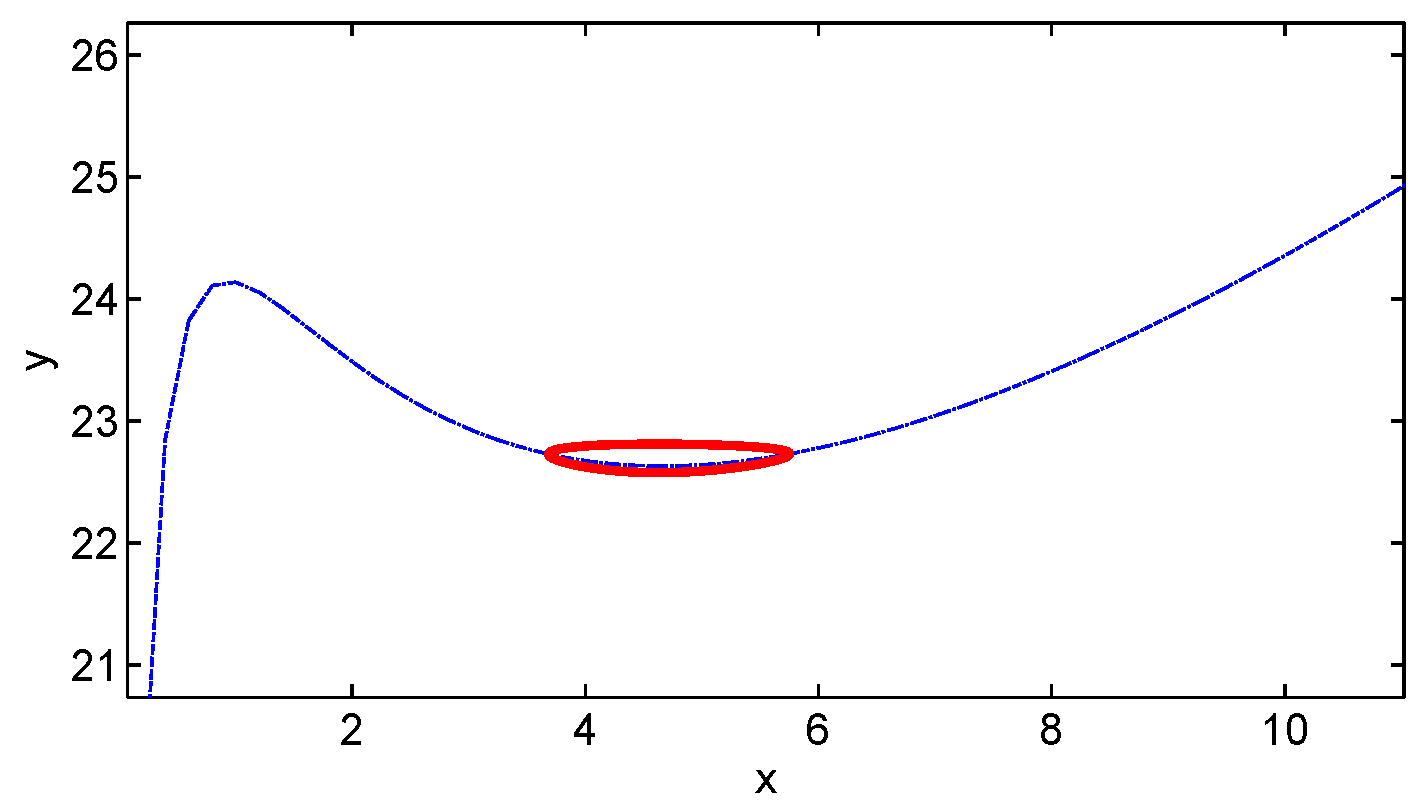}
\label{fig-S1-canard-nohead-R0}
\end{minipage}%
}%
\subfigure[]{
\begin{minipage}[t]{0.32\linewidth}
\centering
\includegraphics[width=4.5cm,height=4.5cm]{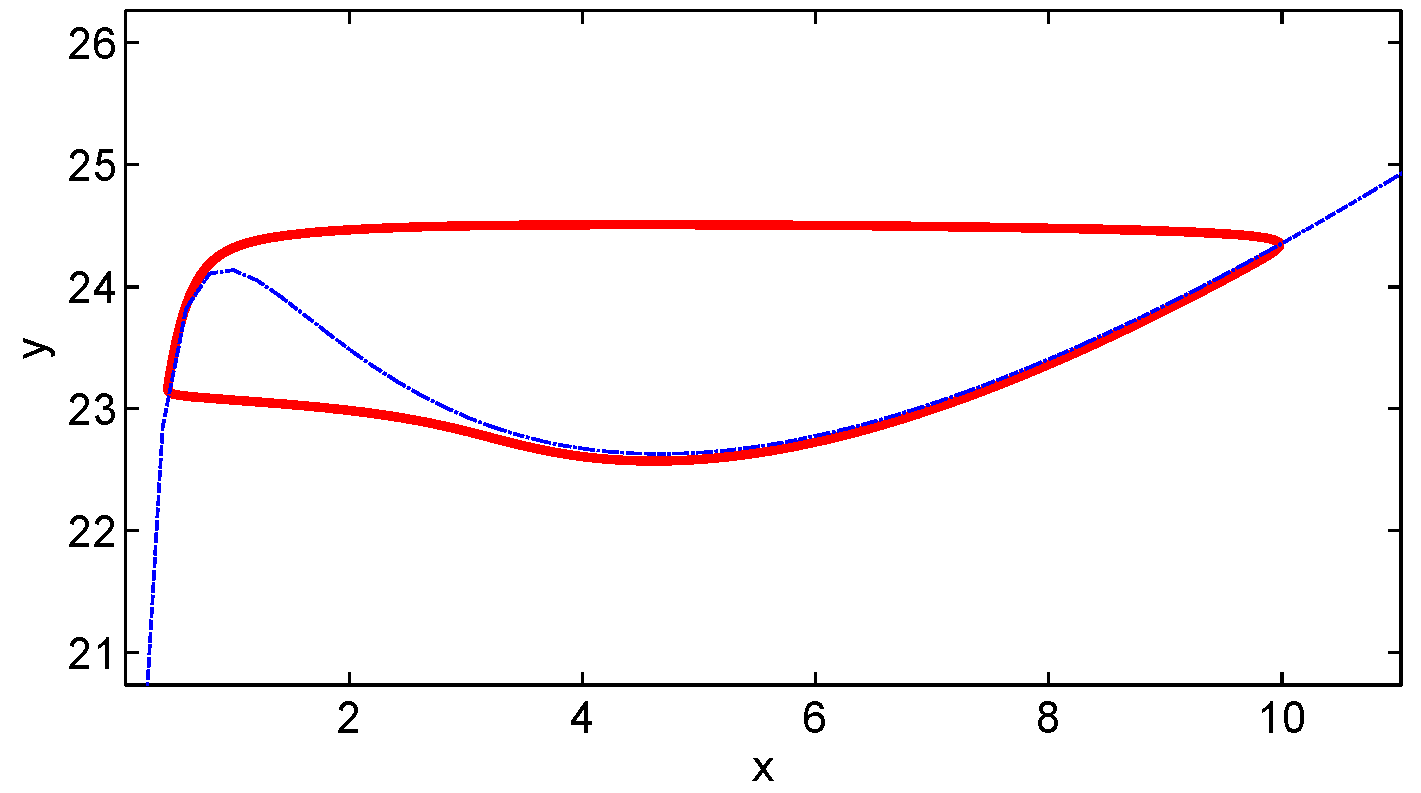}
\label{fig-S1-canard-head-R0}
\end{minipage}%
}%
\subfigure[]{
\begin{minipage}[t]{0.32\linewidth}
\centering
\includegraphics[width=4.5cm,height=4.5cm]{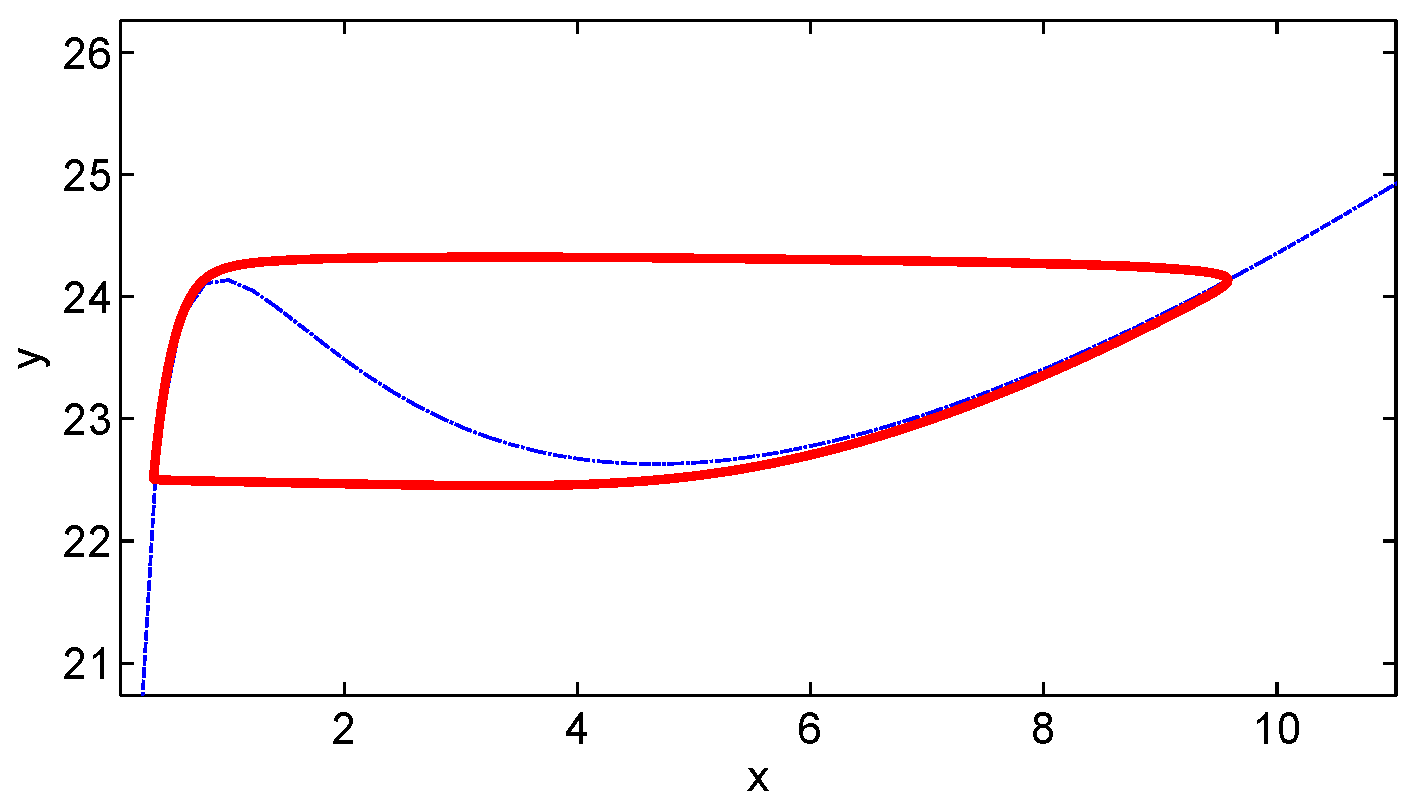}
\label{fig-S1-M-relaxation}
\end{minipage}%
}%
\centering
\caption{Canard explosion in the slow-fast system (\ref{2D-model-14}).
The dashed blue curves are the critical manifolds,
the solid red cycles indicate periodic orbits.
\ref{fig-S1-canard-nohead-R0} Canard cycle without head arises when $v=103$.\
\ref{fig-S1-canard-head-R0} Canard cycle with head arises when $v=98$.\
\ref{fig-S1-M-relaxation}  Relaxation oscillation arises when $v=53$.}
         \label{fig-S1}
\end{figure}
\end{example}

\section{Concluding remarks}
\setcounter{equation}{0}
\setcounter{lemma}{0}
\setcounter{theorem}{0}
\setcounter{remark}{0}

We have investigated the dynamics of the THTN model,
which is a circadian oscillator model based on
the dimerization and proteolysis of PER and TIM proteins in Drosophila.
After giving a classification of all possible distributions of the equilibria,
we obtain the existence of a bounded attractor in the first quadrant,
the nonexistence of periodic solutions in the high degradation rate case,
and the global dynamics in the low degradation rate case.
These results are helpful for understanding  the effects of
the rates of mRNA degradation and synthesis on the periodic oscillations in the THTN model.

More concretely, Theorem \ref{prop-dynamics-1} shows that the circadian oscillation disappears
when the rate $k_{m}$ of mRNA degradation is sufficiently high.
As a result, the oscillatory behavior requires  the rate $k_{m}$ of mRNA degradation to be bounded.
As stated in  {\bf (vi)} of  Theorem \ref{thm-one}, under some parameter conditions
there exists the configuration of a big limit cycle enclosing a small one in the THTN model
and a numerical example is presented in Figure \ref{fig-S1-coexistence}.
This interesting phenomenon suggests that depending on
the different biological environments,
the circadian oscillator exhibits different periodic behaviors.
Relaxation oscillations and canard cycles are widely found in many circadian oscillator models
(see, for instance, \cite{Forger-17,Keener-Sneyd-98}).
Theorems \ref{thm-one}, \ref{thm-two} and \ref{thm-three} show that
these oscillations could also appear in the THTN model
and Figure \ref{fig-S1}  gives several concrete examples.
For example, as shown in Figure \ref{fig-S1-M-relaxation},
when the concentration of mRNA is high,
the total amount of PER protein increases to a high level in a short time.
After that the concentration of mRNA  decreases until it reaches a low level,
and as a consequence the total amount of PER protein  quickly decreases.
Then the concentration of mRNA increases to a high level again.
This process leads to the occurrence of  a relaxation oscillation.
Theorems \ref{thm-one}, \ref{thm-two} and \ref{thm-three} also give the nonexistence of periodic solutions,
and the existence of several complex oscillations including canard explosion
and periodic solutions bifurcating from homoclinic  orbits and  heteroclinic orbits as the parameter $v$ varies.
These results suggest that the periods and the amplitudes of the circadian oscillations could be
affected by the ratio of the rate of mRNA degradation to the rate of mRNA synthesis.

It is also possible  to understand the dynamics of the the THTN model with the general rate $k_{m}$.
In fact, by some changes the THTN model can be transformed in
a Li\'enard-like equation
\begin{eqnarray*}
\begin{aligned}
\frac{d x}{d t} &= y-\left((\varepsilon+1)(x^{2}+2x)+\frac{b_{2}x^{2}+2(b_{1}+b_{2})x}{x^{2}+2x+a}\right),
\\
\frac{d y}{d t} &= 2\varepsilon(x+1)\left(\frac{v}{x^{4}+c}-\frac{b_{2}x^{2}+2(b_{1}+b_{2})x}{x^{2}+2x+a}-x^{2}-2x\right),
\end{aligned}
\end{eqnarray*}
where the parameters are defined as in system (\ref{2D-model-12}).
Then the results on Li\'enard equations (see, for instance, \cite{Dumortieretal06,ZZF-etal})
can be applied to obtain the global dynamics of the THTN model in the general case.
The Li\'enard-like structure for the THTN model  could be  helpful to
investigate the effects of the model parameters on the periods of circadian oscillations.

\section*{Appendix A: Proof of Lemma \ref{distribution-inter}}

Before proving Lemma \ref{distribution-inter},
we give the next auxiliary lemma.
\renewcommand\thelemma{A}
\begin{lemma}\label{lm-clm-01}
There exist positive parameters $c$ and $v$ such that
the graph of  $\psi_{2}$ passes
a pair of  points $(\omega_{1},y_{1})$ and $(\omega_{2}, y_{2})$ with $\omega_{1}<\omega_{2}$ in $\mathbb{R}^{2}_{+}$
if and only if the following properties hold:
\renewcommand\theequation{A.1}
\begin{eqnarray}
\label{ineq-02}
\frac{(\omega_{1}-\phi(\omega_{1}))^{2}}{(\omega_{2}-\phi(\omega_{2}))^{2}}<\frac{y_{2}}{y_{1}}<1.
\end{eqnarray}
\end{lemma}
{\bf Proof.}
If the graph of  $\psi_{2}$ passes
 points $(\omega_{1},y_{1})$ and $(\omega_{2}, y_{2})$ with $\omega_{1}<\omega_{2}$,
then $0<y_{2}<y_{1}$ and
\renewcommand\theequation{A.2}
\begin{eqnarray}
v=y_{i}(c+(\omega_{i}-\phi(\omega_{i}))^{2}), \ \ \ \ i=1,2.
\label{eq-002}
\end{eqnarray}
Clearly, the above equations have a unique solution $(c,v)$ in the form
\begin{eqnarray*}
c=\frac{y_{2}(\omega_{2}-\phi(\omega_{2}))^{2}-y_{1}(\omega_{1}-\phi(\omega_{1}))^{2}}{y_{1}-y_{2}},\ \ \
v=\frac{y_{1}y_{2}((\omega_{2}-\phi(\omega_{2}))^{2}-(\omega_{1}-\phi(\omega_{1}))^{2})}{y_{1}-y_{2}}.
\end{eqnarray*}
Since $c>0$ and $v>0$, then (\ref{ineq-02}) holds.
Thus, the sufficiency is proved.

If two points $(\omega_{i},y_{i})$ satisfy  $\omega_{1}<\omega_{2}$  and (\ref{ineq-02}),
then  these equations in (\ref{eq-002}) have a unique solution $(c,v)$ with $c>0$ and $v>0$.
Thus, the necessity is proved.
This finishes the proof.
\hfill$\Box$

Now we prove Lemma \ref{distribution-inter} by the above lemma.\\
{\bf Proof of Lemma \ref{distribution-inter}.}
By the monotonicity of the functions $\psi_{1}$ and $\psi_{2}$,
we obtain that all possible combinations of intersection point sequences are as follows:
$L$, $M$, $R$, $LM$, $MM$, $MR$, $LMM$, $LMR$, $MMM$ and $MMR$.
To complete the proof,
it is only necessary to prove that all types shown in this lemma can be realized.
Let the parameters $b_1=b\widetilde{b}_{1}$ and $b_{2}=b\widetilde{b}_{2}$,
and the function $\varphi$ be defined by $\varphi(x)=(\widetilde{b}_{1}\phi(x)+\widetilde{b}_{2}x)/(a+x)$ for $x\geq 0$.
Then $\psi_{1}(x)=b\varphi(x)+x$.
By Lemma \ref{lm-psi-1-prpty} there exist some $a^{*}$, $b^{*}$ and $\widetilde{b}_{i}^{*}$ such that
for some $x^{*}>0$,
\renewcommand\theequation{A.3}
\begin{eqnarray}
\frac{d\psi_{1}}{dx}(x^{*})=b^{*}\frac{d\varphi}{dx}(x^{*})+1=0, \ \ \
\frac{d^{2}\psi_{1}}{dx^{2}}(x^{*})=b^{*}\frac{d^{2}\varphi}{dx^{2}}(x^{*})=0,
\ \ \ \frac{d^{3}\psi_{1}}{dx^{3}}(x^{*})=b^{*}\frac{d^{3}\varphi}{dx^{3}}(x^{*})>0.
\label{eqs-1}
\end{eqnarray}
By the second equation, we observe that $x^{*}$ is independent of $b$
and only depends on the constants $a$ and $\widetilde{b}_{i}$.
Taking $c=c^{*}:=(6u^{5}_{*}+5u^{4}_{*})/(2u_{*}+3)$, $u_{*}:=\sqrt{1+x^{*}}-1$
and $v=v^{*}:=(c^{*}+(x^{*}-\phi(x^{*}))^{2})\psi_{1}(x^{*})$,
 by Lemma \ref{lm-psi-0-2-prpty} we have
 \renewcommand\theequation{A.4}
\begin{eqnarray}
\frac{d^{2}\psi_{2}}{dx^{2}}(x^{*})=0, \ \ \ \psi_{1}(x^{*})=\psi_{2}(x^{*}).
\label{eqs-2}
\end{eqnarray}
Let  the parameters $a=a^{*}$, $\widetilde{b}_{i}=\widetilde{b}_{i}^{*}$ and $c=c^{*}$ be fixed.
Consider the following equations
\renewcommand\theequation{A.5}
\begin{eqnarray}
\frac{\partial \psi_{1}}{\partial x}(x,b,v)=0, \ \ \
\psi(x,b,v)=\psi_{1}(x,b,v)-\psi_{2}(x,b,v)=0.
\label{eqs-3}
\end{eqnarray}
By (\ref{eqs-1}) and (\ref{eqs-2}) we have $(x,b,v)=(x^{*},b^{*},v^{*})$ is a solution of (\ref{eqs-3}).
Since $b^{*}>0$ and $\frac{\partial \psi_{2}}{\partial x}(x^{*},b^{*},v^{*})<0$,
then the matrix
\begin{eqnarray*}
\left(
\begin{array}{cc}
\frac{\partial^{2}\psi_{1}}{\partial x^{2}} & \frac{\partial^{2}\psi_{1}}{\partial b\partial x}\\
\frac{\partial \psi}{\partial x} & \frac{\partial \psi}{\partial b}
\end{array}
\right)_{(x^{*},b^{*},v^{*})}
=\left(
\begin{array}{cc}
0  & -\frac{1}{b^{*}}\\
\frac{\partial \psi_{2}}{\partial x}(x^{*},b^{*},v^{*}) & \varphi(x^{*},b^{*},v^{*})
\end{array}
\right)
\end{eqnarray*}
is nonsingular.
Thus by the Implicit Function Theorem, there exist two $C^{\infty}$ functions
\begin{eqnarray*}
x(v) =  x^{*}+\alpha_{1}(v-v^{*})+O((v-v^{*})^{2}),\ \ \
b(v) =  b^{*}+\alpha_{2}(v-v^{*})^{2}+O((v-v^{*})^{3})
\end{eqnarray*}
such that $\frac{\partial \psi_{1}}{\partial x}(x(v),b(v),v)=0$
and $\psi(x(v),b(v),v)=0$ for small $|v-v^{*}|$,
where the constants
\begin{eqnarray*}
\alpha_{1} = -\frac{1}{\frac{\partial \psi_{2}(x^{*},b^{*},v^{*})}{\partial x}(c^{*}+(x^{*}-\phi(x^{*}))^{2})}>0,\ \ \
\alpha_{2} = (\alpha_{1}b^{*})^{2}\frac{\partial^{3} \varphi(x^{*},b^{*},v^{*})}{\partial x^{3}}>0.
\end{eqnarray*}
For sufficiently small $|v-v^{*}|>0$ we have $b(v)>b^{*}$.
By the first equation in (\ref{eqs-1}) we obtain that
$\frac{\partial \psi_{1}}{\partial x}(x^{*},b^{*},v^{*})=-1/b^{*}<0$,
which implies that $\frac{\partial \psi_{1}}{\partial x}(x^{*},b(v),v)=1-b(v)/b^{*}<0$ for sufficiently small $|v-v^{*}|>0$.
From  Lemma \ref{lm-psi-1-prpty} it follows that the function $\frac{\partial \psi_{1}}{\partial x}(\cdot,b(v),v)$
has exactly two positive zeros $x_{m}(v)$ and $x_{M}(v)$ with $0<x_{m}(v)<x^{*}<x_{M}(v)$.
Since the constant $\alpha_{1}$ satisfies $\alpha_{1}>0$,
then $x(v)$ satisfies $x(v)=x_{m}(v)$ for $v<v^{*}$
and $x(v)=x_{M}(v)$ for $v>v^{*}$.
By continuity we obtain that  for sufficiently small $|v-v^{*}|>0$,
there is  a constant $\varrho_{2}>0$ such that $x^{*}-\varrho_{2}<x_{m}(v)<x^{*}<x_{M}(v)<x^{*}+\varrho_{2}$
and
\begin{eqnarray*}
\frac{\partial \psi_{2}}{\partial x}(x,b,v)\leq -2\varrho_{2}
<-\varrho_{2}<\frac{\partial \psi_{1}}{\partial x}(x,b,v)\leq 0 \ \  &&\mbox{ for }\  \  x_{m}(v)\leq x \leq x_{M}(v).
\end{eqnarray*}
Thus for small $v^{*}-v>0$ (resp. $v-v^{*}>0$),
equation $\psi(x,b(v),v)=0$ with respect to $x$ has exactly one positive root $x=x_{m}(v)$ (resp. $x=x_{M}(v)$).
Hence, the sequences $L^{0}$ and $R^{0}$ exist.
Under the assumption that the sequence $L^{0}$ appears,
let $(\omega_{2}, y_{2})=(x_{M},\psi_{2}(x_{M}))$ and $\omega_{1}=x_{m}$ be fixed.
By varying $y_{1}$,
we obtain $L^{1}$ by decreasing $y_{1}$ slightly from $y_{1}=\psi_{1}(x_{m})$,
and  $M$  by increasing $y_{1}$ slightly.
Similarly, we can get  the sequence $R^{1}$.
Thus, the proof for  {\bf (i)} is obtained.

Take the parameters such that the intersection point sequence $L^{0}$ appears.
Let two points $(\omega_{1}, y_{1})$ and $(\omega_{2}, y_{2})$
satisfy $(\omega_{1}, y_{1})=(x_{m},\psi_{2}(x_{m}))$ and $(\omega_{2}, y_{2})=(x_{M},\psi_{2}(x_{M}))$.
Then by Lemma \ref{distribution-inter},
\begin{eqnarray*}
(x_{m}-\phi(x_{m}))^{2}/(x_{M}-\phi(x_{M}))^{2}<\psi_{2}(x_{M})/\psi_{2}(x_{m})<1,
\end{eqnarray*}
which implies that for fixed $(\omega_{1}, y_{1})=(x_{m},\psi_{2}(x_{m}))$ and $\omega_{2}=x_{M}$,
the inequalities in (\ref{ineq-02}) hold for each $y_{2}$ with $\psi_{2}(x_{M})\leq y_{2}<\psi_{1}(x_{m})=\psi_{2}(x_{m})$.
In particular, set $y_{2}=\psi_{1}(x_{M})$.
Then by Lemma \ref{distribution-inter}
there exist some parameters $c$ and $v$ such that $\psi_{1}(x_{m})=\psi_{2}(x_{m})$ and $\psi_{1}(x_{M})=\psi_{2}(x_{M})$.
Note that $\frac{d \psi_{1}}{dx}(x_{i})=0>\frac{d \psi_{2}}{dx}(x_{i})$, $i=m, M$,
and $\psi$ has at most three positive zeros,
then there exists exactly one point $x_{3}\in (x_{m},x_{M})$ such that $\psi_{1}(x_{3})=\psi_{2}(x_{3})$.
Thus the sequence $L^{0}MR^{0}$ appears
and $\psi_{1}$ transversally intersects with  $\psi_{2}$ at three different points.
Varying $y_{2}$ slightly,
we get the sequences $L^{0}MM$ for $y_{2}-\psi_{1}(x_{M})<0$ and $L^{0}MR^{1}$ for $y_{2}-\psi_{1}(x_{M})>0$.
By decreasing $y_{2}$ again, the sequence $L^{0}M$ can be obtained.
Hence, the sequences $L^{0}M$, $L^{0}MM$, $L^{0}MR^{0}$ and $L^{0}MR^{1}$ exist for suitable parameters.
Similarly,
we can obtain the sequences $MR^{0}$, $MMR^{0}$ and $L^{1}MR^{0}$ starting from $R^{0}$,
the sequences $L^{1}M$, $L^{1}MM$ and $L^{1}MR^{1}$  from $L^{1}$,
the sequences $MM$, $MMM$ and $MMR^{1}$ from  $M$,
and the sequence $MR^{1}$ from $R^{1}$.
Thus, we give the proof for {\bf (ii)} and {\bf (iii)}.
Therefore, the proof is now complete.
\hfill$\Box$

\section*{Appendix B: Proof of (\ref{system-CM})}

We write system \eqref{SN-normal-form} as the form
\renewcommand\theequation{B.1}
\begin{eqnarray}
\label{eq-app-1}
\begin{aligned}
\frac{dx}{dt}&= \bar{X}_{2}(x+y,v),
\\
\frac{dy}{dt}&= \mu_{1}y+\bar{Y}_{2}(x+y,v),\\
\frac{dv}{dt}&= 0,
\end{aligned}
\end{eqnarray}
where $\bar{X}_{2}$ and $\bar{Y}_{2}$ denote $X_{2}$ and $Y_{2}$ with  $v^{0}$ replaced by $v+v^{0}$, respectively.
Since system \eqref{eq-app-1} has two zero eigenvalues and one nonzero eigenvalue $\mu_{1}$ at the origin,
then by the Center Manifold Theory \cite[Section 1.3]{Carr-81},
system \eqref{eq-app-1} has a $C^{3}$ center manifold $y=\tilde{y}(x,v)$ for sufficiently small $|x|$ and $|v|$.
By a direct computation,
the restriction of \eqref{eq-app-1} to the center manifold has the expansion
\renewcommand\theequation{B.2}
\begin{eqnarray}
\label{eq-app-2}
\begin{aligned}
\frac{dx}{dt}=&\, \frac{\varepsilon}{D_{1}\psi_{1}(x_{0},\lambda^{0})+\varepsilon}
        \times \left(\frac{1}{c^{0}+(x_{0}-\phi(x_{0}))^{2}}v\right.\\
        & +(D_{13}\psi_{2}(x_{0},\lambda^{0},v^{0}))(x+\tilde{y}(x,v))v\\
       &+\frac{1}{2}\left( D_{11}\psi_{2}(x_{0},\lambda^{0},v^{0})-D_{11}\psi_{1}(x_{0},\lambda^{0})\right)(x+\tilde{y}(x,v))^{2}\\
       &\left.
       +\frac{1}{2}\left( D_{33}\psi_{2}(x_{0},\lambda^{0},v^{0})\right)v^{2} \right)+O(|(x,v)|^{3}),\\
\frac{dv}{dt}=&\, 0.
\end{aligned}
\end{eqnarray}
Note that $D_{33}\psi_{2}(x_{0},\lambda^{0},v^{0})=0$ and
$\tilde{y}(x,v)=O(|(x,v)|^{2})$ for  sufficiently small $|x|$ and $|v|$.
Then  we can write (\ref{eq-app-2}) as the form (\ref{system-CM}).
This finishes the proof.

%


\end{document}